\documentclass[preprint,12pt]{elsarticle}

\usepackage[T1]{fontenc}
\usepackage{lmodern}
\usepackage{amsmath,amssymb,amsthm,mathtools}
\usepackage{array,longtable}
\usepackage{microtype}
\usepackage{enumitem}
\usepackage[hidelinks]{hyperref}
\usepackage{bookmark}
\newtheorem{theorem}{Theorem}[section]
\newtheorem{lemma}[theorem]{Lemma}
\newtheorem{proposition}[theorem]{Proposition}
\newtheorem{corollary}[theorem]{Corollary}
\theoremstyle{definition}
\newtheorem{definition}[theorem]{Definition}
\newtheorem{example}[theorem]{Example}
\theoremstyle{remark}
\newtheorem{remark}[theorem]{Remark}

\newcommand{\HR}{\mathsf H_R}
\newcommand{\KR}{\mathcal K_R}
\newcommand{\Std}{\mathsf{Std}}
\newcommand{\seq}[1]{\Rightarrow #1}

\newcommand{\DM}{\operatorname{DM}}
\DeclareMathOperator{\Un}{U}
\setlist{nosep}
\hypersetup{pdftitle={Strong Standard Completeness of a Hypersequent Calculus for Cyclic Involutive Residuated Lattices},pdfsubject={Manuscript revision 1.1.5; companion Lean development v1.1.0}}

\journal{Annals of Pure and Applied Logic}

\begin{document}

\begin{frontmatter}

\title{Strong Standard Completeness of a Hypersequent Calculus for
Cyclic Involutive Residuated Lattices}
\author[aff1]{SanMin Wang\corref{cor1}}
\ead{wangsanmin@hotmail.com}
\cortext[cor1]{Corresponding author.}
\affiliation[aff1]{organization={Faculty of Science, Zhejiang Sci-Tech University},
  city={Hangzhou}, postcode={310018}, country={P.R. China}}

\begin{abstract}
We present a hypersequent calculus $\HR$ for cyclic involutive
residuated lattices and prove strong standard completeness for arbitrary
theories in a fixed countable propositional language. Standard models
are chains on $[0,1]$; neither commutativity nor equality of the monoid
unit and the dualizing constant is assumed. The chain semantics is
the class $\KR$ defined by unit cyclicity and all ordered finite-history
conditions. Every countable member of $\KR$ embeds into a standard
chain preserving the entire signature. The construction combines
central-support reconstruction, gap insertion and relative contraction.
A relative comparison construction supplies the required chain
countermodels. We also prove admissibility of the unrestricted density
rule for theoremhood and effective elimination from finite proofs
without additional assumption leaves. Elimination uses terminating
enumeration of density-free proofs; ordinary Cut is retained, and no
complexity bound is claimed. Further constructions give bounded
standard embeddings for nontrivial countable chains of involutive
uninorm logic and pseudo-uninorm logic with weak commutativity.
Using ordinary chain completeness for arbitrary theories, these yield
strong standard completeness for both logics. A companion Lean~4
development formalizes the main results, with its scope and external
inputs stated explicitly.
\end{abstract}


\end{frontmatter}

\section{Introduction}
\label{sl:introduction}
We introduce an explicit hypersequent calculus $\HR$ in the
setting of cyclic involutive residuated lattices and prove its
strong standard completeness. The standard models are chains
whose underlying ordered set is the real unit interval $[0,1]$;
neither commutativity nor identification of the monoid unit with
the dualizing constant is assumed. The completeness theorem
concerns consequence from arbitrary theories in a fixed countable
propositional language, with all premises interpreted under the
same valuation.

The chain semantics of $\HR$ is the class $\KR$ specified in
Section~\ref{sl:logic:calculus}: cyclic involutive residuated chains
satisfying unit cyclicity and every ordered finite-history
condition. These requirements hold in every densely ordered
cyclic involutive residuated chain
(Lemma~\ref{sl:logic:dense-histories}), hence in every standard
model.

The main algebraic result is a full standard representation
theorem for every countable member of $\KR$
(Theorem~\ref{sl:rep:standard}). It preserves multiplication,
both residuals, involution and both designated constants
simultaneously. A comparison construction relative to arbitrary
theories supplies countable countermodels in $\KR$
(Theorem~\ref{sl:logic:theory-chain-completeness}); a single full
embedding transports each countermodel and all its premises to
a standard chain. This yields strong standard completeness for
$\HR$ (Theorem~\ref{sl:logic:strong-standard}). We also obtain
effective density elimination for finite proofs without
additional assumption leaves, retaining ordinary Cut
(Theorem~\ref{sl:logic:den-elimination}).

Sections~\ref{sl:iul:section} and~\ref{sl:psul:section} apply these
constructions to IUL and $\mathsf{HpsUL}^{*}$. The applications
require, respectively, commutative specializations and an
involutive extension satisfying all histories, together with
preservation of the original endpoints. The resulting bounded
standard embeddings yield strong standard completeness for
these logics using their ordinary chain-completeness theorems
for arbitrary theories in the respective countable languages.

The connection between density elimination and standard completeness is
established in the theory of substructural fuzzy logics.  Metcalfe and
Montagna developed this approach for uninorm logic and several of its
extensions~\cite{MetcalfeMontagna2007}.  Ciabattoni and Metcalfe subsequently
gave general sufficient conditions for density elimination by
substitution~\cite{CiabattoniMetcalfe2008}.  These proof-theoretic criteria
require hypotheses on the underlying calculus; they do not establish
admissibility for an arbitrary extension by additional structural rules.
Our use of ordinary Cut and the entire family of ordered history rules
therefore requires a separate argument.

The algebraic counterpart is densification: embedding a chain into a
dense chain while preserving its operations and the defining conditions
of the class.  Metcalfe and Tsinakis show that every dense residuated
chain is $e$-cyclic, and relate densification of semilinear varieties to
generation by their dense chains as quasivarieties
\cite[Theorem~1 and Lemma~1]{MetcalfeTsinakis2017}.
Baldi and Terui obtain densification through residuated frames for
specified integral varieties and commutative varieties
\cite{BaldiTerui2016}; Galatos and Hor\v{c}\'ik give an explicit
one-gap construction for commutative residuated chains
\cite[Theorems~3.5--3.6]{GalatosHorcik2022}.
Preservation of an existing residuated reduct does not by itself make
the old designated element dualizing on the extension.  This is a
separate requirement in the present problem.

The closest proof-theoretic precedents include Wang's published
density-elimination arguments for involutive uninorm logic and for the
noncommutative system $\mathbf{GpsUL}^{*}$
\cite{WangIUL2019,WangPseudo2019}.
The latter system adds weak commutativity at the monoid unit to
pseudo-uninorm logic; its language and rules do not impose a cyclic
dualizing element.  Its density-elimination and standard-completeness
statements are Theorems~4 and~5 of the published article.
The former paper concerns commutative involutive logic.  Jenei's later
discussion~\cite[Introduction]{JeneiFixedPoint2023} records reservations
about the verification of its proof.  We neither use that proof as an
input nor resolve those reservations here.

For commutative involutive chains, Jenei's group representation
\cite{JeneiGroups2022,JeneiGroupsCorrection2022,JeneiGroupsCorrection2023}
and subsequent densification theorem give a particularly close
algebraic comparison.  Every countable odd commutative involutive
chain embeds into a countable dense chain of the same class
\cite[Corollary~17]{JeneiDensification2025}; the bounded version preserves
the endpoints, and yields strong standard completeness of
$\mathbf{IUL}^{fp}$~\cite[Corollary~19 and Theorem~20]{JeneiDensification2025}.
Here odd means that the monoid unit and the dualizing constant coincide.

\begin{table}[tbp]
\centering
\small
\renewcommand{\arraystretch}{1.15}
\begin{tabular}{@{}>{\raggedright\arraybackslash}p{.17\textwidth}>{\raggedright\arraybackslash}p{.34\textwidth}>{\raggedright\arraybackslash}p{.41\textwidth}@{}}
\hline
Aspect & Jenei's densification result & Present representation theorem \\
\hline
Input class & Odd commutative involutive chains & $\KR$, allowing $d\ne e$ and noncommutative products \\
Construction & Representation by layer groups & Central support, complete fibres, and relative return collapse \\
Embedding result & Dense extension in the same class; bounded endpoints retained & Full standard embedding; bounded endpoints treated separately in the IUL and pseudo-uninorm transfers \\
Stated logical result & Strong standard completeness of $\mathbf{IUL}^{fp}$ & Strong standard completeness of $\HR$, IUL and $\mathsf{HpsUL}^{*}$;
the latter two use ordinary theory-chain completeness. Den admissibility for $\HR$ \\
\hline
\end{tabular}
\caption{Scope of the comparison with \cite[Corollaries~17, 19 and Theorem~20]{JeneiDensification2025}.}
\label{sl:tab:jenei-comparison}
\end{table}

Table~\ref{sl:tab:jenei-comparison} distinguishes the input
classes, constructions and preservation requirements.
Our construction retains the full signature for the explicitly
axiomatized class $\KR$, allowing noncommutative products and $d\ne e$;
its input conditions are essential to the representation statement.

The proof separates the chain by a closed central support term.
The difficult lower branch is reduced to complete group fibres.
Canonical completion attains their internal self-support maxima;
pullback and a countable elementary substructure retain the exact
original group. A lexicographic enlargement then makes that group
dense. These preparation steps make it possible to insert points in
all remaining gaps while preserving associativity, both actual
residuals and all histories. The upper branch requires protected
support intervals and a separate reconstruction. Detailed auxiliary
constructions and the full chain-completeness proof are supplied in
Appendices~\ref{sl:logic:chain-supplement}--\ref{sl:upper:section},
with explicit cross-references.

\paragraph{Machine-checked formalization}
A companion Lean~4 development provides machine-checked proofs
of the main representation and completeness results established
here, together with an executable procedure for Den elimination.
The source code and reproducibility materials are archived
on Zenodo~\cite{WangLeanArchive}.
Section~\ref{sl:formal:scope} specifies the verified scope,
external inputs, and computational guarantees.

\section{The calculus and its chain semantics}
\label{sl:logic:calculus}

A cyclic involutive residuated lattice is a lattice-ordered monoid
$(L,\wedge,\vee,\cdot,e)$ with residuals satisfying
\[
 xy\le z\quad\Longleftrightarrow\quad y\le x\backslash z
 \quad\Longleftrightarrow\quad x\le z/y,
\]
and a distinguished element $d$ such that
$x'=x\backslash d=d/x$ and $x''=x$.
We use the derived operation $x+y=(x'y')'$; in particular,
\begin{equation}\label{sl:logic:terms}
 e'=d,\qquad
 x\backslash y=(y'x)'=y+x',\qquad
 y/x=(xy')'=x'+y.
\end{equation}
A \emph{chain} has a total lattice order.  A \emph{standard chain}
has underlying ordered set $[0,1]$.  Neither $e$ nor $d$ is required
to be an endpoint, and multiplication need not be commutative.
The class of all standard cyclic involutive residuated chains is
denoted by $\Std$.

Formulas are built from countably many variables, their formal duals,
$e,d$, and $\cdot,+,\wedge,\vee$.  Involution is extended as a
syntactic operation by $A''=A$, $e'=d$, $d'=e$, and
\[
 (AB)'=A'+B',\quad (A+B)'=A'B',\quad
 (A\wedge B)'=A'\vee B',\quad
 (A\vee B)'=A'\wedge B'.
\]
Thus the displayed order is not reversed under involution.
Residual formulas abbreviate the terms in \eqref{sl:logic:terms}.
All substitution instances use this convention.

A word $\Gamma$ is a finite ordered list of formulas, possibly empty.
A sequent is $\seq\Gamma$, and a hypersequent is a finite multiset
of sequents, written with external bars.  External permutation is
implicit; internal exchange is not.  Capital Latin letters $G,H,J$
denote hypersequents, and Greek capitals denote words.  A missing
external background denotes the empty multiset.  All words in the
following rules may be empty.

The initial sequents and logical rules are
\[
 \overline{\seq{A,A'}}\;\mathrm{Id}
 \qquad
 \overline{\seq e}\;\mathrm{Unit},
\]
\begin{equation}\label{sl:logic:logical-rules}
\begin{gathered}
 \frac{G\mid\seq{A,\Gamma}\qquad H\mid\seq{B,\Delta}}
      {G\mid H\mid\seq{AB,\Gamma,\Delta}}\;(\cdot),
 \\[6pt]
 \frac{G\mid\seq{\Gamma,A,B,\Delta}}
      {G\mid\seq{\Gamma,A+B,\Delta}}\;(+)
 \qquad
 \frac{G\mid\seq{\Gamma,\Delta}}
      {G\mid\seq{\Gamma,d,\Delta}}\;(d),
 \\[6pt]
 \frac{G\mid\seq{\Gamma,A,\Delta}\qquad
       G\mid\seq{\Gamma,B,\Delta}}
      {G\mid\seq{\Gamma,A\wedge B,\Delta}}\;(\wedge),
 \\[6pt]
 \frac{G\mid\seq{\Gamma,A,\Delta}}
      {G\mid\seq{\Gamma,A\vee B,\Delta}}\;(\vee_1)
 \qquad
 \frac{G\mid\seq{\Gamma,B,\Delta}}
      {G\mid\seq{\Gamma,A\vee B,\Delta}}\;(\vee_2).
\end{gathered}
\end{equation}
The context of $(\wedge)$ is shared literally.  The structural rules
and ordinary Cut are
\begin{equation}\label{sl:logic:structural-rules}
\begin{gathered}
 \frac{G\mid\seq{\Gamma,\Delta}}
      {G\mid\seq{\Delta,\Gamma}}\;\mathrm{Cyc}
 \qquad
 \frac{G}{G\mid S}\;\mathrm{EW}
 \qquad
 \frac{G\mid S\mid S}{G\mid S}\;\mathrm{EC},
 \\[6pt]
 \frac{G\mid\seq{\Gamma_1,\Gamma_2}\qquad
       H\mid\seq{\Delta_1,\Delta_2}}
      {G\mid H\mid\seq{\Gamma_1,\Delta_2}
               \mid\seq{\Delta_1,\Gamma_2}}\;\mathrm{COM},
 \\[6pt]
 \frac{G\mid\seq{\Gamma,\Delta,e}}
      {G\mid\seq{\Delta,\Gamma,e}}\;\mathrm{WCM}
 \qquad
 \frac{G\mid\seq{\Gamma,A}\qquad H\mid\seq{A',\Delta}}
      {G\mid H\mid\seq{\Gamma,\Delta}}\;\mathrm{Cut}.
\end{gathered}
\end{equation}
Here $S$ is an arbitrary sequent.  Cyclically rotating the two active
premises and the two outputs of COM gives its equivalent presentation
that exchanges a middle block in each premise.  The displayed WCM
retains its literal final formula $e$; it is not unrestricted exchange.

For every finite word $w$ over $\{S,T\}$ (with empty word $\epsilon$)
define closed formulas by
\begin{equation}\label{sl:logic:history-terms}
\begin{gathered}
 (c_\epsilon,h_\epsilon)=(e,e),\\
 (c_{wS},h_{wS})=(dc_w,e+h_w),\qquad
 (c_{wT},h_{wT})=(e+c_w,dh_w).
\end{gathered}
\end{equation}
The remaining rules, one for every such $w$, are exactly
\begin{equation}\label{sl:logic:history-rule}
 \frac{G\mid\seq{c_w,\Sigma,\Pi}\qquad
       H\mid\seq{h_w,\Omega,\Lambda}\qquad
       J\mid\seq{\Gamma,\Delta}}
      {G\mid H\mid J\mid
       \seq{\Pi,\Sigma,\Gamma,\Lambda,\Omega,\Delta}}\;R_w.
\end{equation}
We denote this calculus by $\HR$.  Its proofs are finite trees.
There is no initial sequent $\seq d$, no internal weakening or
contraction, and no density rule in $\HR$.  Ordinary Cut is retained.

Under a valuation in a chain, put
$\Gamma^\sharp=A_1'\cdots A_n'$ for
$\Gamma=A_1,\ldots,A_n$, with $\varnothing^\sharp=e$.
The sequent $\seq\Gamma$ means $\Gamma^\sharp\le d$,
and an external bar means disjunction.  The empty hypersequent
is therefore false.  Let $\KR$ consist of the cyclic involutive
residuated chains satisfying
\begin{equation}\label{sl:logic:unit-cyclicity}
 ab\le e\quad\Longrightarrow\quad ba\le e
\end{equation}
and, for every $w\in\{S,T\}^*$,
\begin{equation}\label{sl:logic:history-condition}
 ab\le c_w,\quad uv\le h_w
 \quad\Longrightarrow\quad bazvu\le z
 \qquad(a,b,u,v,z\in L).
\end{equation}
The history terms in these conditions are evaluated by
\eqref{sl:logic:history-terms}.  In particular, membership in $\KR$
requires every finite history, without a fixed bound on its length.

\begin{lemma}[Rotation on a dense chain]\label{sl:logic:rotation}
Let $L$ be a densely ordered chain.  Suppose that $f,g:L\to L$
are monotone, $g$ has a right adjoint $g^*$, and $f(g(t))\le t$
for all $t$.  Then $g(f(t))\le t$ for all $t$.
\end{lemma}
\begin{proof}
If $g(f(x))>x$, adjunction implies $g^*(x)<f(x)$.
Choose $g^*(x)<t<f(x)$.  Then $g(t)>x$, whence
$f(g(t))\ge f(x)>t$, a contradiction.
\end{proof}

\begin{lemma}\label{sl:logic:dense-histories}
Every densely ordered cyclic involutive residuated chain belongs to
$\KR$.
\end{lemma}
\begin{proof}
If $ab\le e$, apply Lemma~\ref{sl:logic:rotation} to
$f(t)=at$ and $g(t)=bt$.  Evaluating $g(f(t))\le t$ at $e$
gives $ba\le e$.  Condition \eqref{sl:logic:unit-cyclicity}, with
$a,b$ interchanged, says precisely $x\backslash e=e/x$.
By \eqref{sl:logic:terms},
\[
 (dx)'=x\backslash e=e/x=(xd)',
\]
so $d$ is central.  Consequently $e+x=(dx')'=d\backslash x$.

Induction on $w$ now gives
\begin{equation}\label{sl:logic:unrotated-bound}
                         c_wzh_w\le z\qquad(z\in L).
\end{equation}
The empty history is immediate.  For the two induction steps use
centrality of $d$ and $d(d\backslash x)\le x$ to obtain
\begin{align*}
 (dc_w)z(d\backslash h_w)
   &=c_wz\,d(d\backslash h_w)\le c_wzh_w\le z,\\
 (d\backslash c_w)z(dh_w)
   &=d(d\backslash c_w)zh_w\le c_wzh_w\le z.
\end{align*}
No centrality of $c_w$ or $h_w$ is assumed.

If $ab\le c_w$ and $uv\le h_w$, then
$abtuv\le c_wth_w\le t$.  Apply
Lemma~\ref{sl:logic:rotation} to $f(t)=atv$ and $g(t)=btu$;
the latter has right adjoint $g^*(s)=(b\backslash s)/u$.
The resulting inequality is $bazvu\le z$, as required.
\end{proof}

\begin{theorem}[Soundness]\label{sl:logic:soundness}
Every rule of $\HR$ preserves validity in $\KR$.
In particular $\HR\vdash H$ implies $\Std\models H$.
\end{theorem}
\begin{proof}
An active component $\seq{A,\Gamma}$ is true exactly when
$\Gamma^\sharp\le A$.  Thus the two premises of $(\cdot)$
give $\Gamma^\sharp\Delta^\sharp\le AB$.  The $(+)$ rule
replaces $A'B'$ by the equal value $(A+B)'$, and the $(d)$ rule
inserts the unit $d'=e$.  Multiplication preserves binary joins in
each argument, so $(A\wedge B)'=A'\vee B'$ verifies $(\wedge)$.
The join rules follow from monotonicity and order reversal of
involution.  Id and Unit hold because $AA'\le d$ and $e'=d$.
Cyc follows from $xy\le d$ if and only if $yx\le d$.

For Cut, write $\gamma=\Gamma^\sharp$ and
$\delta=\Delta^\sharp$.  The active premises give
$\gamma A'\le d$, hence $\gamma\le A$, and $A\delta\le d$;
therefore $\gamma\delta\le d$.
In WCM the active premise is $\gamma\delta d\le d$,
equivalently $\gamma\delta\le e$, since $d/d=e$.
Condition \eqref{sl:logic:unit-cyclicity} gives the conclusion.
For COM, put $\gamma_i=\Gamma_i^\sharp$ and
$\delta_i=\Delta_i^\sharp$.  If $\gamma_2\le\delta_2$,
then $\delta_1\gamma_2\le\delta_1\delta_2\le d$;
if $\delta_2\le\gamma_2$, then
$\gamma_1\delta_2\le\gamma_1\gamma_2\le d$.
One conclusion component is therefore true.

For $R_w$, put
$(a,b,u,v,r,s)=(\Sigma^\sharp,\Pi^\sharp,
\Omega^\sharp,\Lambda^\sharp,\Gamma^\sharp,\Delta^\sharp)$.
The three active premises give $ab\le c_w$, $uv\le h_w$,
and $rs\le d$.  Condition \eqref{sl:logic:history-condition}
gives $barvus\le rs\le d$, exactly the ordered conclusion.

If a premise background is true, its inherited copy makes the
conclusion true.  For shared $(\wedge)$ this observation uses its
one common background.  EW and EC preserve external disjunction.
This proves soundness in $\KR$; standard soundness follows from
Lemma~\ref{sl:logic:dense-histories}.
\end{proof}

\subsection{Countable chain completeness}
\label{sl:logic:chain-completeness-section}

\begin{definition}[Derivability and satisfaction with assumptions]
\label{sl:logic:theory-consequence}
For any set $T$ of formulas, write $T\vdash_{\HR}H$ when the
finite hypersequent $H$ has a finite $\HR$ proof allowing
$\seq A$, for $A\in T$, as additional leaves.
These are specified sequents, not substitution schemata;
repeated uses are allowed. No internal structural rule is added.
For a valuation $v$ in a chain $C$, write $v\models T$ if
$e_C\le v(A)$ for every $A\in T$.
This agrees with the sequent interpretation, since
$v(A)'\le d_C$ if and only if $e_C\le v(A)$.
For a class $\mathcal C$ of chains, write $T\models_{\mathcal C}H$
if every valuation in every member of $\mathcal C$ that satisfies
$T$ also satisfies $H$.
For a formula conclusion we abbreviate
$T\vdash_{\HR}\seq A$ by $T\vdash_{\HR}A$.
All theories in this paper belong to the fixed countable
propositional language; they need not be finite or decidable.
\end{definition}

\begin{theorem}[Countable chain completeness with assumptions]
\label{sl:logic:theory-chain-completeness}
Let $T$ be any set of formulas and let $H$ be a finite hypersequent.
If $T\nvdash_{\HR}H$, there are a countable chain $C\in\KR$
and a valuation $v$ in $C$ such that $v\models T$ and
$v\not\models H$.
\end{theorem}
\begin{proof}[Proof outline]
Start with the specified leaves $\mathcal T_0=\{\seq A:A\in T\}$.
The continuation lemma lifts a proof from an additional leaf $S$
over any background $J$ for which $J\mid S$ is proved. COM on two
identities proves $\seq{A',B}\mid\seq{B',A}$. If adjoining either
component made $H$ derivable, two continuation substitutions and EC
would already prove $H$. Thus one orientation of each formula pair
can be adjoined while preserving underivability.

Enumerating the pairs yields a total comparison preorder whose
formula quotient is a countable cyclic involutive residuated chain.
WCM and the actual $R_w$ instances verify unit cyclicity and every
ordered history. Each assumption leaf $\seq A$ yields
$\seq{d,A}=\seq{e',A}$, so $[e]\le[A]$ in the quotient.
A true component of $H$ would unfold by Cut and the logical rules
to a proof of that component, hence of $H$ by EW. The canonical
valuation therefore satisfies $T$ and falsifies $H$.
Appendix~\ref{sl:logic:chain-supplement} supplies the
full continuation, quotient, and truth arguments.
\end{proof}

\begin{theorem}[Countable chain completeness]\label{sl:logic:chain-completeness}
If a finite hypersequent $H$ is not derivable in $\HR$, then it is
false under a valuation in a countable member of $\KR$.
\end{theorem}
\begin{proof}
Take $T=\varnothing$ in Theorem~\ref{sl:logic:theory-chain-completeness}.
\end{proof}

\section{Completion and preservation of the history conditions}
\label{sl:completion-section}

An embedding is called \emph{full} if it preserves the lattice operations,
multiplication, both residuals, involution, and the two designated
constants. For chains, an order embedding preserving multiplication,
involution and the constants is full, by the residual formulas
\eqref{sl:logic:terms}.

We use the involutive Dedekind--MacNeille completion theorem
\cite[Corollary~4.6]{GilFerezEtAl2017Preprint}; the published treatment is
\cite{GilFerezEtAl2020}. It supplies a complete cyclic involutive
residuated lattice and a full embedding of the original algebra.
We write $\DM(C)$ for the resulting completion of $C$, equipped
with these extended algebraic operations. In the chain case $\DM(C)$
is again a chain. The additional history conditions require the
following verification.

\begin{theorem}[Full completion]\label{sl:completion}
If $C\in\KR$, then $\DM(C)\in\KR$, and the natural embedding is full.
\end{theorem}
\begin{proof}
Identify $C$ with its image in $\widehat C=\DM(C)$. Every $x\in\widehat C$
is the join of its old lower approximants:
\[
 x=\bigvee\{x_0\in C:x_0\le x\}.
\]
Multiplication in a complete residuated lattice preserves arbitrary joins
separately in each coordinate. If $ab\le e$, then old $a_0\le a$ and
$b_0\le b$ satisfy $a_0b_0\le e$, hence $b_0a_0\le e$. Taking the two
independent joins gives $ba\le e$.

Fix a finite history $w$, and suppose $ab\le c_w$ and $uv\le h_w$ in
$\widehat C$. Fullness preserves the values of the closed guards.
For old lower approximants of $a,b,u,v,z$, the original history gives
\[
 b_0a_0z_0v_0u_0\le z_0\le z.
\]
Taking the five independent joins yields $bazvu\le z$. If an approximating
set is empty, its join is the bottom element, and the corresponding
product is bottom by complete join preservation. Thus the same argument
also covers that case. Each finite history is preserved, with no bound
on its length.
\end{proof}

\begin{corollary}[Standard completion]\label{sl:standard-completion}
Every nontrivial countable densely ordered member of $\KR$ has a full
embedding into a member of $\Std\cap\KR$.
\end{corollary}
\begin{proof}
Its MacNeille completion is a complete densely ordered chain with
endpoints and a countable order-dense subset. It is therefore
order-isomorphic to $[0,1]$. Transport the operations along such an
isomorphism and apply Theorem~\ref{sl:completion}.
\end{proof}

Completing an arbitrary chain does not remove its adjacent pairs.
Consequently the essential task below is to construct full dense
extensions before applying this corollary.

\section{Central support and reconstruction}
\label{sl:struct:section}

We write $e$ for the monoid identity and $d$ for the cyclic
dualizer. Thus $x'=x\backslash d=d/x$ is an antitone involution,
and
\[
 x\backslash y=(y'x)',\qquad y/x=(xy')'.
\]
An embedding is \emph{full} if it preserves this entire signature,
including the two constants. Preservation of order, multiplication,
prime and the constants already implies fullness. We use the
history constants and the class $\KR$ defined in the preceding
section. Only unit cyclicity is needed until histories are explicitly
invoked in the reconstruction argument.

\begin{lemma}[Centrality and the history bounds]
\label{sl:struct:centrality}
Unit cyclicity $xy\le e\Rightarrow yx\le e$ is equivalent to
centrality of $d$. Under these equivalent conditions,
\[
 e+x=d\backslash x,
 \qquad c_wzh_w\le z\quad(w\in\{S,T\}^*,\ z\in L).
\]
\end{lemma}
\begin{proof}
Unit cyclicity in both directions says $x\backslash e=e/x$.
The residual formulas identify these elements with $(dx)'$ and
$(xd)'$, respectively, proving the equivalence. Centrality gives
$e+x=(dx')'=(x'd)'=d\backslash x$. If $c_wzh_w\le z$, then
\[
 (dc_w)z(d\backslash h_w)
   =c_wz\,d(d\backslash h_w)\le c_wzh_w,
\]
and
\[
 (d\backslash c_w)z(dh_w)
   =d(d\backslash c_w)zh_w\le c_wzh_w.
\]
These are the two induction steps; the empty history is equality.
\end{proof}

\begin{theorem}[The closed central support]
\label{sl:struct:support}
In a unit-cyclic cyclic involutive residuated chain put
\[
 t=d^3\backslash d^2,\qquad p=d^2t^2,\qquad
 q=pd=d^2t,\qquad Q=pt=d^2t^3.
\]
Then $p,q,Q$ are central and
\[
\begin{gathered}
 d^3t=d^2,\qquad p^2=p,\qquad I:=pL=Lp=d^2L,\\
 qQ=Qq=p,\qquad Qd=p,\qquad q^2=d^2.
\end{gathered}
\]
The set $I$ is a monoid with identity $p$. Multiplication by $d$
on $I$ is an order automorphism, with inverse multiplication by $t$.
In particular $q,Q$ are mutually inverse central units of $I$.
\end{theorem}
\begin{proof}
Set $h=(d^2)'=d\backslash e$. The residual formulas and centrality
give $t=(d^3h)'$. For
\[
 b=d^3\backslash(d^3h)=(d^3t)'
\]
we have $h\le b$ and $d^3b=d^3h$. Also $dh\le e$.
If $db\le e$, then $b\le h$, and hence $b=h$. If $e\le db$,
then
\[
 d^2\le d^3b=d^3h=d^2(dh)\le d^2.
\]
Consequently $t=(d^3h)'=h$ and $b=(d^3t)'=h$ in this case too.
Involution now yields $d^3t=d^2$.

Multiplying this equality by $dt$ gives $d^4t^2=d^2$; thus
$pd^2=d^2$ and $p^2=d^4t^4=p$. Both inclusions in $pL=d^2L$
follow from these equations. On this set, $dt(d^2a)=d^2a$, so
the monotone maps given by $d$ and $t$ are mutual inverses.
In particular multiplication by $d^2$ is injective there.
Since $d$ is central by Lemma~\ref{sl:struct:centrality}, so is $d^2$.
For any $a\in L$ we have
\[
 pa=d^2(t^2a),\qquad ap=a(d^2t^2)=d^2(at^2),
\]
so both $pa,ap$ belong to $d^2L$. Moreover,
$d^2p=pd^2=d^2$, and centrality of $d^2$ gives
\[
 d^2(pa)=d^2a=d^2(ap).
\]
Thus $pa=ap$, proving centrality of $p$, and then of $q=pd$.
The displayed equations for $q,Q,d$ follow by multiplication.
Finally $Qa,aQ\in I$ and
$q(Qa)=pa=ap=q(aQ)$. Since $q$ acts invertibly on $I$,
$Qa=aQ$. All cancellation in this proof occurs inside $I$.
\end{proof}

The following formulation also applies to the alternative interior
support arising in the upper branch.

\begin{proposition}[Interior projection and exact fibres]
\label{sl:struct:interior}
Let $r\le e$ be a central idempotent of a cyclic involutive
residuated chain $L$. The map
\[
 \pi:L\longrightarrow J=rL,\qquad \pi(x)=rx
\]
is a full surjective homomorphism when $J$ is given identity $r$,
dualizer $q=rd$, and involution $a^*=ra'$. In particular,
\[
 (rx)\backslash_J(ry)=r(x\backslash y),\qquad
 (ry)/_J(rx)=r(y/x),\qquad (rx)^*=rx'.
\]
Every fibre has the exact form
\begin{equation}\label{sl:struct:fibres}
 F_a:=\pi^{-1}(a)=[a,(a^*)']_L.
\end{equation}
If $q$ is a central unit of $J$, every nonunit of $J$ has a
singleton fibre.
\end{proposition}
\begin{proof}
For $z\in J$ we have $(rx)z=xz$. If this is at most $ry$, then
$xz\le y$, whence $z=rz\le r(x\backslash y)$. Conversely,
$z\le r(x\backslash y)$ gives $xz\le ry$. This proves the left
residual formula; reversing the factor positions proves the right
one. Products, minima, maxima and the constants are preserved.
The residual formulas give $(rx)^*=rx'$, so $*$ is involutive and
cyclic. If $rx=a$, then $a\le x$ and $a^*=rx'\le x'$, proving
the two bounds in \eqref{sl:struct:fibres}. Conversely its two
endpoints project to $a$, and monotonicity proves the equality.

Here is the singleton assertion in a form useful below. For any
surjective full homomorphism $\varphi:B\to A$ of such chains,
$aa'<d_A$ forces $\varphi^{-1}(a)$ to be a singleton: if $u,w$
belong to that fibre, then
$\varphi(uw')=aa'<\varphi(d_B)$, so $uw'<d_B$ and $u\le w$;
interchange $u,w$. A nontrivial fibre over $a\in J$ has a
nontrivial complementary fibre. Therefore $aa^*=a^*a=q$.
If $q^{-1}$ is the local inverse of $q$, then $q^{-1}a^*$ is
a two-sided inverse of $a$ in $J$. This proves the assertion
without any assumption concerning one-sided inverses.
\end{proof}

\begin{lemma}[The complete group core]
\label{sl:struct:core}
Under the hypotheses of Proposition~\ref{sl:struct:interior},
suppose $q$ is a central unit. Write $U(J)$ for the group of units
of the monoid $J$ with its local identity $r$, and put
\[
 G=U(J),\qquad K=\pi^{-1}(G).
\]
Then $K$ is a full subalgebra of $L$,
\[
 rK=G,\qquad L=J\cup K,\qquad J\cap K=G,
\]
and, for $a\in J$ and $x\in L$,
\begin{equation}\label{sl:struct:floor-actions}
 ax=a\pi(x),\qquad xa=\pi(x)a.
\end{equation}
If $L\in\KR$, then $K\in\KR$.
\end{lemma}
\begin{proof}
The group $G$ is a full subalgebra of $J$ with its local constants:
$g^*=qg^{-1}$, $g\backslash_J h=g^{-1}h$ and
$h/_Jg=hg^{-1}$. Minima and maxima select an argument. Thus its
inverse image is full and contains $e,d$. Singleton nonunit fibres
give the union and intersection assertions. Since $a=ra=ar$,
centrality gives \eqref{sl:struct:floor-actions}. Universal
quasi-identities, including all histories, pass to subalgebras.
\end{proof}

\begin{definition}[Complete group core]
\label{sl:struct:complete-core-definition}
Under the hypotheses of Lemma~\ref{sl:struct:core}, the
\emph{complete group core of $L$ at $r$} is
\[
 K_r(L)=\pi^{-1}\bigl(U(rL)\bigr),\qquad \pi(x)=rx.
\]
Thus every fibre above a unit of the local algebra $rL$ is included,
with its original operations.  The word \emph{complete} refers to this
inclusion of entire fibres; it does not mean Dedekind completeness.
Lemma~\ref{sl:struct:core} says that this inverse image is a full
subalgebra, rather than a choice of one representative from each fibre.
\end{definition}

For the lower branch $p\le e$, the support and interior results apply
with $r=p$ and local dualizer $q=pd$, whose inverse is $Q$.
Consequently the core at the closed support is exactly
$K_p(L)=\{x\in L:px\in U(pL)\}$.

\begin{lemma}[Central unit dualizers and histories]
\label{sl:struct:unit-histories}
Every cyclic involutive residuated chain with central invertible
dualizer $q$ belongs to $\KR$. With local identity $r$ its history
constants are
\[
 c_w=q^{k(w)},\qquad h_w=q^{-k(w)},\qquad
 k(w)=\#S(w)-\#T(w),
\]
where $\#S(w)$ and $\#T(w)$ count the occurrences of $S$ and $T$
in $w$, and powers use identity $r$.
\end{lemma}
\begin{proof}
The operation $x^\circ=q^{-1}x^*$ makes $r$ a cyclic dualizer,
because
\[
 xy\le r\ \Longleftrightarrow\ x(qy)\le q
 \ \Longleftrightarrow\ y\le q^{-1}x^*.
\]
For any central unit $a$, $(ax)^*=a^{-1}x^*$, directly from
residuation. Hence $(x^\circ)^\circ=x$. The unit threshold is
cyclic, and multiplication by $q^{-n}$ gives
\[
 xy\le q^n\ \Longleftrightarrow\ yx\le q^n.
\]
On the central units $a+b=q^{-1}ab$, so induction gives the
history constants displayed above. Rotating both cap premises
then gives $bazvu\le q^{k(w)}zq^{-k(w)}=z$.
\end{proof}

\begin{theorem}[Simultaneous core and support gluing]
\label{sl:struct:gluing}
Let $L\in\KR$ have a central idempotent $r\le e$ such that
$J=rL$ has central unit dualizer $q=rd$. Put $G=U(J)$ and
$K=\{x\in L:rx\in G\}$. Suppose:
\begin{enumerate}
\item $K\hookrightarrow\widehat K\in\KR$ is full, $r$ remains
central in $\widehat K$, and $r\widehat K=G$ as the old embedded
group;
\item $J\hookrightarrow T$ is full for the local constants $r,q$,
and $q$ remains a central unit in $T$.
\end{enumerate}
Then $L$ has a full extension $B\in\KR$ with underlying set
\[
 B=\widehat K\mathbin{\dot\cup}(T\setminus G).
\]
The order replaces every $g\in G$ in $T$ by its entire
$\widehat K$ fibre and leaves all other $T$ points single.
Multiplication by $r$ projects $B$ onto the floor copy of $T$.
Both the given core extension and all its internal inserted points
are retained. Countability is preserved.
\end{theorem}
\begin{proof}
Apply Proposition~\ref{sl:struct:interior} in $\widehat K$.
Its projection $\widehat\pi(k)=rk$ is full onto $G$. Define
$\Pi:B\to T$ by $\widehat\pi$ on $\widehat K$ and the identity
outside. Order first by $\Pi$, then by the old fibre order.
The floor inclusion $i:T\to B$ sends $g\in G$ to its least
fibre point $g$ and is the identity elsewhere; hence
\begin{equation}\label{sl:struct:floor-inequality}
 i\Pi(x)\le x.
\end{equation}
Keep multiplication on $\widehat K$ and otherwise put
\begin{equation}\label{sl:struct:glued-product}
 xy=i(\Pi(x)\Pi(y)).
\end{equation}
For $g\in G$ and $k\in\widehat K$, centrality of $r$ gives
$gk=g\widehat\pi(k)$ and $kg=\widehat\pi(k)g$. Consequently
the uniform floor identities are
\begin{equation}\label{sl:struct:glued-floor}
 i(a)x=i(a\Pi(x)),\qquad xi(a)=i(\Pi(x)a).
\end{equation}
Every product containing an outside factor is therefore the
floor of the product of all its projections, independently of
bracketing. This proves associativity; products wholly in
$\widehat K$ use its associativity. It also proves that the old
$e$ is the identity, since $\Pi(e)=r$.

For a fixed floor factor, monotonicity follows in $T$. A remaining
factor $k$ belongs to $\widehat K$ and has unit projection $g$.
If the two varying arguments have distinct projections, left or
right multiplication by $g$ preserves their strict order. If their
projections agree, the arguments are equal outside $G$, and inside
$G$ monotonicity is that of $\widehat K$. Thus multiplication is
monotone across every boundary as well as within each piece.

Keep prime on $\widehat K$, and use the local prime $*$ of $T$
on $T\setminus G$. Both $G$ and its complement are $*$-invariant:
$g^*=qg^{-1}$ belongs to $G$ and $*$ is an involution. The full
projection $\widehat\pi$ shows that these definitions give an
antitone involution on $B$. For a product with an outside factor,
\[
 xy\le d\ \Longleftrightarrow\ \Pi(x)\Pi(y)\le q
 \ \Longleftrightarrow\ \Pi(y)\le\Pi(x)^*.
\]
The first equivalence uses that a floor $i(a)$ is at most $d$
exactly when $a\le q$. In the last comparison, if both bases
are outside $G$ they refer to singleton fibres; if just one is
outside $G$ equality is impossible. In either case it is exactly
$y\le x'$. For two inputs in $\widehat K$ this duality already
holds. Associativity and cyclic duality now give both residuals
by the involutive formulas. They also show that $\Pi$ is full
for the local constants of $T$.

Unit cyclicity holds inside $\widehat K$. Otherwise both ordered
products are floors, and the assertion follows from unit cyclicity
in $T$ and $r\le e$. All closed history constants project to
their corresponding constants in $T$. If all five factors in a
history instance lie in $\widehat K$, use that instance there.
Otherwise its conclusion product is a floor, and the projected
instance, valid by Lemma~\ref{sl:struct:unit-histories}, gives
\[
 bazvu=i(\Pi(b)\Pi(a)\Pi(z)\Pi(v)\Pi(u))
       \le i\Pi(z)\le z.
\]
Thus $B\in\KR$ for every finite history.

The old chain is $K\cup(J\setminus G)$. Its mixed products agree
with \eqref{sl:struct:glued-product} by
\eqref{sl:struct:floor-actions}, and its outside prime agrees
with the local prime because those fibres are singletons.
Its inclusion is therefore full. Finally
$rx=i\Pi(x)$ by \eqref{sl:struct:glued-floor}, proving all
remaining assertions.
\end{proof}

Taking $T=J$ gives fixed-floor core gluing. Taking
$\widehat K=K$ gives lifting of every full support extension
that retains the central unit dualizer. Neither specialization
requires that nonunits form an ideal: two outside elements may
multiply to an old group element, and the resulting point is
its floor, as prescribed by \eqref{sl:struct:glued-product}.

\begin{theorem}[Densification at a central unit dualizer]
\label{sl:struct:unit-dense}
Every countable cyclic involutive residuated chain with central
invertible dualizer has a countable full dense extension without
endpoints, retaining its dualizer as a central unit. Both chains
belong to $\KR$.
\end{theorem}
\begin{proof}
Normalize its prime to $x^\circ=d^{-1}x'$ as in
Lemma~\ref{sl:struct:unit-histories}. The resulting chain is odd.
Corollary~\ref{sl:struct:odd-dense}, proved by the explicit
construction in Appendix~\ref{sl:struct:odd-appendix}, gives a
countable dense odd extension without endpoints, in which the
old $d$ remains central with its inverse. Restore prime by
$x'=dx^\circ$. The same normalization calculation in reverse
proves cyclic duality and involutivity. Multiplication, order
and both residuals are unchanged, and the restored prime agrees
with the original prime on the old chain. The embedding is full.
Lemma~\ref{sl:struct:unit-histories} gives all histories in the
source and target.
\end{proof}

\begin{example}[A noncommutative, nonodd central-unit chain]
\label{sl:struct:example-heisenberg}
Let $H=\mathbb Z^3$ with lexicographic order, in which the first
coordinate is compared first, and multiplication
\[
 (a,b,c)(r,s,t)=(a+r,b+s,c+t+as).
\]
This is a countable group with
\[
 e=(0,0,0),\qquad (a,b,c)^{-1}=(-a,-b,ab-c).
\]
Associativity follows by comparing the third coordinates of
both bracketings: for a third factor $(u,v,w)$, both are
$c+t+w+as+av+rv$. The order is invariant under both left and
right translations. Indeed, translation preserves the first
unequal coordinate; if the first two coordinates coincide,
the same quantity is added to the third coordinate of both
elements. Thus $H$ is a totally ordered residuated group, with
\[
 x\backslash y=x^{-1}y,\qquad y/x=yx^{-1}.
\]

Choose the central element $d=(0,0,1)$ as dualizer. Then
\[
 (a,b,c)'=d(a,b,c)^{-1}=(-a,-b,1+ab-c)
\]
is an antitone involution, and
$x\backslash d=x^{-1}d=dx^{-1}=d/x$. Hence this is a cyclic
involutive residuated chain, and it is nonodd because $d\ne e$.
Its support terms satisfy
\[
 t=d^{-1}=(0,0,-1),\qquad p=e,\qquad q=d,\qquad Q=d^{-1}.
\]
It belongs to $\KR$ by Lemma~\ref{sl:struct:unit-histories}, with
$c_w=d^{k(w)}$, $h_w=d^{-k(w)}$ and $k(w)=\#S(w)-\#T(w)$.

For $x=(1,0,0)$ and $y=(0,1,0)$, direct calculation gives
\[
 xy=(1,1,1),\qquad yx=(1,1,0),
\]
and
\[
 x\backslash y=(-1,1,-1),\qquad y/x=(-1,1,0).
\]
Thus the left and right multiplication actions, and the two
residual operations, differ even in the central-unit case
$p=e$. Both local identities of every group element are still
$x/x=x\backslash x=e$.
\end{example}

\begin{corollary}[Reconstruction from dense group fibres]
\label{sl:struct:reconstruction}
Let $L\in\KR$ be countable, and let $r\le e$ be a central
idempotent such that $q=rd$ is a central unit of $J=rL$.
Put $G=U(J)$ and $K=\{x:rx\in G\}$. Suppose $K$ has a
countable full extension $\widehat K\in\KR$ in which $r$
is central, $r\widehat K=G$, and every fibre of multiplication
by $r$ is densely ordered, with singleton fibres allowed.
Then $L$ has a countable full dense extension without endpoints,
and hence a full standard embedding.
\end{corollary}
\begin{proof}
Theorem~\ref{sl:struct:unit-dense} gives a countable dense local
extension $J\hookrightarrow T$ without endpoints, retaining $q$
as a central unit. Apply Theorem~\ref{sl:struct:gluing} to
$\widehat K$ and $T$. For two resulting points with different
projections, choose an intermediate projection in $T$ and take
its floor. Two points with the same projection lie in one of the
given dense fibres; a singleton fibre contains no strict pair.
There are no endpoints because $T$ has none. This proves density
and the remaining countability assertion. The full standard
completion is supplied by Corollary~\ref{sl:standard-completion}.
\end{proof}

\begin{corollary}[Extracting a fixed-floor extension from a standard core]
\label{sl:struct:core-standard-reconstruction}
In the setting of Corollary~\ref{sl:struct:reconstruction},
suppose the complete core $K$ has a full embedding into a standard
chain $S\in\KR$ in which the old $r$ is central. Then $L$ has
a full standard embedding. In particular this applies to the lower
branch $r=p\le e$, where centrality of the closed term $p$ in
$S$ is automatic.
\end{corollary}
\begin{proof}
Identify $K$ with its image. The old inverse equations make
$q=rd$ a central unit of $rS$. The group $G$ is a full local
subalgebra of $rS$, because its residuals are group division and
its prime is $g^*=qg^{-1}$. By
Proposition~\ref{sl:struct:interior},
\[
 E=\{s\in S:rs\in G\}
\]
is a full subalgebra containing $K$, with $rE=G$. Its fibre over
$g$ is the entire interval $[g,(g^*)']_S$. Both endpoints belong
to the old core. Hence every strict pair in a fibre of $E$ has
an intermediate point in that fibre, by density of the standard
chain. Starting with $K$, choose one such point for every strict
same-fibre pair, take the full subalgebra generated by all these
choices inside $E$, and repeat countably many times. The union
$\widehat K$ is countable, full, has $r\widehat K=G$, and all
its fibres are densely ordered. Apply
Corollary~\ref{sl:struct:reconstruction}. For $r=p$, centrality
in $S$ follows from Theorem~\ref{sl:struct:support}, since the
embedding preserves the defining closed term.
\end{proof}

\section{Standard representation of countable lower complete cores}
\label{sl:lower:section}

\begin{definition}[Lower complete core]
\label{sl:lower:core-definition}
A member $K\in\KR$ is a \emph{lower complete core} if its closed support
$p$ satisfies $p\le e$ and $pK$ is a group with identity $p$.
Equivalently, in this case $K=K_p(K)$.
Countability is a separate hypothesis in the representation results.
\end{definition}

We prove the representation theorem for a countable lower complete core
$K$, and write $G=pK$.  Throughout, embeddings preserve the entire
original signature. We use the structural identities
\[
 p^2=p\le e,\qquad p,d\text{ central},\qquad q=pd\in Z(G),
 \qquad d^2=q^2,
 \tag{L1}\label{sl:lower:structure}
\]
and the full interior projection $\pi(x)=px$. Its fibres are convex
and the fibre of $g\in G$ has least point $g$. On the group $G$, the
identity is $p$, the dualizer is $q$, and $g^*=qg^{-1}$; consequently
\[
\begin{gathered}
 \pi(xy)=\pi(x)\pi(y),\qquad \pi(x^-)=\pi(x)^{-1},\\
 \pi(x')=q\pi(x)^{-1},\qquad x^-=x\backslash e=e/x.
\end{gathered}
 \tag{L2}\label{sl:lower:projection}
\]
These facts follow from the structural results
\ref{sl:struct:centrality}, \ref{sl:struct:support}, and
\ref{sl:struct:interior}. We also use support-base lifting
(Theorem~\ref{sl:struct:gluing}) and full completion
(Theorem~\ref{sl:completion}). All the gap
insertion and preservation arguments needed here are proved explicitly.

We first construct one simultaneous gap insertion, assuming attained
self-supports and internal gaps. The self-support maxima determine actual
comparison bounds and hence the greatest cyclic product on the new
carrier. This comparison construction alone does not prove associativity.
The relative return collapse of Lemma~\ref{sl:lower:collapse} excludes
the active polar pairs that would be forced by a failed associator in
Theorem~\ref{sl:lower:associativity}. After proving that these
properties survive insertion, we supply the preparation that makes them
available for an arbitrary complete core and iterate the construction.

\subsection{The old actions on the split carrier}

For an $e$-cyclic cyclic involutive residuated chain $A$, let $D(A)$ be
obtained by inserting exactly one point $X_{a,b}$ in each adjacent pair
$a\prec b$ of $A$. Extend prime by $X_{a,b}'=X_{b',a'}$. Write
$r(X_{a,b})=b$ and $r(a)=a$ for $a\in A$. The construction in this
subsection requires neither endpoints nor a support group.

An old labelled context $f=(u,v)$ means $f(x)=uxv$. An admissible
return is a context $R=(s,t)$ such that
\[
 su,us,vt,tv\le e.
 \tag{L3}\label{sl:lower:return}
\]
Then $Rf\le\mathrm{id}$ and $fR\le\mathrm{id}$. We call $f$ active
at $(a,b)$ if $Rf(b)=b$ for an admissible return. Equivalently,
\[
                 u^-u b v v^-=b.
 \tag{L4}\label{sl:lower:activity}
\]
Indeed any admissible $s,t$ satisfy $s\le u^-$ and $t\le v^-$,
while $(u^-,v^-)$ is itself admissible.

\begin{lemma}[The prescribed old actions]
\label{sl:lower:actions}
If $f$ is active at $(a,b)$, then $f(b)$ has an immediate predecessor
$c$. The formulas
\[
 f(X_{a,b})=
 \begin{cases}X_{c,f(b)},&f\text{ is active at }(a,b),\\
 f(b),&f\text{ is inactive at }(a,b),\end{cases}
 \tag{L5}\label{sl:lower:action-definition}
\]
together with the old values define monotone actions on $D(A)$.
They are monotone also in each coefficient, preserve composition of
labelled contexts, and satisfy $r(fX)=f(rX)$. Each action has an
actual right adjoint, whose value on an old target is the old right
adjoint. In particular, with $L_uX=uX$ and $R_uX=Xu$,
\[
\begin{gathered}
 L_uL_v=L_{uv},\quad R_vR_u=R_{uv},\quad L_uR_v=R_vL_u,\\
 L_u^*(Y)=(Y'u)',\qquad R_u^*(Y)=(uY')'.
\end{gathered}
 \tag{L6}\label{sl:lower:action-laws}
\]
The old unit fixes every point.
\end{lemma}
The proof is given in Appendix~\ref{sl:lower:appendix-actions}.

For a left action, activity is exactly $u^-ub=b$; for a right
action it is $buu^-=b$. A two-sided admissible return for a
one-sided action can always be replaced by a one-sided return:
if $subt=b$ and $t\le e$, then $b\le sub\le b$.

\subsection{Attainment of coefficient comparisons}

For a gap $J$ define its positive self-support set by
\[
 S_J=\{u\in A:u\ge e,\ uX_J=X_J\}.
 \tag{L7}\label{sl:lower:self-support}
\]
We write (PC) for the condition that each $S_J$ has a greatest
member $s_J$. In old operations this condition says that
$\{u\ge e:ub=b,\ u^-b=b\}$ has a maximum for every gap
$J=(a,b)$. Indeed $ub=b$ changes the activity equation
$u^-ub=b$ into $u^-b=b$.

\begin{lemma}[Transport and comparison maxima]
\label{sl:lower:comparison}
Assume (PC). Every cut $\{u\in A:uX\le Y\}$, for $X,Y\in D(A)$,
has a greatest member $m(X,Y)$. For old $x,y$ and gaps
$J=(a,b)$, $K=(c,h)$, it is given by
\[
\begin{array}{c|c}
\text{case}&m(X,Y)\\\hline
Y=y\in A&y/r(X)\\
X=x\in A,\ Y=X_K&c/x\\
X=X_J,\ Y=X_K,\ \text{no coefficient sends }X_J\text{ to }X_K&c/b\\
X=X_J,\ Y=X_K,\ uX_J=X_K&(u^-)\backslash s_J=s_K/u^-.
\end{array}
 \tag{L8}\label{sl:lower:comparison-table}
\]
In the last case the maximum itself sends $X_J$ to $X_K$.
Moreover $s_J^2=s_J$.
\end{lemma}
The proof is given in Appendix~\ref{sl:lower:appendix-comparison}.

By order reversal, all the following minima now exist:
\[
 k(X,Y)=\min\{z\in A:z'X\le Y'\}=m(X,Y')'.
 \tag{L10}\label{sl:lower:k}
\]
Conversely, their existence implies (PC) by taking $(X,Y)=(X_J,X_J')$.
Right-coefficient cuts are attained as well, since
$Xu\le Y\iff uY'\le X'$ by the old-action adjunctions.

\subsection{The greatest cyclic multiplication}

\begin{lemma}[Cyclic selection and associativity with old factors]
\label{sl:lower:selection}
Under (PC), define on $D(A)$
\[
 P(X,Y)=\min\{k(X,Y),\ r(X)Y,\ Xr(Y)\}.
 \tag{L11}\label{sl:lower:P}
\]
This operation is monotone, unital, and retains every prescribed
product having an old factor. For $z\in A$ it satisfies
\[
 P(X,Y)\le z\iff z'X\le Y',\qquad r(P(X,Y))=k(X,Y).
 \tag{L12}\label{sl:lower:old-target}
\]
It satisfies all cyclic adjunctions and has actual residuals
\[
 X\backslash Y=P(Y',X)',\qquad Y/X=P(X,Y')'.
 \tag{L13}\label{sl:lower:residuals}
\]
Furthermore, for $u\in A$,
\[
 P(uX,Y)=uP(X,Y),\quad P(X,Yu)=P(X,Y)u,
 \quad P(Xu,Y)=P(X,uY),
 \tag{L14}\label{sl:lower:old-assoc}
\]
and the two bracketings of every triple have the same old upper
anchor:
\[
 r(P(P(X,Y),Z))=r(P(X,P(Y,Z))).
 \tag{L15}\label{sl:lower:anchor-assoc}
\]
\end{lemma}
The proof is given in Appendix~\ref{sl:lower:appendix-selection}.

\subsection{Active coefficients and relative return collapse}

\begin{lemma}[Closure of active coefficients]
\label{sl:lower:active-closure}
The coefficients having a nonempty active left domain are exactly
those having a nonempty active right domain. Their common set
$\Gamma(A)$ is closed under multiplication and under $u\mapsto u^-$;
if $A$ has a gap, it contains $e$.
\end{lemma}
\begin{proof}
Suppose $uX=Y$ is active, with $X,Y$ new and anchors $b,h$.
The actual right adjoint $L_u^*(Y)=(Y'u)'$ is at least $X$ but
strictly below $b$, because $ub=h>Y$. The only point in $[X,b)$
is $X$, so $Y'u=X'$, an active right action. The opposite
argument gives the converse. For the negative return put $v=u^-$; the activity equations give
$vh=b$ and $uvh=h$, so $v$ is active on the image gap and
returns its midpoint to $X$. This proves closure under $u\mapsto u^-$.
The unit is active at every gap.

Let $u,v\in\Gamma(A)$, and put $\alpha=uu^-\le e$,
$\beta=v^-v\le e$. An image upper endpoint for $u$ is fixed
by $\alpha$, and an active input upper endpoint for $v$ is
fixed by $\beta$. If $\alpha\le\beta$, an image upper endpoint
$h$ for $u$ satisfies $h=\alpha h\le\beta h\le h$, so $v$
is active at that image gap and $vu$ is active by composition.
If $\beta\le\alpha$, take an active input upper endpoint $c$
for $v$. Then $c=\beta c\le\alpha c\le c$, so $uu^-c=c$.
The inequalities $uu^-\le e$ and $u^-u\le e$, the latter by
unit cyclicity, make $u$ an admissible return for the left action
of $u^-$. Thus $u^-$ is active at this gap by
\eqref{sl:lower:return} and the definition of activity.
Its image gap has upper endpoint $b_0=u^-c$, and
\[
 ub_0=c,\qquad u^-u b_0=u^-c=b_0.
\]
Consequently $u$ is active at the image gap by
\eqref{sl:lower:activity}, and sends it back to the active input
gap of $v$. Composition gives an active action of $vu$.
This argument uses $u$ as an admissible return and requires no
identity $u^{--}=u$. Here $u^-=u\backslash e=e/u$, whereas
prime denotes the given cyclic involution.
This proves multiplicative closure.
\end{proof}

Say that a lower complete core has internal gaps if
\[
              a\prec b\quad\Longrightarrow\quad pb<b.
 \tag{L19}\label{sl:lower:NF}
\]
This is equivalent to $pa=pb$ for every gap. Indeed if $pa<pb$,
then $a<pb\le b$ and adjacency forces $pb=b$. Conversely,
$pa=pb=b$ would contradict $pa\le a<b$. A densely ordered
floor group guarantees internal gaps: if $pb=b$ were the upper
endpoint of a gap, a floor strictly between $pa$ and $b$ would
lie strictly between $a$ and $b$. The same property holds for a
trivial floor group, whose sole point is least in the whole chain.

\begin{lemma}[Relative return collapse]
\label{sl:lower:collapse}
In a lower complete core with internal gaps, every $z\le d$ with
$pz=q$ has empty left and right active domains.
\end{lemma}
\begin{proof}
Put $t=z^-$. By \eqref{sl:lower:projection}, $pt=q^{-1}$.
Centrality of $d,p$ and $d^2=q^2$ give
\[
 (tz)^2\le(td)^2=t^2d^2=t^2q^2=(pt)^2q^2=p.
 \tag{L20}\label{sl:lower:relative-collapse}
\]
The insertion of $p$ is valid because $pq^2=q^2$; no interchange
of $t$ and $z$ is used. If $z$ were active on the left at a
gap with upper endpoint $b$, then $tzb=b$, so
$b=(tz)^2b\le pb\le b$, contrary to \eqref{sl:lower:NF}.
For a right action, use $(zt)^2\le(dt)^2=p$ and $bzt=b$.
\end{proof}

\begin{theorem}[Associativity of the internal-gap insertion]
\label{sl:lower:associativity}
If a lower complete core $A\in\KR$ has internal gaps and satisfies
(PC), the multiplication $P$ of \eqref{sl:lower:P} is associative.
\end{theorem}
\begin{proof}
We first record what a strict correction means. If
$k=k(X,Y)<r(X)r(Y)$, both $X,Y$ are new, since with an old
factor their prescribed product has anchor $r(X)r(Y)$.
Write their gaps as $(a,b),(c,h)$. The inequality $k'X\le Y'$
cannot hold by an inactive action: that would give the old value
$k'b\le h'$ and hence $bh\le k$, a contradiction. If it is
active with upper anchor below $c'$, that anchor is at most $h'$
and gives the same contradiction. It therefore has anchor $c'$
and equals the unique midpoint $Y'$. In particular $k'$ is active.

For arbitrary $X,Y,Z$ put $b=r(X)$, $h=r(Y)$, $t=r(Z)$,
$a=k(X,Y)$, $c=k(Y,Z)$, and
\[
 L=P(P(X,Y),Z),\qquad R=P(X,P(Y,Z)),\qquad K=r(L)=r(R).
\]
The equality of anchors follows from Lemma~\ref{sl:lower:selection}.
Expanding \eqref{sl:lower:P}, using the old-factor equations and
preservation of minima by old actions, gives
\[
 \begin{split}
 L&=\min\{K,aZ,bYt,X(ht)\},\\
 R&=\min\{K,(bh)Z,bYt,Xc\}.
 \end{split}
 \tag{L21}\label{sl:lower:associator-expansion}
\]
The discarded terms are $at\ge aZ$ and $bc\ge Xc$.

If $L\ne R$, the one-midpoint carrier and their common anchor
force one to be $K$ and the other its immediate predecessor.
Suppose $L<K=R$. The second line implies $(bh)Z,bYt,Xc\ge K$;
since $c\le ht$, also $X(ht)\ge K$. Thus the first line forces
$L=aZ<K$, which is new. Hence $a$ is active. Necessarily
$a<bh$, so the preceding strict-correction argument makes $a'$
active as well. By Lemma~\ref{sl:lower:active-closure}, $aa'$
is active. But cyclic duality gives $aa'\le d$, and
\[
 \pi(aa')=\pi(a)q\pi(a)^{-1}=q.
\]
Lemma~\ref{sl:lower:collapse} excludes this. If $R<K=L$, apply
the same argument in the opposite monoid with the triple reversed;
it gives an active right polar pair, and the equality of active
coefficient sets gives the identical contradiction. Thus $L=R$.
\end{proof}

\begin{lemma}[Fullness and all finite histories]
\label{sl:lower:history-transfer}
The associative insertion in Theorem~\ref{sl:lower:associativity} belongs
to $\KR$, and $A\hookrightarrow D(A)$ is full.
\end{lemma}
The proof is given in Appendix~\ref{sl:lower:appendix-history-transfer}.
It treats an arbitrary finite word $w$, rather than reducing the history
scheme to the particular short words used in the structural lemmas.

\begin{lemma}[The insertion invariants]
\label{sl:lower:invariants}
If $A$ is a lower complete core with internal gaps and (PC), then
$D(A)$ is again such a core, with exactly the same floor group
$pD(A)=pA$, the same $p,q$, internal gaps, and (PC).
\end{lemma}
The proof is given in Appendix~\ref{sl:lower:appendix-invariants}.

\begin{example}[One internal gap and its comparison correction]
\label{sl:lower:example-three-point}
Let $A=\{0<1\}$, with $e=1$, $d=0$, $xy=\min(x,y)$,
and $0'=1$, $1'=0$. Its residuals are
\[
 x\backslash y=y/x=
 \begin{cases}1,&x\le y,\\y,&x>y.\end{cases}
\]
This is an integral commutative cyclic involutive residuated chain.
It belongs to $\KR$: unit cyclicity is immediate, and
$bazvu\le z$ holds for all five elements by integrality and
monotonicity, independently of the history premises. The closed
support terms are
\[
 t=0\backslash0=1,\qquad p=q=Q=0.
\]
Thus $G=pA=\{0\}$ is the trivial floor group and $A$ is a lower
complete core. Its sole gap $J=(0,1)$ is internal, because
$p1=0<1$.

Write $X=X_{0,1}$, so that $D(A)=\{0<X<1\}$, $X'=X$, and
$r(X)=1$. Every old coefficient has $u^-=u\backslash e=1$.
Consequently the left and right activity equations at $J$ reduce
to $u=1$. More generally, an old context $(u,v)$ is active exactly
when $u=v=1$. The prescribed actions therefore give
\[
 0X=X0=0,\qquad 1X=X1=X.
\]
Since $1$ is the only old coefficient at least $e$, one has
$S_J=\{1\}$ and $s_J=1$; in particular, (PC) holds. The defining
comparison minimum and selected product are
\[
 k(X,X)=\min\{z\in A:z'X\le X'\}=0,\qquad
 P(X,X)=\min\{0,X,X\}=0.
\]
Here $0$ qualifies in the first minimum because $0'X=1X=X=X'$.
The complete multiplication table is therefore
\[
\begin{array}{c|ccc}
 P&0&X&1\\\hline
 0&0&0&0\\
 X&0&0&X\\
 1&0&X&1
\end{array}
\]
and the residual formulas give
\[
 X\backslash0=P(1,X)'=X,\qquad
 X\backslash X=P(X,X)'=1.
\]
The table is associative: $0$ is absorbing, $1$ is the identity,
and any product with at least two occurrences of $X$ is $0$.
It is commutative and integral, so every finite history holds
by the same inequality $bazvu\le z$. The old residuals are
unchanged, and the old prime and constants are retained, giving
a full embedding $A\hookrightarrow D(A)$. This example exhibits
the role of the comparison term explicitly: $k(X,X)=0$, whereas
$r(X)r(X)=1$.
\end{example}

\subsection{Preparing an arbitrary countable lower core}

We attain (PC) on a changed carrier by canonical extension.
The precise completion facts and their source locators are stated
in Appendix~\ref{sl:lower:appendix-canonical}.

\begin{lemma}[Canonical completion and attained gap supports]
\label{sl:lower:canonical}
If $A\in\KR$ is bounded, then its canonical extension
$A^\delta\in\KR$ is a full chain extension in which every gap has
greatest positive self-support.
\end{lemma}
The proof is given in Appendix~\ref{sl:lower:appendix-canonical}.

\begin{lemma}[Attained internal supports with the original floor group]
\label{sl:lower:internal-PC}
Every countable lower complete core $K$ has a countable full
extension $C\in\KR$ with $pC=pK=G$ such that every gap whose
endpoints have equal degree has greatest positive self-support.
\end{lemma}
\begin{proof}
Apply the full completion theorem to $K$, and then
Lemma~\ref{sl:lower:canonical}, obtaining $K\subseteq B\in\KR$
with support maxima at every gap. The old closed term $p$ stays
idempotent and at most $e$. It stays central in MacNeille completion
by join density, and in canonical completion by the closed-product
formula followed by join density. Thus the interior-support
projection $B\to pB$ is full. The old $G$ is a full subalgebra
of $pB$ for its local constants $p,q$. To see the residual claim
explicitly, for $g,h\in G$ and $z\in pB$, multiplication by $g$
has order-preserving inverse multiplication by the retained
$g^{-1}$. Hence $gz\le h$ exactly when $z\le g^{-1}h$,
and the actual residuals are $g\backslash h=g^{-1}h$ and
$h/g=hg^{-1}$. Its local prime is $g^{-1}q$; lattice operations
select one argument. Therefore
\[
                         E=\{x\in B:px\in G\}
 \tag{L27}\label{sl:lower:pullback}
\]
is a full original-signature subalgebra containing $K$, with $pE=G$.

If $a\prec_E b$ and $pa=pb=g$, the whole interval $[a,b]_B$
is in $E$, so it is also a gap of $B$. Its greatest self-support
$s$ satisfies $sb=b$. Applying the projection and multiplying by
the retained inverse of $g$ gives $ps=p$, so $s\in E$. It is
therefore the greatest self-support in $E$ as well.

This internal-support property is first-order in the original
finite signature: use
\[
 \begin{split}
 \operatorname{Gap}(a,b)&\equiv a<b\ \&\
           \neg\exists z\,(a<z<b),\\
 \operatorname{Fix}_b(u)&\equiv e\le u\ \&\ ub=b\ \&\ (u\backslash e)b=b,
 \end{split}
\]
and the sentence asserting that, whenever
$\operatorname{Gap}(a,b)$ and $pa=pb$, the set defined by
$\operatorname{Fix}_b$ has a greatest member. A countable elementary
substructure $C\preccurlyeq E$ containing $K$ exists by downward
L\"owenheim--Skolem. It is a full actual chain, satisfies this
sentence and every history, and has $pC=G$, since
$K\subseteq C\subseteq E$.
\end{proof}

\begin{lemma}[Dense floor preparation]
\label{sl:lower:preparation}
Every countable lower complete core has a nontrivial countable full
extension $A\in\KR$ with a densely ordered floor group, internal
gaps, and (PC).
\end{lemma}
\begin{proof}
Take $C$ from Lemma~\ref{sl:lower:internal-PC} and let
$H=G\overrightarrow\times\mathbb Q$, with componentwise group
product and lexicographic order. Embed $g$ as $(g,0)$. This is a
full group embedding for local identity $p$ and central dualizer
$q=(q,0)$. The group $H$ is countable, dense, without endpoints
and nontrivial. Apply Theorem~\ref{sl:struct:gluing} to $C$ and $H$
to form an extension $A$: the fibre of each old $g\in G$ is unchanged, while
every $h\in H\setminus G$ has the singleton fibre $\{h\}$.
Old products and prime are retained; a product involving a new
floor is the floor at the product of the degrees, and its prime is
$h\mapsto qh^{-1}$. The lifting theorem supplies both actual
residuals and every finite history for these precise rules.

Distinct endpoint degrees cannot form a gap, since an intermediate
degree in $H$ supplies an intermediate floor. Thus every gap lies
inside a single unchanged old fibre, and is a gap of $C$.
For a positive self-support $u$ of such a gap with upper degree
$g$, the equation $ub=b$ gives $\pi(u)g=g$, hence $\pi(u)=p$.
The neutral fibre is unchanged, so every possible self-support is
old. The greatest old one exists by the preceding lemma, and its
defining equations are preserved. Therefore every gap has (PC),
and all gaps are internal by their equal endpoint degrees.
\end{proof}

\begin{theorem}[Representation of lower complete cores]
\label{sl:lower:representation}
Every countable lower complete core $K\in\KR$ has a nontrivial
countable full dense extension in $\KR$ and a full embedding into
a standard chain.
\end{theorem}
\begin{proof}
Start with $A_0$ supplied by Lemma~\ref{sl:lower:preparation} and
put $A_{n+1}=D(A_n)$ with multiplication \eqref{sl:lower:P}.
Theorems~\ref{sl:lower:associativity} and
Lemmas~\ref{sl:lower:history-transfer}--\ref{sl:lower:invariants}
apply at every stage. Each stage is countable, each inclusion is
full, and the exact floor group is the same $H$ throughout.
The union $A_\infty=\bigcup_{n<\omega}A_n$ is a countable
member of $\KR$: each finite tuple, all its operations, and any
specified finite-history instance lie in one stage. For $x<y$
in this union, choose a common stage. Either it already has an
intermediate point or its next stage inserts one. Thus the union
is dense and nontrivial. Corollary~\ref{sl:standard-completion} gives
its full standard embedding. Composing the inclusions proves
the assertion for $K$.
\end{proof}

\begin{corollary}[Preservation of old central units]
\label{sl:lower:central-units}
The extensions and standard embedding in
Theorem~\ref{sl:lower:representation} can simultaneously preserve
the centrality of every old central global unit of $K$.
\end{corollary}
\begin{proof}
Let $u$ be any such unit. Fullness retains its old inverse. In a
MacNeille completion, multiplication by $u$ on either side agrees
on old points and preserves joins, so join density makes $u$
central everywhere. In a canonical extension,
\eqref{sl:lower:closed-product} gives
$uk=\bigwedge\{ua:a\in A,k\le a\}
 =\bigwedge\{au:a\in A,k\le a\}=ku$ for closed $k$;
join density again gives centrality on all points. The pullback
and elementary substructure inherit this identity.

In the floor enlargement, $\pi(u)$ is central in $G$, hence
$(\pi(u),0)$ is central in $G\overrightarrow\times\mathbb Q$.
Thus $u$ commutes with each new singleton floor by the explicit
product rule, and with each old point as before. At a midpoint
insertion, $u^-=u^{-1}$ and $u^-u=uu^-=e$. Both left and
right actions are therefore active at every gap; their equal old
upper values $ub=bu$ determine the same midpoint. Thus $u$
remains central at every stage, in the union, and in the final
MacNeille completion. These arguments apply to each old central
unit without adding any choices or countability requirements.
\end{proof}

\section{The upper-support construction}
\label{sl:upper:main-section}
A full embedding retains the comparison between $e$ and the closed
support $p$, so the lower construction cannot handle $e<p$.
We summarize the two upper constructions here; their complete proofs
are in Appendix~\ref{sl:upper:section}.

\begin{theorem}[Upper-support representation]
\label{sl:upper:main}
\label{sl:upper:representation}
Every countable $C\in\KR$ with $e\le p$ has a countable full dense
extension and a full standard embedding. The target retains $e\le p$
and every original finite history.
\end{theorem}
\begin{proof}
If $p=e$, then $d=q$ is a central global unit by
Theorem~\ref{sl:struct:support}. Apply
Theorem~\ref{sl:struct:unit-dense} and standard completion.
This includes a singleton input: then $d=e$, and the first odd-clone
step of Theorem~\ref{sl:struct:odd-clone} adds a point on each side
of its old point. Thus the dense extension is nontrivial before
standard completion. Fullness here concerns the original signature;
bounded embeddings preserving both endpoints are treated separately
in the application sections.

For $p>e$ put $I=pC$, $q=pd$, $f=p'$ and $n=q'$.
Lemma~\ref{sl:upper:corner} gives the local identity $p$ and
local dualizer $f$, the protected cover $n<p$ in $I$, and the
exhaustive alternatives $n\in U(I)$ and $n^2=n$.

In the invertible case, put $G=U(I)$ and
$K=\{x\in C:px\in G\}$, the complete core of
Lemma~\ref{sl:upper:exact-bands}. Theorem~\ref{sl:upper:neutral-maximum}
prepares a countable full extension of $K$ with exactly this group and a
greatest neutral live coefficient. Theorem~\ref{sl:upper:graded-dense}
then densifies its protected bands while retaining $G$ and their old
endpoints. Lemma~\ref{sl:upper:protected-glue} restores the remaining
support points. Corollary~\ref{sl:upper:invertible-core} assembles
these steps into a countable full dense extension of the entire input.

In the idempotent case, multiplication by $n$ is a full interior
projection with central unit local dualizer $nd=f$. In its complete
group core, Lemma~\ref{sl:upper:separate-bands} unfolds and refolds
the active bands, making them dense while retaining their two-point
common corner. Lemma~\ref{sl:upper:relative-corner} densifies the
remaining corner while fixing that common part. Theorem~\ref{sl:upper:idempotent-core}
glues these extensions and reconstructs the entire input, again giving
a countable full dense extension.

Standard completion supplies the claimed embeddings in both branches.
The dense extensions belong to $\KR$, so every finite history holds.
Fullness preserves the closed terms $e,p$ and hence their comparison
in the standard target.
\end{proof}

\section{Full standard representation}
\label{sl:representation-section}

The roles of completion and density are distinct.  Table~\ref{sl:tab:preservation}
records which support group is retained at each stage and when countability
and density are obtained.  All displayed extensions preserve the original
signature and the finite history scheme; this does not assert preservation
of arbitrary existing suprema or infima.

\begingroup
\small
\renewcommand{\arraystretch}{1.16}
\begin{longtable}{@{}>{\raggedright\arraybackslash}p{.22\textwidth}>{\raggedright\arraybackslash}p{.46\textwidth}>{\raggedright\arraybackslash}p{.24\textwidth}@{}}
\caption{Preservation and changes of support in the construction.  The
symbols $G$ and $H$ refer to the input and prepared groups of the
respective branch.}
\label{sl:tab:preservation}\\
\hline
Construction & Support and algebraic information & Size and density \\
\hline
\endfirsthead
\multicolumn{3}{l}{\tablename\ \thetable\ (continued)}\\
\hline
Construction & Support and algebraic information & Size and density \\
\hline
\endhead
\hline
\endfoot
\hline
\endlastfoot
MacNeille completion
(\ref{sl:completion})
& Fullness and all histories are retained.  The support and its unit
group need not remain unchanged.
& May be uncountable; adjacent pairs are not removed. \\
Canonical preparation
(\ref{sl:lower:canonical})
& Self-support maxima are attained in the preparatory extension.
No equality with the original support group is asserted here.
& The intermediate extension need not be countable or dense. \\
Lower pullback and elementary subchain
(\ref{sl:lower:internal-PC})
& Pullback over the original $G$, followed by a countable elementary
subchain, gives exactly $pC=G$ and retains the required maxima.
& Countability is obtained at this step; density is not yet asserted. \\
Lexicographic preparation
(\ref{sl:lower:preparation})
& The base is deliberately enlarged from $G$ to
$H=G\mathbin{\overrightarrow\times}\mathbb Q$; old fibres are retained.
& $H$ is countable and dense.  Internal fibre gaps may remain. \\
Iterated midpoint insertion
(\ref{sl:lower:invariants}, \ref{sl:lower:representation})
& Every insertion and their union have the same prepared floor group $H$.
& The union is countable, dense, and has no endpoints. \\
Upper filter preparation
(\ref{sl:upper:neutral-maximum})
& After completion, pullback and generated-subalgebra selection restore
exactly the original $G$; a greatest neutral coefficient is obtained.
& Intermediate completions may be uncountable; the selected extension
is countable. \\
Upper band densification
(\ref{sl:upper:graded-dense})
& The support group remains exactly $G$ and the protected endpoints
are retained.
& The resulting extension is countable and every protected closed
band is dense. \\
Whole-chain reconstruction
(\ref{sl:struct:reconstruction}, \ref{sl:upper:protected-glue})
& The local support is extended.  Interior reconstruction need not
keep its unit group unchanged; the relative protected construction
does keep the old unit group.
& The reconstructed chain is countable, dense, and has no endpoints. \\
Standard completion
(\ref{sl:standard-completion})
& The final embedding is full; equality of the final support group
with an earlier group is not claimed.
& A nontrivial countable dense input yields a chain ordered as $[0,1]$. \\
\end{longtable}
\endgroup

Table~\ref{sl:tab:dependencies} gives the main dependency interfaces.
In particular, the lower-core insertion proof is independent of the
upper construction.  The odd-chain results supply support densification
for reconstruction and for the relative upper construction, rather than
an additional hypothesis on the lower input.

\begin{table}[tbp]
\centering
\small
\renewcommand{\arraystretch}{1.16}
\begin{tabular}{@{}>{\raggedright\arraybackslash}p{.25\textwidth}>{\raggedright\arraybackslash}p{.67\textwidth}@{}}
\hline
Result & Inputs used by the proof \\
\hline
Lower-core representation
(\ref{sl:lower:representation})
& Full completion, canonical preparation, support-base lifting, and
the lower insertion lemmas.  Final standardization uses
Corollary~\ref{sl:standard-completion}. \\
Upper invertible branch
(\ref{sl:upper:neutral-maximum}--\ref{sl:upper:protected-glue})
& Proper-filter preparation and intrinsic band densification, followed
by the relative odd-chain support extension
(\ref{sl:struct:relative-unit}) and protected gluing. \\
Upper idempotent branch
(\ref{sl:upper:idempotent-core})
& Unfolding into the upper invertible construction, the relative lower
corner (\ref{sl:upper:relative-corner}), and interior reconstruction.
The relative lower corner uses lower-core representation. \\
Full representation
(\ref{sl:rep:standard})
& The upper theorem when $e\le p$.  When $p<e$, the lower core is
represented first, then the whole chain is restored by
Corollary~\ref{sl:struct:core-standard-reconstruction}; this restoration
uses odd-chain support densification and full completion. \\
Logical consequences
(\ref{sl:logic:strong-standard}, \ref{sl:logic:den-elimination})
& Relative countable comparison completeness and full representation
give strong standard completeness; density soundness gives Den admissibility.
The IUL and pseudo-uninorm transfers additionally use the stated
signature-preservation results and ordinary theory-chain completeness. \\
\hline
\end{tabular}
\caption{Dependencies between the algebraic construction and its
logical applications.}
\label{sl:tab:dependencies}
\end{table}

\begin{theorem}\label{sl:rep:standard}
Every countable member of $\KR$ embeds fully into a standard cyclic
involutive residuated chain. In particular, multiplication, both
residuals, involution and both designated constants are preserved
simultaneously.
\end{theorem}
\begin{proof}
Let $L\in\KR$ be countable and let $p$ be its closed central support.
If $e\le p$, apply Theorem~\ref{sl:upper:representation}.
Suppose $p<e$. The full interior projection $x\mapsto px$ has image
$I=pL$, whose central local dualizer is $q=pd$. Let $G=U(I)$ and
let $K=\{x\in L:px\in G\}$. Lemma~\ref{sl:struct:core} makes
$K$ a countable full subalgebra with $pK=G$. It is therefore a lower
complete core, to which Theorem~\ref{sl:lower:representation} applies.
Corollary~\ref{sl:struct:core-standard-reconstruction} reconstructs the
entire chain $L$ from that standard core embedding and a dense extension
of the support. All the embeddings used in this reconstruction are full.
\end{proof}

The distinction between the two comparisons with $p$ is exhaustive
because the order is total. The preparatory properties used in the
lower construction, such as attained self-supports and a dense support
group, are properties of constructed extensions. They impose no
additional hypotheses on $L$.

\section{Strong standard completeness and elimination of density}
\label{sl:logic:consequences}

We now use the full representation theorem proved in this paper:
every countable member of $\KR$ embeds, in the entire signature,
into a standard chain (Theorem~\ref{sl:rep:standard}).
An embedding preserves the values of all terms under one valuation,
including residuals, involution, and all closed history terms.

\begin{corollary}[Finite strict gaps]\label{sl:logic:finite-gaps}
Suppose a finite set $\Sigma$ of term comparisons, together with
$a<b$, is realized in a member of $\KR$.
If $\rho$ is a fresh variable, then
\[
                  \Sigma\cup\{a<\rho,\ \rho<b\}
\]
is realized in a standard member of $\KR$.
Here comparisons may be equations, weak inequalities, or strict
inequalities, and every term is evaluated under the same valuation.
\end{corollary}
\begin{proof}
Take the full subalgebra generated by the finitely many values of
variables in the pattern.  It is countable, since the signature is
finite, and is a nontrivial member of $\KR$: the defining conditions
are universal and the strict pair is retained.
Apply Theorem~\ref{sl:rep:standard}.  The embedding preserves all
the old term values and strict comparisons.  Choose the value of
$\rho$ strictly between the images of $a$ and $b$ in $[0,1]$.
Freshness leaves the old pattern unchanged, and the standard target
belongs to $\KR$ by Lemma~\ref{sl:logic:dense-histories}.
\end{proof}

\begin{theorem}[Strong standard completeness]
\label{sl:logic:strong-standard}
For every set $T$ of formulas in the fixed countable propositional
language and every finite hypersequent $H$,
\[
 T\vdash_{\HR}H
 \quad\Longleftrightarrow\quad T\models_{\Std}H.
\]
In particular, for every formula $A$,
\[
 T\vdash_{\HR}A
 \quad\Longleftrightarrow\quad
 \forall S\in\Std\ \forall v,\quad
 \bigl[(\forall B\in T,\ e_S\le v(B))
       \Longrightarrow e_S\le v(A)\bigr].
\]
\end{theorem}
\begin{proof}
The rule checks in Theorem~\ref{sl:logic:soundness} are pointwise
in a valuation. Induction on a proof therefore proves soundness
with the specified assumption leaves.

Conversely, suppose $T\nvdash_{\HR}H$.
Theorem~\ref{sl:logic:theory-chain-completeness} supplies a countable
chain $C\in\KR$ and a valuation $v$ satisfying $T$ and falsifying
$H$. By Theorem~\ref{sl:rep:standard}, choose a full embedding
$j:C\hookrightarrow S$ into a standard chain.
For $v'=j\circ v$, induction on formulas gives $v'(A)=j(v(A))$.
Since $j(e_C)=e_S$ and $j$ preserves and reflects order,
$v'$ satisfies every member of $T$.
Fullness likewise preserves and reflects the truth of each
component of $H$, so $v'$ falsifies $H$.
Thus $T\not\models_{\Std}H$.
All premises are retained by this single embedding; no
preservation of infinite suprema or infima is used.
\end{proof}

\begin{theorem}[Standard theoremhood]\label{sl:logic:standard}
For every finite hypersequent $H$,
\[
                         \HR\vdash H
            \quad\Longleftrightarrow\quad \Std\models H.
\]
\end{theorem}
\begin{proof}
Take $T=\varnothing$ in Theorem~\ref{sl:logic:strong-standard}.
The empty hypersequent remains false and is underivable by soundness.
\end{proof}

\begin{remark}[Scope of strong completeness]
The strong result concerns the calculus $\HR$ and the class
$\KR$; it does not remove the unit-cyclicity or history conditions
from the input representation theorem. All standard cyclic
involutive residuated chains belong to $\KR$ by
Lemma~\ref{sl:logic:dense-histories}, which justifies the use of
$\Std$ here. Standard consequence quantifies over the class of
standard algebras, not over one fixed algebra on $[0,1]$.
\end{remark}

Consider the density rule
\begin{equation}\label{sl:logic:den-rule}
 \frac{G\mid\seq{\Gamma,\rho}\mid\seq{\rho',\Delta}}
      {G\mid\seq{\Gamma,\Delta}}\;\mathrm{Den},
 \qquad \rho\notin G,\Gamma,\Delta.
\end{equation}
The freshness condition includes occurrences of either polarity
and occurrences inside formulas.

\begin{lemma}\label{sl:logic:den-soundness}
Den preserves validity over every densely ordered cyclic involutive
residuated chain.
\end{lemma}
\begin{proof}
If its conclusion is false under a valuation, every component of
$G$ is false and $\gamma\delta>d$, where
$\gamma=\Gamma^\sharp$ and $\delta=\Delta^\sharp$.
By residuation this means $\delta'<\gamma$.
Choose $\delta'<r<\gamma$ and set $v(\rho)=r$.
The two new active components would mean, respectively,
$\gamma r'\le d$, equivalently $\gamma\le r$, and
$r\delta\le d$, equivalently $r\le\delta'$.
Both are false.  Freshness retains the falsity of $G$ and the old
word values, so the premise is false under the extended valuation.
\end{proof}

\begin{theorem}[Unrestricted density elimination]\label{sl:logic:den-elimination}
The rule \eqref{sl:logic:den-rule} is admissible in $\HR$.
Moreover, from any finite $\HR+\mathrm{Den}$ proof one can
effectively obtain a finite Den-free $\HR$ proof of the same
end hypersequent.  Ordinary Cut remains available.
\end{theorem}
\begin{proof}
If the Den premise is derivable in $\HR$, it is standard-valid
by soundness.  Lemma~\ref{sl:logic:den-soundness} makes its
conclusion standard-valid, and Theorem~\ref{sl:logic:standard}
makes that conclusion derivable in $\HR$.  This is admissibility
for arbitrary premise proofs, with no restriction on their use of
Cut, EC, COM, shared conjunction, or any $R_w$.

For effectiveness, the normalized formula syntax and all finite
words $w\in\{S,T\}^*$ are computably enumerable.
The recursion \eqref{sl:logic:history-terms} computes the two
guards of each $R_w$.  Finite proof trees annotated by their rule
names, formula data, word splittings, and history indices have a
decidable correctness check.  Thus all finite $\HR$ proofs can
be effectively enumerated, including all choices of Cut formulas.

Proceed upward through a supplied finite $\HR+\mathrm{Den}$
proof.  Retain an ordinary inference after replacing its premise
proofs.  At a Den inference, admissibility guarantees an original
$\HR$ proof of its syntactic conclusion.  Enumerate original
proofs until one with that conclusion appears, and insert it.
Every such search terminates.  There are finitely many Den nodes,
so this produces the required finite proof effectively.
\end{proof}

The Den-elimination theorem above concerns proofs without additional
assumption leaves. Strong completeness with assumptions is stated
separately in Theorem~\ref{sl:logic:strong-standard}.
No claim is made here of effective Den elimination relative to an
arbitrary theory, elimination of ordinary Cut, or a complexity bound
for the proof-enumeration procedure. Relative density elimination
would also require the appropriate freshness conditions on assumptions.

\section{Transfer to involutive uninorm logic}
\label{sl:iul:section}

We use $\mathsf{IUL}$ for involutive uninorm logic and
$\mathsf{IUL}^{fp}$ for its extension by $e\leftrightarrow d$.
An $\mathsf{IUL}$-chain is a bounded commutative involutive residuated
chain in the signature
\[
 (\wedge,\vee,\cdot,\to,e,d,\bot,\top),
 \qquad x'=x\to d,\qquad x''=x.
\]
The designated condition for a formula is $e\le v(\varphi)$; neither
$e$ nor $d$ is required to be a lattice endpoint. In particular, no
equation $e=d$ is imposed. Standard chains have ordered carrier
$[0,1]$, with $\bot=0$ and $\top=1$.

Jenei's finite strong standard completeness theorem concerns
$\mathsf{IUL}^{fp}$~\cite[Theorem~2]{JeneiFixedPoint2023}; his later
densification theorem yields strong standard completeness for the same
fixed-point extension~\cite[Theorem~20]{JeneiDensification2025}.
The distinction is essential for embeddings. If $e\ne d$ in an input
chain, an embedding preserving both constants cannot have an odd target:
the target equation $e=d$ would identify their two images.

We apply the full representation theorem to general $\mathsf{IUL}$.
The representation constructions already supply comparison maxima,
associativity, cyclic duality, actual residuals, and density. For this
application we additionally prove that they preserve commutativity,
and then recover the two lattice-bound constants by restriction to
the old endpoint interval. These steps give the corollary below.
The proof-theoretic argument of~\cite{WangIUL2019} is not used.

\subsection{The history conditions on commutative inputs}

\begin{lemma}[Automatic histories]
\label{sl:iul:histories}
Every commutative involutive residuated chain belongs to $\KR$.
In particular, the reduct of every $\mathsf{IUL}$-chain obtained by
forgetting $\bot,\top$ belongs to $\KR$. No density or oddness
assumption is needed.
\end{lemma}
\begin{proof}
Unit cyclicity is immediate from commutativity. The involutive residual
formula gives $e+x=(dx')'=d\to x$. Thus the closed history terms of
\eqref{sl:logic:history-terms} satisfy
\[
\begin{gathered}
 (c_\epsilon,h_\epsilon)=(e,e),\\
 (c_{wS},h_{wS})=(dc_w,d\to h_w),\qquad
 (c_{wT},h_{wT})=(d\to c_w,dh_w).
\end{gathered}
\]
We prove $c_wh_w\le e$ by induction on $w$. The empty history gives
equality. For the two induction steps, commutativity and residuation give
\[
 \begin{split}
 c_{wS}h_{wS}
   &=c_w\bigl(d(d\to h_w)\bigr)\le c_wh_w,\\
 c_{wT}h_{wT}
   &=\bigl(d(d\to c_w)\bigr)h_w\le c_wh_w.
 \end{split}
\]
If $ab\le c_w$ and $uv\le h_w$, then
\[
 bazvu=(ab)(uv)z\le c_wh_wz\le z.
\]
This is precisely \eqref{sl:logic:history-condition}, for every finite
history and every choice of the five elements.
\end{proof}

\subsection{Commutativity of the extension constructions}

Fullness of an embedding preserves products of old elements, but does
not by itself impose commutativity on new elements. We therefore check
the actual constructions, including products of two new points.

\begin{lemma}[Symmetry of the lower insertion]
\label{sl:iul:lower-symmetry}
Suppose the old chain $A$ in Lemmas~\ref{sl:lower:actions}--\ref{sl:lower:selection}
is commutative and satisfies the hypotheses ensuring the existence of
the minima in \eqref{sl:lower:k}. Then the multiplication $P$ of
\eqref{sl:lower:P} is commutative.
\end{lemma}
\begin{proof}
For an old coefficient $u$ and a gap with upper endpoint $b$, left
activity means $u^-ub=b$, while right activity means $buu^-=b$.
These conditions coincide, and the old upper values $ub$ and $bu$
coincide. The prescribed actions in
Lemma~\ref{sl:lower:actions} therefore satisfy
\[
 uX=Xu\qquad(u\in A,\ X\in D(A)).
\]
For an old $z$, their right-adjoint formulas yield
\[
 \begin{split}
 z'X\le Y'
 &\iff X\le (Yz')'\\
 &\iff Yz'\le X'\\
 &\iff z'Y\le X'.
 \end{split}
\]
Consequently the two sets whose minima define $k(X,Y)$ and $k(Y,X)$
are equal. The assumed attainment gives $k(X,Y)=k(Y,X)$.
Since $r(X),r(Y)$ are old, the other two entries in the minimum are
exchanged when $X,Y$ are exchanged. Hence
\[
 \begin{split}
 P(X,Y)&=\min\{k(X,Y),r(X)Y,Xr(Y)\}\\
       &=\min\{k(Y,X),r(Y)X,Yr(X)\}=P(Y,X).
 \end{split}
\]
This argument proves symmetry; it does not infer attainment or
associativity from symmetry.
\end{proof}

\begin{proposition}[Commutative specialization of the constructions]
\label{sl:iul:commutative-constructions}
For a commutative input, the completion, preparation, insertion, and
reconstruction procedures used in Theorem~\ref{sl:rep:standard}
can be performed with commutative intermediate algebras and a
commutative standard target.
\end{proposition}
\begin{proof}
We verify preservation in the constituent constructions.

\emph{Completions and restrictions.}
In a MacNeille completion, separate preservation of joins and join
density extend the old equation $xy=yx$ to all pairs of completed
elements. In the canonical extension, the closed-product formula
\eqref{sl:lower:closed-product} is symmetric in its two arguments;
expansion as joins of closed elements extends this symmetry to the
whole algebra. In the proper-filter completion of
Lemma~\ref{sl:upper:proper-filter}, old commutativity gives
\[
 F*G=\mathord\uparrow\{ab:a\in F,b\in G\}=G*F,
\]
so the induced product of sets of filters is symmetric as well.
Subalgebras, the indicated inverse images, elementary substructures,
and increasing unions inherit commutativity.

\emph{Odd extensions and support groups.}
For a commutative odd input, the local identities in
Appendix~\ref{sl:struct:odd-appendix} satisfy $\ell(x)=\rho(x)=u$.
Every arrow is therefore an endomorphism at its local identity.
For $z>u$, the absorption identities give
\[
 \alpha_x(z)=xzx'=zxx'=zu'=z,
\]
and $\alpha_x(u)=u$. Thus every tail map is the identity.
The groups of matching isomorphisms and the quotient groups are
Abelian; old matching products are commutative. The new matching
product \eqref{sl:struct:new-matched} uses their commutative product
and the symmetric operation $\min(i,j)$ on marks. In
\eqref{sl:struct:clone-product}, both arguments are lifted to the
same maximum of their two local identities; interchanging the
arguments therefore leaves the product unchanged. The full clone,
its relative restriction in Corollary~\ref{sl:struct:relative-odd},
and their iterations are commutative. Normalizing or restoring a
central unit dualizer changes prime but not multiplication.
The degree groups obtained from commutative inputs are Abelian,
as are the lexicographic enlargements by $\mathbb Q$ or $\mathbb Z$.

\emph{Lower cores and interior gluing.}
The preceding completion and restriction arguments apply to
Lemma~\ref{sl:lower:internal-PC}. The dense floor preparation uses
an Abelian enlarged group. At each subsequent insertion,
Lemma~\ref{sl:iul:lower-symmetry} applies; its attainment hypotheses
are supplied by the cited preparation and invariant lemmas.
In Theorem~\ref{sl:struct:gluing}, a product involving an outside
factor has the form $i(\Pi(x)\Pi(y))$, which is symmetric when the
support algebra is commutative. Products in the core already commute.
This covers the lower representation and its whole-chain
reconstruction, including the central-unit support extension.

\emph{Upper coefficients and middle intervals.}
Use the notation of Lemma~\ref{sl:upper:graded-coefficients}.
Commutativity of $U$ gives $p_u=q_u$ and
\[
 L(u,v)\iff R(v,u).
\]
Its raw product
\[
 (u,r)(v,s)=
 \left(uv,\min\bigl(\{1\}\cup\{r:R(u,v)\}
                         \cup\{s:L(u,v)\}\bigr)\right)
\]
is consequently symmetric. The closed-cut product
$I\otimes K=j(\mathord\downarrow IK)$ remains symmetric, as does
its restriction to the selected degree fibres.

For the middle action extensions, equal old left and right actions
have equal activity tests in Lemma~\ref{sl:upper:action-gap}.
They transport a point to the same labelled point or collapse it to
the same old anchor. The completion in
Lemma~\ref{sl:upper:complete-action} and the joins in
\eqref{sl:upper:cut-actions} therefore preserve $\lambda_A=\rho_A$.
Cyclic adjunction gives, for every old coefficient $u$,
\[
 uX\le Y'\iff Yu\le X'\iff uY\le X'.
\]
If $X,Y$ have degrees $g,k$, the boundary degree
$b=(gk)^{-1}\delta$ is unchanged by interchanging them, because
the degree group is Abelian. Their boundary comparison sets are
equal, including the case in which they are empty. Whenever the
comparison element is defined, it follows that
\[
 \kappa_b(X,Y')=\kappa_b(Y,X').
\]
In \eqref{sl:upper:middle-product}, the old anchor product is
symmetric, the comparison branch has just been shown symmetric,
and the remaining branch is the support top of degree $gk=kg$.
All three cases are therefore symmetric. Coefficient--middle
products commute because $\lambda_A=\rho_A$; coefficient--reverse
products commute because the two residuals of the commutative
coefficient algebra agree. The other prescribed products are
support tops of the same product degree. This verifies
commutativity in Theorem~\ref{sl:upper:graded-dense}.

\emph{Upper reconstruction.}
Adjoining external bounds uses symmetric absorbing-product rules.
The neutral-maximum preparation uses the proper-filter completion
and subalgebras already treated above. Protected top gluing in
Lemma~\ref{sl:upper:protected-glue} evaluates a mixed product as
$\Pi(x)\Pi(y)$ in a commutative support algebra. In the idempotent
upper branch, the unfolded product
\eqref{sl:upper:unfold-2} is symmetric over the Abelian group
$H\overrightarrow\times\mathbb Z$. Its upper-invertible
extension has just been treated; restriction to the retained
positive bands preserves symmetry, and the separate floor rules
\eqref{sl:upper:unfold-5} are symmetric over $H$.
The deletion in Lemma~\ref{sl:upper:relative-corner} is a full
subalgebra restriction. Finally, in
\eqref{sl:upper:glue-3}, the two mixed products are $(pa)x$ and
$x(pa)$ inside the commutative corner extension, so they agree.
This covers both upper branches and their whole-chain gluing.

At every use, the cited results supply the existence, fullness,
associativity, residuals, and density required for that construction;
the preceding calculations supply commutativity. The final
MacNeille completion preserves it by the first paragraph.
\end{proof}

\subsection{Preserving the lattice bounds}

\begin{lemma}[Restriction to the old endpoint interval]
\label{sl:iul:endpoints}
Let $A$ be a nontrivial bounded commutative involutive residuated
chain. Suppose $i:A\hookrightarrow S$ preserves order, multiplication,
prime, $e$, and $d$, where $S$ is a commutative standard chain; the
embedding is not assumed to preserve the bounds of $S$.
Then $A$ has an embedding into a standard $\mathsf{IUL}$-chain
preserving the entire bounded signature.
\end{lemma}
\begin{proof}
Put $a=i(\bot_A)$ and $b=i(\top_A)$. The bottom of a bounded
residuated lattice is absorbing, and $\top_Ae=\top_A$ with
$e\le\top_A$ implies $\top_A^2=\top_A$. Thus
\[
 a^2=a,\qquad b^2=b,\qquad a'=b,\qquad b'=a.
\]
The closed real interval $S_0=[a,b]$ contains $e_S=i(e_A)$ and
$d_S=i(d_A)$. For $x,y\in S_0$, monotonicity gives
\[
 a=a^2\le xy\le b^2=b.
\]
The involution exchanges the endpoints and therefore preserves
$S_0$. Since $x\to y=(xy')'$, the interval is closed under the
ambient residual itself, as well as under the inherited lattice
operations. It is a full subalgebra with bottom $a$ and top $b$;
no truncation or alteration of the residual is involved.

Nontriviality and injectivity give $a<b$. Transport the operations
of $S_0$ along $t\mapsto(t-a)/(b-a)$. The resulting algebra has
carrier $[0,1]$, and the composite embedding preserves all
operations and all four constants $e,d,\bot,\top$.
\end{proof}

\subsection{Standard embeddings and strong standard completeness}

\begin{corollary}[Bounded standard embeddings for IUL]
\label{sl:iul:bounded-embedding}
Every nontrivial countable $\mathsf{IUL}$-chain embeds into a
standard $\mathsf{IUL}$-chain preserving its entire bounded signature.
\end{corollary}
\begin{proof}
Let $A$ be a nontrivial countable $\mathsf{IUL}$-chain.
Lemma~\ref{sl:iul:histories} places its reduct in $\KR$.
Theorem~\ref{sl:rep:standard}, specialized as in
Proposition~\ref{sl:iul:commutative-constructions}, supplies a full
embedding of that reduct into a commutative standard chain.
Lemma~\ref{sl:iul:endpoints} then restores preservation of the
bounded signature.
\end{proof}

For any class of algebras in the signature of this section,
$T\models\varphi$ denotes consequence under the designated
condition $e\le v(\psi)$, with all premises interpreted by the
same valuation. The same convention will be used for pseudo-uninorm
logic in the next section.

\begin{corollary}[Strong standard completeness for IUL]
\label{sl:iul:strong-standard}
In a countable propositional language, for every set $T$ of formulas
and every formula $\varphi$,
\[
 T\vdash_{\mathsf{IUL}}\varphi
 \quad\Longleftrightarrow\quad
 T\models_{\operatorname{Std}(\mathsf{IUL})}\varphi.
\]
No fixed-point axiom $e\leftrightarrow d$ is imposed.
\end{corollary}
\begin{proof}
Use ordinary chain completeness for arbitrary theories
\cite[Theorem~1, with Definitions~2 and~4--6]{WangFinite2018}.
Its constant $f$ is our $d$, and commutativity identifies its two
residuals. If $T\nvdash_{\mathsf{IUL}}\varphi$, there are an
$\mathsf{IUL}$-chain $A$ and a valuation $v$ with
$e_A\le v(\gamma)$ for every $\gamma\in T$, whereas
$v(\varphi)<e_A$.
Let $A_0$ be the subalgebra generated by the values of all
propositional variables and all signature constants, including
$\bot$ and $\top$. It is countable because the language is countable
and the signature is finitary. It remains an $\mathsf{IUL}$-chain,
retains all these inequalities, and is nontrivial since
$v(\varphi)<e_A$.

Corollary~\ref{sl:iul:bounded-embedding} supplies a full bounded
standard embedding $j:A_0\hookrightarrow S$.
For $v_S=j\circ v$, with $v$ now valued in $A_0$, induction on
formulas gives $v_S(\psi)=j(v(\psi))$ for every formula $\psi$.
Thus $e_S\le v_S(\gamma)$ for every $\gamma\in T$, while
$v_S(\varphi)<e_S$. This is a standard countermodel.
The reverse implication follows from soundness.
Only the ordinary theory-chain-completeness result is used from
the cited source, not a standard-completeness result.
\end{proof}

\begin{corollary}[Standard embeddings and standard completeness for IUL]
\label{sl:iul:standard-transfer}
Every nontrivial countable
$\mathsf{IUL}$-chain embeds into a standard $\mathsf{IUL}$-chain
preserving its entire bounded signature. Consequently, using the
ordinary chain completeness of $\mathsf{IUL}$,
\[
 \mathsf{IUL}\vdash\varphi
 \quad\Longleftrightarrow\quad
 \text{$\varphi$ is valid in every standard $\mathsf{IUL}$-chain}.
\]
\end{corollary}
\begin{proof}
The embedding assertion is Corollary~\ref{sl:iul:bounded-embedding};
the theoremhood equivalence is the case $T=\varnothing$ of
Corollary~\ref{sl:iul:strong-standard}.
\end{proof}

\begin{remark}[Scope of the transfer]
\label{sl:iul:scope}
No fixed-point or additional history axiom is imposed. Histories
are automatic by Lemma~\ref{sl:iul:histories}; commutativity of the
target and preservation of its bounds require the separate
construction and endpoint arguments above. The logical transfer
uses ordinary chain completeness for arbitrary theories, retaining
all premises in one countably generated countermodel and one standard
embedding; it does not assume compactness of standard semantics.
\end{remark}

\section{Transfer to pseudo-uninorm logic}
\label{sl:psul:section}

We use the convention
\[
 ab\le c\quad\Longleftrightarrow\quad b\le a\backslash c
 \quad\Longleftrightarrow\quad a\le c/b.
\]
A $\mathsf{psUL}^{*}$-chain is a bounded residuated chain
\[
 A=(A,\wedge,\vee,\cdot,\backslash,/,e,\bot_A,\top_A)
\]
satisfying unit cyclicity
\begin{equation}\label{sl:psul:eq:uc}
 ab\le e\quad\Longleftrightarrow\quad ba\le e.
\end{equation}
In Wang's notation~\cite{WangPseudo2019}, $a\to b=a\backslash b$,
$a\rightsquigarrow b=b/a$, and $t=e$. The added weak commutativity
axiom is equivalent to $e/a\le a\backslash e$ for every $a$.
Indeed, this inequality sends $ba\le e$ to $ab\le e$;
interchanging $a,b$ gives the converse. Thus it is equivalent to
\eqref{sl:psul:eq:uc}. No involution is assumed in $A$.

The empty-theory target is the theoremhood statement of
\cite[Theorem~5]{WangPseudo2019}; Corollary~\ref{sl:psul:strong-standard}
below strengthens it to arbitrary theories:
\begin{equation}\label{sl:psul:eq:target}
 \mathsf{HpsUL}^{*}\vdash\varphi
 \quad\Longleftrightarrow\quad
 v(\varphi)\ge e\text{ in every standard }\mathsf{psUL}^{*}\text{-chain}.
\end{equation}
Here ``standard'' means that the lattice reduct is $[0,1]$, with its
usual bounds. This convention does not require the monoid unit to be
strictly below $1$. The star is essential: the target is $\mathsf{HpsUL}^{*}$,
not the unextended $\mathsf{HpsUL}$.

The representation input is Theorem~\ref{sl:rep:standard}, with
$\KR$ defined by \eqref{sl:logic:unit-cyclicity} and all ordered
history conditions \eqref{sl:logic:history-condition}. The transfer
below supplies the embedding and endpoint arguments needed to apply
that theorem to non-involutive inputs. It uses ordinary chain
completeness, but neither density elimination nor the standard
completeness result of~\cite{WangPseudo2019}.

\subsection{The upper-copy involutive extension}

The following extension is the construction of
Galatos--Raftery~\cite[Section~9, Theorem~9.1]{GalatosRaftery2004}.
We give the operation-level verification needed for the present application.

Take a disjoint order-dual copy $\overline A=\{\bar a:a\in A\}$ and
two new endpoints. Order
\[
 \mathsf{R}(A)=\{\mathbf0_R\}\mathbin{\dot\cup}A
 \mathbin{\dot\cup}\overline A\mathbin{\dot\cup}\{\mathbf1_R\}
\]
by $\mathbf0_R<a<\bar b<\mathbf1_R$ for all $a,b\in A$, and
$\bar a\le\bar b$ if and only if $b\le a$.
Define an order-reversing involution by
\[
 a'=\bar a,\quad(\bar a)'=a,\quad
 \mathbf0_R'=\mathbf1_R,\quad\mathbf1_R'=\mathbf0_R,
 \qquad d=\bar e.
\]
The product of internal elements is
\begin{equation}\label{sl:psul:eq:reflectionproduct}
 ab=(ab)_A,\qquad
 a\bar b=\overline{b/a},\qquad
 \bar a b=\overline{b\backslash a},\qquad
 \bar a\bar b=\mathbf1_R.
\end{equation}
Set $\mathbf0_R x=x\mathbf0_R=\mathbf0_R$ for all $x$, and
$\mathbf1_R x=x\mathbf1_R=\mathbf1_R$ whenever $x\ne\mathbf0_R$.
The zero rule takes precedence, so $\mathbf0_R\mathbf1_R=\mathbf1_R\mathbf0_R=\mathbf0_R$.

\begin{lemma}\label{sl:psul:lem:reflection}
For every residuated chain $A$, $\mathsf{R}(A)$ is a cyclic involutive
residuated chain. The inclusion of $A$ preserves its lattice operations,
multiplication, monoid unit, and both residuals. Furthermore,
\begin{equation}\label{sl:psul:eq:nzd}
 xy=\mathbf0_R\quad\Longleftrightarrow\quad x=\mathbf0_R\text{ or }y=\mathbf0_R.
\end{equation}
\end{lemma}
\begin{proof}
The identity $e$ acts as identity on both copies, since
$a/e=a=e\backslash a$, and on the endpoints by definition.
Equation~\eqref{sl:psul:eq:nzd} follows directly from the product definitions.

We check associativity exhaustively by the types of the three factors.
For three factors in $A$, use old associativity. When exactly one factor
is in $\overline A$, the three positions reduce to
\begin{align*}
 c/(ab)&=(c/b)/a,\\
 c\backslash(b/a)&=(c\backslash b)/a,\\
 (bc)\backslash a&=c\backslash(b\backslash a).
\end{align*}
For example, the middle identity holds because an element $t$ is below
either side exactly when $cta\le b$. The other two follow in the same
way from $tab\le c$ and $bct\le a$, respectively.
When at least two factors are in $\overline A$, both bracketings give
$\mathbf1_R$. A factor $\mathbf0_R$ makes both bracketings $\mathbf0_R$.
If $\mathbf1_R$ occurs without $\mathbf0_R$, both bracketings give $\mathbf1_R$, using
\eqref{sl:psul:eq:nzd}. This covers all triples.

We next check the cyclic comparison law
\begin{equation}\label{sl:psul:eq:cyclic}
 xy\le z'\quad\Longleftrightarrow\quad yz\le x'.
\end{equation}
If $x,y,z$ include $\mathbf0_R$, both comparisons hold. If they include
$\mathbf1_R$ but not $\mathbf0_R$, both comparisons fail, by \eqref{sl:psul:eq:nzd}.
For three factors in $A$, both comparisons hold, since $A<\overline A$.
For at least two factors in $\overline A$, both comparisons fail.
The remaining case, two old elements and one upper-copy element, is
exactly the chain of equivalences
\[
 ab\le c
 \quad\Longleftrightarrow\quad
 \overline{c/b}\le\bar a
 \quad\Longleftrightarrow\quad
 \overline{a\backslash c}\le\bar b.
\]
These are the three cyclic placements of $a,b,\bar c$.

Define
\begin{equation}\label{sl:psul:eq:residuals}
 x\backslash_R y=(y'x)',\qquad y/_R x=(xy')'.
\end{equation}
Applying \eqref{sl:psul:eq:cyclic} to $(x,t,y')$, twice if necessary, gives
\[
 xt\le y\quad\Longleftrightarrow\quad
 t\le(y'x)'\quad\Longleftrightarrow\quad x\le(ty')'.
\]
Thus these are actual left and right residuals, and in particular the
product is monotone in each coordinate. Moreover,
$x\backslash_R d=(ex)'=x'=d/_R x$, so $d$ is cyclic dualizing.
For old $a,b\in A$, the definitions give
\[
 a\backslash_R b=(\bar b a)'=a\backslash_A b,
 \qquad b/_R a=(a\bar b)'=b/_A a.
\]
The total order restricts to the old order, so the lattice operations
are also preserved. The old bounds need not be the new endpoints.
\end{proof}

\subsection{Why every ordered history is satisfied}

\begin{lemma}\label{sl:psul:lem:KR}
If $A$ satisfies unit cyclicity, then $\mathsf{R}(A)\in\KR$.
\end{lemma}
\begin{proof}
For $x,y\ne\mathbf0_R$, the comparison $xy\le e$ can hold only when both
$x,y\in A$: every other product lies in $\overline A$ or equals $\mathbf1_R$.
Hence unit cyclicity in $\mathsf{R}(A)$ follows from unit cyclicity in $A$;
the zero cases are symmetric.

Since $d=\bar e$, we have
\[
 d^2=\mathbf1_R,\qquad e+e=(d^2)'=\mathbf0_R,
 \qquad d\mathbf0_R=\mathbf0_R,\qquad
 e+\mathbf0_R=(d\mathbf1_R)'=\mathbf0_R.
\]
Consequently
\[
 (c_S,h_S)=(d,\mathbf0_R),\qquad
 (c_T,h_T)=(\mathbf0_R,d).
\]
The recursion~\eqref{sl:logic:history-terms} preserves a zero coordinate in either
position. Induction on the word length therefore yields
\begin{equation}\label{sl:psul:eq:zeroguard}
 w\ne\epsilon\quad\Longrightarrow\quad
 c_w=\mathbf0_R\text{ or }h_w=\mathbf0_R.
\end{equation}
If $ab\le c_w$ and $uv\le h_w$ for such a word, one of $ab,uv$
is $\mathbf0_R$. By \eqref{sl:psul:eq:nzd}, at least one of $a,b,u,v$ equals
$\mathbf0_R$, so $bazvu=\mathbf0_R\le z$.
For $w=\epsilon$, unit cyclicity gives
\[
 ab\le e,\ uv\le e\quad\Longrightarrow\quad
 ba\le e,\ vu\le e\quad\Longrightarrow\quad bazvu\le eze=z.
\]
Thus all histories, including the empty one, hold.
\end{proof}

This is the additional check that permits application of
Theorem~\ref{sl:rep:standard}.
The original chain need not carry an involution, a dualizer, or any
extra history axioms. The reflection is countable whenever $A$ is.

\subsection{Recovering the original bounds}

\begin{lemma}\label{sl:psul:lem:interval}
Let $A$ be a nontrivial bounded residuated chain and let $i:A\to S$
be an order embedding into a residuated chain on $[0,1]$ preserving
multiplication, both residuals, and $e$, but not necessarily the bounds.
Then $i$ factors through a full bounded embedding into a residuated
chain on $[0,1]$, obtained by restricting $S$ to the interval between
the images of the two old bounds and rescaling.
\end{lemma}
\begin{proof}
Put $\alpha=i(\bot_A)$ and $\beta=i(\top_A)$.
In every bounded residuated lattice the bottom is absorbing. Since
$e\le\top_A$, the top is idempotent. Preservation of old operations gives
\begin{equation}\label{sl:psul:eq:endpoints}
 \alpha^2=\alpha,\quad\beta^2=\beta,\quad
 \alpha\beta=\beta\alpha=\alpha,\quad
 \alpha\backslash_S\beta=\beta=\beta/_S\alpha.
\end{equation}
The real interval $J=[\alpha,\beta]$ contains $i(e)$. For $x,y\in J$,
\[
 \alpha=\alpha^2\le xy\le\beta^2=\beta.
\]
Moreover, both ambient residuals stay in $J$:
\begin{equation}\label{sl:psul:eq:intervalres}
 \alpha\le x\backslash_S y\le\alpha\backslash_S\beta=\beta,
 \qquad
 \alpha\le y/_S x\le\beta/_S\alpha=\beta.
\end{equation}
The upper bounds use the order properties of residuals. The lower
bounds follow respectively from
$x\alpha\le\beta\alpha=\alpha\le y$ and
$\alpha x\le\alpha\beta=\alpha\le y$.
Thus restriction retains the actual residuals; no truncation is used.
The lattice operations restrict automatically, and the new bounds are
$\alpha,\beta$. Nontriviality and injectivity imply $\alpha<\beta$.
Transport all operations along
\[
 \theta:J\longrightarrow[0,1],\qquad
 \theta(t)=\frac{t-\alpha}{\beta-\alpha}.
\]
Then $\theta\circ i$ preserves the entire bounded residuated signature.
\end{proof}

\begin{lemma}\label{sl:psul:lem:denseuc}
Every densely ordered residuated chain is unit-cyclic.
\end{lemma}
\begin{proof}
Suppose $ab\le e$ and $ba>e$. Residuation gives $b\backslash e<a$.
Choose $b\backslash e<r<a$ by density. Then $br>e$, and hence
\[
 a\le a(br)=(ab)r\le r<a,
\]
a contradiction. Interchange $a,b$ for the converse.
\end{proof}

\begin{corollary}[Full standard embedding for pseudo-uninorm chains]\label{sl:psul:thm:embedding}
Every nontrivial countable bounded unit-cyclic residuated
chain embeds into a standard $\mathsf{psUL}^{*}$-chain, preserving
\[
 \wedge,\vee,\cdot,\backslash,/,e,\bot,\top.
\]
\end{corollary}
\begin{proof}
By Lemmas~\ref{sl:psul:lem:reflection} and~\ref{sl:psul:lem:KR}, $A$ embeds as a
residuated reduct into the countable member $\mathsf{R}(A)$ of $\KR$.
Apply Theorem~\ref{sl:rep:standard} to obtain a full embedding
$\mathsf{R}(A)\to S$, where
$S$ has underlying order $[0,1]$. The composite preserves the old
residuated signature. Lemma~\ref{sl:psul:lem:interval} makes it preserve the
two old bounds as well. The resulting standard residuated chain is
unit-cyclic by Lemma~\ref{sl:psul:lem:denseuc}, and is therefore a standard
$\mathsf{psUL}^{*}$-chain.
\end{proof}

\begin{corollary}[Countable densification for pseudo-uninorm chains]\label{sl:psul:cor:dense}
Every nontrivial countable
$\mathsf{psUL}^{*}$-chain has a full bounded embedding into a
countable densely ordered $\mathsf{psUL}^{*}$-chain.
\end{corollary}
\begin{proof}
In the standard target of Corollary~\ref{sl:psul:thm:embedding}, take the
subalgebra generated by the image of $A$ together with
$\mathbb Q\cap(0,1)$. Both endpoints already belong to the image,
since the embedding preserves the bounds. The signature is finite,
so this subalgebra is countable. The included rational points make
its order dense. All defining operations, bounds, and unit cyclicity
are retained.
\end{proof}

\subsection{Strong standard completeness}

For the theory-level statement, the consequence relation is that
of the calculus in \cite[Definitions~2 and~6, Theorem~1]{WangFinite2018}.
For the displayed signature it is restricted to theories and
conclusions not containing the additional constant $f$ of that
presentation. The result also covers its full language with $f$:
interpret $f$ as an arbitrary distinguished element, include its
value among the generators below, and interpret it in the target
as its image under the embedding. It is not identified with the
dualizer of the auxiliary involutive extension.

\begin{corollary}[Strong standard completeness for $\mathsf{HpsUL}^{*}$]
\label{sl:psul:strong-standard}
In a countable propositional language, for every theory $T$ and
formula $\varphi$ in the chosen presentation just specified,
\[
 T\vdash_{\mathsf{HpsUL}^{*}}\varphi
 \quad\Longleftrightarrow\quad
 T\models_{\operatorname{Std}(\mathsf{psUL}^{*})}\varphi.
\]
In the presentation with $f$, the standard models in this formula
include all interpretations of that additional constant.
\end{corollary}
\begin{proof}
Ordinary chain completeness for arbitrary theories
\cite[Theorem~1, with Definitions~2 and~4--6]{WangFinite2018}
gives, from $T\nvdash_{\mathsf{HpsUL}^{*}}\varphi$, a
$\mathsf{psUL}^{*}$-chain $A$ and a valuation $v$ satisfying
$e_A\le v(\gamma)$ for all $\gamma\in T$ and
$v(\varphi)<e_A$.
Take the subalgebra $A_0$ generated by all variable values and all
constants of the chosen signature, including both bounds and,
when present, $f_A$. It is countable, nontrivial, remains
unit-cyclic, and retains all these premise and conclusion values.

By Corollary~\ref{sl:psul:thm:embedding}, its residuated reduct
has a standard embedding $j:A_0\hookrightarrow S$ preserving
both residuals, multiplication, lattice operations, the unit,
and the two bounds. If $f$ is present, set $f_S=j(f_A)$.
Induction on formulas then gives $v_S(\psi)=j(v(\psi))$ for
every formula $\psi$, where $v_S=j\circ v$.
In particular, all members of $T$ remain designated and
$v_S(\varphi)<e_S$. This refutes standard consequence.
The converse is soundness. No density-elimination or
standard-completeness theorem is used as an input.
\end{proof}

\begin{corollary}[Standard completeness for $\mathsf{HpsUL}^{*}$]\label{sl:psul:thm:complete}
The standard theoremhood equivalence \eqref{sl:psul:eq:target} holds.
\end{corollary}
\begin{proof}
Take $T=\varnothing$ in Corollary~\ref{sl:psul:strong-standard},
in the displayed signature without the additional constant $f$.
\end{proof}

\begin{remark}[Scope]
The reflection construction is due to \cite{GalatosRaftery2004};
the history and endpoint checks above make the representation theorem
applicable to it. Corollary~\ref{sl:psul:cor:dense} supplies the lower
densification link in Wang's Figure~4, while
Corollary~\ref{sl:psul:thm:embedding} gives a direct standard embedding.
The logical conclusions use ordinary theory-chain completeness,
not Wang's density-elimination or standard-completeness result.
The monoid unit may equal $1$, as it must for an integral input
whose bounds are preserved; this is allowed by the standard-algebra
convention in \cite[Theorem~5]{WangPseudo2019}.
\end{remark}

\newcounter{maintextequations}
\setcounter{maintextequations}{\value{equation}}
\appendix
\renewcommand{\thesection}{\Alph{section}}
\counterwithout{equation}{section}
\setcounter{equation}{\value{maintextequations}}
\renewcommand{\theequation}{\arabic{equation}}
\section*{Technical appendices}

The following appendices contain the detailed auxiliary proofs used
in the main text. Numbered references point to the main sections;
lettered references point to the appendices. All constructions use
the original signature and the class $\KR$ defined in
Section~\ref{sl:logic:calculus}.
No result from an unpublished working manuscript is used as a theorem.

\section{The countable comparison construction}
\label{sl:logic:chain-supplement}

For a set $\mathcal T$ of specified single-component sequents,
write $\mathcal T\vdash H$ if $H$ has a finite $\HR$ proof
allowing members of $\mathcal T$ as additional leaves.
An additional leaf is its specified sequent, not a substitution
schema.  Repeated uses of it are allowed.  The original axioms
retain all their instances.

\begin{lemma}[External continuation]\label{sl:logic:continuation}
For a specified sequent $S$ and finite hypersequents $J,K$,
\[
 \mathcal T\cup\{S\}\vdash K,\qquad
 \mathcal T\vdash J\mid S
 \quad\Longrightarrow\quad
 \mathcal T\vdash J\mid K.
\]
\end{lemma}
\begin{proof}
Induct on the first supplied proof, maintaining the entire multiset
$J$ in the constructed external background.  Replace each extra
leaf $S$ by the second supplied proof.  At every other leaf add
$J$ by EW.  Unary logical rules, Cyc, WCM, and EW retain $J$.
At EC contract the two occurrences selected by the old inference,
leaving the extra background intact.

At a product, Cut, or COM inference, the lifted premises supply two
copies of $J$; replay that inference and use EC on literally
identical components until the multiplicities of $J$ are restored.
The same construction contracts the three copies supplied at $R_w$.
For shared conjunction, the lifted premises have the identical
background $J\mid G$ and the identical ordered contexts; apply
the displayed $(\wedge)$ rule with that background.  No unequal
components are contracted or aligned.  These are all the rules.
The resulting proof is finite, since both supplied proofs are finite.
\end{proof}

\begin{lemma}[Comparison orientation]\label{sl:logic:orientation}
If $\mathcal T\nvdash H$, then for every pair of formulas $A,B$
at least one of
\[
 \mathcal T\cup\{\seq{A',B}\}\nvdash H,
 \qquad
 \mathcal T\cup\{\seq{B',A}\}\nvdash H
\]
holds.
\end{lemma}
\begin{proof}
Set $S=\seq{A',B}$ and $T=\seq{B',A}$.
COM applied to $\seq{A',A}$ and $\seq{B',B}$ proves $S\mid T$.
If both enlarged sets proved $H$, first lift the continuation from
$S$ over the full background $T$, obtaining $\mathcal T\vdash T\mid H$.
Lift the continuation from $T$ over the full background $H$ to
obtain $\mathcal T\vdash H\mid H$.
EC then gives $\mathcal T\vdash H$, a contradiction.
This uses both alternatives through their continuations; it never
projects one component from a proved hypersequent.
\end{proof}

\begin{lemma}[The comparison algebra]\label{sl:logic:comparison-algebra}
For any set $\mathcal T$ of specified single-component sequents, define
\begin{equation}\label{sl:logic:preorder}
 A\preccurlyeq B
 \quad\Longleftrightarrow\quad
 \mathcal T\vdash\seq{A',B}.
\end{equation}
This is a preorder whose symmetric part is a congruence for all
formula operations.  Its quotient is a cyclic involutive residuated
lattice with unit $[e]$ and dualizer $[d]$.
If every pair of formulas is comparable, the quotient is a chain.
\end{lemma}
\begin{proof}
Id proves reflexivity and Cut proves transitivity.  Cyc gives
\begin{equation}\label{sl:logic:order-reversal}
                 A\preccurlyeq B
        \quad\Longleftrightarrow\quad B'\preccurlyeq A'.
\end{equation}
We give the remaining calculations to specify the ordered quotient
without invoking a completeness theorem for another calculus.

The $(+)$ rule has a derivable inverse in every ordered context:
\begin{equation}\label{sl:logic:unfold}
 \mathcal T\vdash\seq{\Gamma,A+B,\Delta}
 \quad\Longleftrightarrow\quad
 \mathcal T\vdash\seq{\Gamma,A,B,\Delta}.
\end{equation}
Indeed, product applied to $\seq{A',A}$ and $\seq{B',B}$ proves
$\seq{(A+B)',A,B}$.  Rotate the supplied left component to
$\seq{\Delta,\Gamma,A+B}$, cut against that proof, and rotate
back.  The reverse implication is $(+)$.  No internal exchange
occurs in this calculation.

If $A\preccurlyeq B$ and $C\preccurlyeq D$, cyclically display
their proofs as $\seq{B,A'}$ and $\seq{D,C'}$.
Product gives $\seq{BD,A',C'}$; fold $A',C'$ and rotate to
obtain $AC\preccurlyeq BD$.
The two product bracketings yield actual proofs
\[
 \seq{(AB)C,A',B',C'},\qquad
 \seq{A(BC),A',B',C'}
\]
from three identities.  Folding the negative word using the opposite
bracketing proves both comparisons between $(AB)C$ and $A(BC)$.
For the unit laws, product with $\seq e$ proves
$A\preccurlyeq Ae$ and $A\preccurlyeq eA$.
Conversely insert $d$ in $\seq{A',A}$, and fold $A',d$ or
$d,A'$; this proves $Ae\preccurlyeq A$ and $eA\preccurlyeq A$.

Apply a join injection to the negative formula of an identity to
obtain $A\wedge B\preccurlyeq A$ and $A\wedge B\preccurlyeq B$.
If $C\preccurlyeq A$ and $C\preccurlyeq B$, shared conjunction
with context $C'$ gives $C\preccurlyeq A\wedge B$.
Join injections give $A\preccurlyeq A\vee B$ and
$B\preccurlyeq A\vee B$; its least-upper-bound property follows
by order reversal from the meet property.  These properties also
give monotonicity of both lattice operations.  The operations
therefore respect symmetric equivalence, and the quotient has the
claimed lattice-ordered monoid structure.

By the definitions of the residual formulas and
\eqref{sl:logic:unfold},
$C\preccurlyeq A\backslash B$ means that
$\seq{C',B,A'}$ has a proof.  Rotating and folding gives
$\seq{A'+C',B}$, namely $AC\preccurlyeq B$.
Every step is reversible.  Similarly,
$C\preccurlyeq B/A$ unfolds to $\seq{C',A',B}$ and folds
to $CA\preccurlyeq B$.  Thus both residuation laws hold.
Involution and its order reversal follow from $A''=A$ and
\eqref{sl:logic:order-reversal}; finally
\[
 A\backslash d=(eA)'\equiv A',\qquad
 d/A=(Ae)'\equiv A'.
\]
This verifies the common cyclic dualizer and completes the congruence
assertion.  Totality passes directly to the quotient order.
\end{proof}

\begin{proof}[Proof of Theorem~\ref{sl:logic:theory-chain-completeness}]
Enumerate all pairs of formulas as $(A_n,B_n)$, $n\ge0$.
Starting with $\mathcal T_0=\{\seq A:A\in T\}$, use
Lemma~\ref{sl:logic:orientation} to adjoin one of
$\seq{A_n',B_n}$ and $\seq{B_n',A_n}$ while keeping $H$
underivable; call the result $\mathcal T_{n+1}$.
Put $\mathcal T=\bigcup_n\mathcal T_n$.
A proof is finite, so $\mathcal T\nvdash H$: any contrary
proof would use extra leaves from a single $\mathcal T_n$.
This argument does not require $\mathcal T_0$ to be finite.
Every pair is now comparable.  Lemma~\ref{sl:logic:comparison-algebra}
gives a countable cyclic involutive residuated chain $C$ of formula
classes.

We verify all additional conditions on this quotient.
If $[A][B]\le[e]$, the comparison proof unfolds to
$\seq{A',B',e}$.  WCM gives $\seq{B',A',e}$, which folds
to $[B][A]\le[e]$.
Fix any finite history $w$.  If
$[A][B]\le[c_w]$ and $[U][V]\le[h_w]$, unfold and rotate
their comparison proofs to
\[
 \mathcal T\vdash\seq{c_w,A',B'},\qquad
 \mathcal T\vdash\seq{h_w,U',V'}.
\]
For an arbitrary formula $Z$, use the identity $\seq{Z',Z}$
as the third premise of $R_w$, with
\[
 (\Sigma,\Pi,\Omega,\Lambda,\Gamma,\Delta)
       =(A',B',U',V',Z',Z).
\]
The conclusion and its first-five-formula folding give
\[
 \mathcal T\vdash\seq{B',A',Z',V',U',Z},\qquad
 [B][A][Z][V][U]\le[Z].
\]
Every element has a formula
representative and every closed history formula has its quotient
term value.  This proves \eqref{sl:logic:history-condition} for
all elements and every finite $w$, so $C\in\KR$.

Take the canonical valuation $v(P)=[P]$. Induction on formulas gives
$v(A)=[A]$. If $A\in T$, the assumption leaf $\seq A$ and the
dualizer rule give $\mathcal T\vdash\seq{d,A}$.
Since $d=e'$, the quotient order gives $e_C=[e]\le[A]=v(A)$.
Hence $v\models T$.

Suppose that a component
$\seq{A_1,\ldots,A_k}$ of $H$ were true.  Put
$W=A_1'\cdots A_k'$ with any fixed bracketing, or $W=e$
when $k=0$.  Truth says $[W]\le[d]$, so
$\mathcal T\vdash\seq{W',d}$.  Cut the last $d$ against
$\seq e$ to obtain $\mathcal T\vdash\seq{W'}$.
For $k>0$, repeated unfolding in \eqref{sl:logic:unfold} gives
the original ordered component.  For $k=0$, $W'=d$, and a
second Cut against $\seq e$ gives the empty sequent.
In either case EW adds the remaining components of $H$, contrary
to $\mathcal T\nvdash H$.  Every component is therefore false.
For the empty hypersequent falsity is immediate.
\end{proof}

The orientation construction is an existence argument, not a
decision procedure for choosing the underivability-preserving branch.
The proof uses all external backgrounds and allows arbitrary ordinary
Cut, external contraction, shared conjunction, and finite histories.

\section{A two-ended construction for odd chains}
\label{sl:struct:odd-appendix}

Throughout this section $A$ is an odd cyclic involutive residuated
chain: its dualizer is its identity $e$. We construct a full extension
in which every old point has a new point on each side. The use of
two local identities is essential in the noncommutative case.

The construction treats an element as an arrow between its right and
left local identities. We first show how to lift arrows to common
identities and compare them there. We then split the local layers while
retaining those transports. This produces three ordered copies of the
old chain, with the middle copy preserving all original operations.

\begin{lemma}[Local identities, matching products and tails]
\label{sl:struct:odd-local}
Put
\[
 \ell(x)=x/x=(xx')',\qquad \rho(x)=x\backslash x=(x'x)',
 \qquad E=\{u\ge e:u^2=u\},
\]
and $H_{v,u}=\{x:\ell(x)=v,\ \rho(x)=u\}$.
Then $\ell(x),\rho(x)\in E$, and
\begin{gather}
 \ell(x)x=x=x\rho(x),\quad xx'=\ell(x)',\quad
 x'x=\rho(x)',                                      \label{sl:struct:local-eqs}\\
 \ell(x')=\rho(x),\quad \rho(x')=\ell(x),\quad
 x\rho(x)'=\ell(x)'x.                               \label{sl:struct:intertwine}
\end{gather}
For $u<v$ in $E$,
\begin{equation}\label{sl:struct:absorption}
 uv=vu=v,\quad uu'=u'u=u',\quad vu'=u'v=v.
\end{equation}
Products at matching identities satisfy
$H_{v,u}H_{u,w}\subseteq H_{v,w}$.
For $x\in H_{v,u}$ the maps
\[
 \alpha_x:[u,\infty)\cap E\longrightarrow[v,\infty)\cap E,
 \qquad \alpha_x(u)=v,\quad \alpha_x(z)=xzx'\ (z>u)
\]
and $\alpha_{x'}$ are inverse order isomorphisms. Furthermore
\begin{equation}\label{sl:struct:tail-eqs}
 \alpha_x(z)x=xz,\qquad
 \ell(xz)=\alpha_x(z),\quad \rho(xz)=z.
\end{equation}
For $z>u$, the element $xz$ has a matching inverse $zx'$:
their products are $\alpha_x(z)$ and $z$.
\end{lemma}
\begin{proof}
Residual evaluation and $e\le\ell(x),\rho(x)$ give the local
identity equations. Since $x\rho(x)^2=x$, residuation gives
$\rho(x)^2\le\rho(x)$; the reverse inequality follows from
$\rho(x)\ge e$. The left version is identical. The remaining
equations follow from the involutive formulas and associativity.
For $a\ge e$, the identities imply $\ell(a),\rho(a)\le a$.
Thus a positive idempotent $u$ has $\ell(u)=\rho(u)=u$ and
$uu'=u'u=u'$. Two positive idempotents multiply to their maximum
by monotonicity. For $u<v$, put $b=vu'$. Cyclic duality gives
$b>e$, while $b\le v$ and $vb=b$. Hence
$v\le\ell(b)\le b\le v$, so $vu'=v$. The proof for $u'v$
uses right identities.

For matching $x\in H_{v,u}$ and $y\in H_{u,w}$, put
$a=\ell(xy)\ge v$. Multiplying $axy=xy$ by $y'$ gives
$axu'=xu'$, hence $av'x=v'x$. If $a>v$, absorption gives
$av'=a$ and consequently $ax=v'x\le x$. Thus $ax=x$ and
$a\le\ell(x)=v$, a contradiction. Hence $\ell(xy)=v$.
The symmetric calculation gives $\rho(xy)=w$.

For $z>u$, cyclic duality gives $a=xzx'>e$, because
$xzx'\le e$ would be equivalent to $xz\le x$, or $z\le u$.
Absorption yields
$a^2=xzu'zx'=a$ and $ax=xzu'=xz$. Since $xz>x$, this forces
$a>v$. Also
\[
 (xz)(zx')=a,\qquad (zx')(xz)=z.
\]
If $b=\ell(xz)$, then $b\ge a$ and $b(xz)=xz$; multiplying
by $zx'$ gives $ba=a$. Positive idempotents multiply to their
maximum, so $b=a$. Similarly $\rho(xz)=z$. Finally
$x'(xzx')x=u'zu'=z$. Interchanging $x,x'$ proves the inverse
assertion. Monotonicity and the endpoint values now give the stated
tail order isomorphisms.
\end{proof}

We regard $x\in H_{v,u}$ as an arrow $v\leftarrow u$.
An arrow is an \emph{isomorphism} when it has an inverse for
matching multiplication with positive identities $v,u$; this
inverse is denoted $x^\dagger$. It need not equal $x'$.

\begin{lemma}[The two kinds of local layer]
\label{sl:struct:odd-layers}
For $u\in E$ write $f_u=u'$. The endofibre $H_{u,u}$ is
closed under product, prime and both residuals, with identity $u$,
dualizer $f_u$, and $xx'=x'x=f_u$ for each of its elements.
The element $f_u$ is central in this endofibre. If $u>e$, then
$f_u<u$ is a cover within it. If $f_u^2<f_u$, the endofibre is
a group. If $f_u^2=f_u$, then
\begin{equation}\label{sl:struct:identity-fibre}
 \{x\in H_{u,u}:f_ux=f_u\}=\{f_u,u\}.
\end{equation}
For $u=e$ the endofibre is the global group of units.
\end{lemma}
\begin{proof}
Matching products and prime preserve the endofibre, so the residual
formulas do too. Its local duality follows from
$x\backslash f_u=(ux)'=x'$. Centrality is
$f_ux=(xx')x=x(x'x)=xf_u$. If $f_u<z<u$ in the endofibre,
then $z\ge e$ would imply $u=\ell(z)\le z$. If $z\le e$,
then $z'\ge e$ and $u=\ell(z')\le z'$, giving $z\le f_u$.
Both alternatives are impossible.

Suppose $f_u^2<f_u$. The element $h=(f_u^2)'$ belongs to the
endofibre and satisfies $h>u$ and $f_uh=hf_u\le u$.
If $f_uh\le f_u$, then
$h\le f_u\backslash f_u=(uf_u)'=u$, a contradiction.
The cover just proved gives $f_uh=hf_u=u$. Thus $f_u$ is a
central unit; $f_u^{-1}x'$ is a two-sided inverse of each $x$.
If instead $f_u^2=f_u$, an element with $f_ux=f_u$ satisfies
$x\le f_u\backslash f_u=u$. It also satisfies $x\ge f_u$,
because $f_u\le e$ and $f_ux\le x$. The cover proves
\eqref{sl:struct:identity-fibre}. At $u=e$, the equations
$xx'=x'x=e$ identify exactly the global units.
\end{proof}

Partition $E$ into $\{e\}$ and
\[
 E_I=\{u>e:f_u^2=f_u\},\qquad
 E_J=\{u>e:f_u^2<f_u\}.
\]
There are no other cases, since $f_u\le e$.

\begin{lemma}[Types and quotient arrows]
\label{sl:struct:odd-quotient}
A nonempty $H_{v,u}$ has either $u,v\in E_J$ or
$u,v\in\{e\}\cup E_I$. Every arrow of the first kind is
an isomorphism. No isomorphism connects $e$ to another object.
On the second kind put
\[
 \sigma(x)=f_vx=xf_u.
\]
The images, with matching multiplication, form a groupoid $P$,
whose identity at $u$ is $f_u$ and whose inverse is
$\sigma(x)^{-1}=\sigma(x')$. The image of the isomorphisms is
a subgroupoid. Every fibre of $\sigma$ has one or two points.
It has two points exactly over an isomorphism image with endpoints
in $E_I$; these are the lower point $\sigma(x)$ and the unique
upper isomorphism $x$. If an endpoint is $e$, the fibre is single.
\end{lemma}
\begin{proof}
If $f_u^2=f_u$, the intertwining equation gives
$f_v^2x=f_vx$ and hence $f_v^3=f_v^2$ after multiplying by
$x'$. If $v\in E_J$, Lemma~\ref{sl:struct:odd-layers} allows
cancellation and yields $f_v=v$, a contradiction. Interchanging
the endpoints proves type agreement. When both endpoints belong
to $E_J$, the arrow
$f_u^{-1}x'=x'f_v^{-1}$ is a matching inverse of $x$.
If $gx=e$ and $xg=v$, unit cyclicity gives $v\le e$, so $v=e$.

On the remaining objects all $f_u$ are idempotent. Matching
composition gives
\[
 \sigma(x)\sigma(y)=\sigma(xy),\quad
 \sigma(x)\sigma(x')=f_v,\quad
 \sigma(x')\sigma(x)=f_u.
\]
Also $\sigma^2=\sigma$. These equations prove the groupoid
assertion. Isomorphisms compose and invert, proving the subgroupoid
assertion. Suppose $x<y$ have equal image in $H_{v,u}$, with
$u,v\in E_I$. The products $yx'$ and $x'y$ project to $f_v,f_u$.
Moreover $yx'>f_v$, for otherwise
$y\le f_v/x'=(x'v)'=x$; likewise $x'y>f_u$.
Equation~\eqref{sl:struct:identity-fibre} therefore gives
$yx'=v$ and $x'y=u$. Thus $y$ is an isomorphism and
$x=xx'y=f_vy=\sigma(y)$.
Conversely an isomorphism $y$ satisfies $\sigma(y)<y$:
equality, multiplied by $y^\dagger$, would give $f_v=v$.
The preceding argument proves uniqueness of the upper point.
If an endpoint is $e$, the identity $f_e=e$ and intertwining
make $\sigma$ the identity on that hom set.
\end{proof}

Write $C_u=\{x:\rho(x)=u\}$ for an entire right column;
its arrows can have different left endpoints. All column orders
below are the inherited chain orders, not independently chosen
orders on hom sets.

\begin{lemma}[Column order and lifts]
\label{sl:struct:odd-columns}
On a column of type $\{e\}\cup E_I$, the map $\sigma$ is
monotone, is strictly increasing on isomorphisms, and has adjacent
split pairs in the entire column. For $z>\rho(x)$,
\begin{equation}\label{sl:struct:quotient-lift}
 \sigma(x)z=xz
\end{equation}
when $\sigma(x)$ is defined. If $\rho(x)<\rho(y)$, then
\begin{equation}\label{sl:struct:directed-order}
 x<y\quad\Longleftrightarrow\quad x\rho(y)\le y.
\end{equation}
Matching right multiplication by an isomorphism is an order
isomorphism between the entire corresponding columns. The same
holds for matching multiplication in the quotient groupoid $P$.
For every $u>e$ there is no point of $C_u$ strictly between
$f_u$ and $u$. For an isomorphism $a\in C_u$, the point $af_u$
is its immediate predecessor in $C_u$.
\end{lemma}
\begin{proof}
On $C_u$, the projection is the single monotone multiplication
$x\mapsto xf_u$. Equal projected values are equal elements of
the old chain and therefore have the same left endpoint; the
preceding fibre classification applies. It proves strictness on
isomorphisms and adjacency of every split pair in the whole column.
Equation~\eqref{sl:struct:quotient-lift} is absorption $f_uz=z$.
For \eqref{sl:struct:directed-order}, put $v=\rho(y)>\rho(x)$.
We have $x<xv$, because equality would imply $v\le\rho(x)$.
Thus $xv\le y$ implies $x<y$. Conversely $x<y$ gives
$xv\le yv=y$.

Matching multiplication by $a$ and by $a^\dagger$ are monotone
mutual inverses on the whole columns, by the local identity
equations. For a quotient arrow $b$, the same reasoning uses
$bb^{-1}=f_{\ell(b)}$ and $b^{-1}b=f_{\rho(b)}$; every
quotient arrow is fixed by its corresponding $f$.
The proof of the cover $f_u<u$ in
Lemma~\ref{sl:struct:odd-layers} uses only the right identity:
if $f_u<z<u$, either $e\le z$ gives $u=\rho(z)\le z$, or
$z\le e$ gives $u=\ell(z')\le z'$ and $z\le f_u$.
Thus it proves the cover in $C_u$. Translate this cover by the
appropriate matching isomorphism to obtain $af_u<a$.
\end{proof}

For later calculations put $R_zx=xz$ when $z\ge\rho(x)$.
Its target is $\alpha_x(z)$ by
Lemma~\ref{sl:struct:odd-local}. Its left-lift counterpart is
$L_zx=zx=R_{\alpha_x^{-1}(z)}x$. For arbitrary $x,y$, putting
$k=\max(\rho(x),\ell(y))$ gives
\begin{equation}\label{sl:struct:old-total}
 xy=(R_kx)(R_{\alpha_y^{-1}(k)}y),\quad
 \ell(xy)=\alpha_x(k),\quad
 \rho(xy)=\alpha_y^{-1}(k).
\end{equation}
Indeed $xy=(xk)(ky)$ because $k$ is one of the matching old
identities, and matching composition then gives the endpoints.
On its domain tail,
\begin{equation}\label{sl:struct:old-coherence}
 \alpha_{xy}=\alpha_x\alpha_y,\qquad
 R_z(xy)=(R_{\alpha_y(z)}x)(R_zy)
\end{equation}
whenever $\rho(x)=\ell(y)$ and $z\ge\rho(y)$.
At the initial object these are the local identity equations.
Above it, both factors are strict lifts, and their targets are
computed by Lemma~\ref{sl:struct:odd-local}; associativity of
the old product proves the formulas. Also $R_wR_zx=R_wx$ for
$w\ge z$, because positive idempotents multiply to their maximum.

\subsection{The new objects and arrows}

Form a new ordered object set $E^B$ by inserting $u^-$ immediately
below each $u>e$, and $u^+$ immediately above each
$u\in\{e\}\cup E_I$. Different old bases retain their order.
Every new object has type $I$. Retain all old arrows and add:
\begin{enumerate}
\item for every old isomorphism $a:v\leftarrow u$ with $u,v>e$,
two arrows $a^-_0,a^-_1:v^-\leftarrow u^-$;
\item for every quotient arrow $a:v\leftarrow u$ of $P$, two
arrows $a^+_0,a^+_1:v^+\leftarrow u^+$.
\end{enumerate}
Here $0<1$ gives the order of the copies of one arrow. Each new
column is ordered by its entire underlying old or quotient column,
then by this binary mark. At matching new objects define
\begin{equation}\label{sl:struct:new-matched}
 a^\varepsilon_i b^\varepsilon_j
       =(ab)^\varepsilon_{\min(i,j)},\qquad
 (a^\varepsilon_i)'=(a^{-1})^\varepsilon_{1-i}
 \quad(\varepsilon\in\{-,+\}).
\end{equation}
The inverse is the matching groupoid inverse: for minus arrows,
$a^{-1}=a^\dagger$. The upper identity copy is the new object's
identity. Old matched products and prime remain unchanged.

Represent an old arrow by itself, a minus arrow by its old
isomorphism, and a plus arrow by its old quotient representative;
write $\bar x$ for this representative, forgetting its mark.
The lift to the arrow's own source is itself. For a strict lift
to an object $s$ with old base $z$ define
\begin{equation}\label{sl:struct:new-lifts}
 R^B_sx=
 \begin{cases}
 \bar xz,&s=z,\\
 (\bar xz)^-_1,&s=z^-,\\
 (\sigma(\bar xz))^+_1,&s=z^+.
 \end{cases}
\end{equation}
If $x$ is old and $s=\rho(x)^+$, the last expression is
$(\sigma(x))^+_1$. Every other strict same-base lift begins at
a minus object, whose representative is already an isomorphism.
At a higher old base, $\bar xz$ is a strict old lift and hence
an isomorphism. Thus every expression is defined in its indicated
groupoid. Its target defines a tail map $\alpha^B_x$. At a
higher base $z$ this is $\alpha_{\bar x}(z)$ with the same
available mark $-,0,+$. At the source it is the target of $x$;
same-base transitions retain their marks and change only the base.

\begin{lemma}[Coherent transports]
\label{sl:struct:clone-coherence}
Each $\alpha^B_x$ is an order isomorphism of the upper tails of
its source and target. The maps $R^B_s$ are monotone on columns
and satisfy
\begin{gather}
 R^B_tR^B_sx=R^B_tx\quad(t\ge s\ge\rho(x)),
                                                        \label{sl:struct:clone-coh1}\\
 \alpha^B_{R^B_sx}=\alpha^B_x\big|_{[s,\infty)},
                                                        \label{sl:struct:clone-coh2}\\
 R^B_s(xy)=(R^B_{\alpha^B_y(s)}x)(R^B_sy),\quad
 \alpha^B_{xy}=\alpha^B_x\alpha^B_y
                                                        \label{sl:struct:clone-coh3}
\end{gather}
for matching $x,y$ and $s\ge\rho(y)$. Every strict lift is an
isomorphism, and
\begin{equation}\label{sl:struct:clone-inverse}
 (R^B_sx)^\dagger=R^B_{\alpha^B_x(s)}(x')
 \quad(s>\rho(x)).
\end{equation}
\end{lemma}
\begin{proof}
At a higher base, old tail isomorphisms preserve the available
types: the corresponding old strict lift is an isomorphism and
therefore connects two equal types. At an initial old object, the
same-base tail is $0<+$ for type $\{e\}\cup E_I$, and only $0$
for type $E_J$. At an initial minus object the remaining same-base
tail is $-<0<+$ in type $I$ and $-<0$ in type $J$. A plus object
has no higher same-base object. These lists coincide at source
and target, so the tail maps are order isomorphisms.

Higher-base lifts are old monotone multiplications, optionally
followed by $\sigma$ and insertion as an upper copy. Same-base
lifts are inclusion of the isomorphism column, the monotone
quotient map, or their composite. Forgetting a binary mark is
monotone. This proves column monotonicity.

Coherence at different bases follows from old coherence and
$\sigma(x)z=xz$ for a strict higher-base lift. At one base the
only possible composite is inclusion followed by $\sigma$ along
$-<0<+$, which is exactly the prescribed $-<+$ map. This proves
\eqref{sl:struct:clone-coh1}; taking targets proves
\eqref{sl:struct:clone-coh2}. In
\eqref{sl:struct:clone-coh3}, equality at the source is matched
composition. Otherwise both right-hand factors undergo strict
lifts. At higher bases the equality is
\eqref{sl:struct:old-coherence}; at the same base it follows
because inclusion, $\sigma$, and mark-forgetting each preserve
matching multiplication. Both output marks are $1$, whose product
is again $1$. Taking targets proves the transport identity.
Finally \eqref{sl:struct:clone-inverse} at a higher base is the
old strict-lift inverse formula. For minus-to-old it is the
definition of $\dagger$; for old-to-plus it is
$\sigma(x')=\sigma(x)^{-1}$; minus-to-plus is their composite.
These exhaust the strict same-base cases.
\end{proof}

For arbitrary arrows define, with $k=\max(\rho(x),\ell(y))$,
\begin{equation}\label{sl:struct:clone-product}
 x\odot y=(R^B_kx)
          (R^B_{(\alpha^B_y)^{-1}(k)}y).
\end{equation}
The right side is a matching product. Order equal right columns
as specified above, and different columns by
\begin{equation}\label{sl:struct:clone-order}
 x<y\quad\Longleftrightarrow\quad
 \begin{cases}
 R^B_{\rho(y)}x\le y,&\rho(x)<\rho(y),\\
 x<R^B_{\rho(x)}y,&\rho(x)>\rho(y).
 \end{cases}
\end{equation}

\begin{lemma}[Associativity and the precise new order]
\label{sl:struct:clone-order-proof}
The product \eqref{sl:struct:clone-product} is associative and has
identity $e$. The relation \eqref{sl:struct:clone-order} is a total
order extending that of $A$. As an ordered set with involution,
the extension is
\begin{equation}\label{sl:struct:three-copies}
 B(A)\cong A\mathbin{\overrightarrow\times}\{-1,0,1\},
 \qquad (x,i)'=(x',-i),
\end{equation}
where $(x,0)$ is the old point.
\end{lemma}
\begin{proof}
The old identity lifts to each new upper identity, proving the
unit laws. For associativity write the arrows as
$x:a\leftarrow b$, $y:c\leftarrow d$, $z:f\leftarrow g$ and put
\[
 \beta_0=(\alpha^B_y)^{-1}(\max(b,c)),\quad
 \beta=\max(\beta_0,f),\quad
 \gamma=(\alpha^B_z)^{-1}(\beta).
\]
The coherence identities reduce either bracketing to the matching
product of
\[
 R^B_{\alpha^B_y(\beta)}x,\qquad
 R^B_\beta y,\qquad R^B_\gamma z.
\]
For the second bracketing use that a tail order isomorphism
preserves maxima. Matching products are associative, including
their binary marks, so the total product is associative.

Trichotomy of \eqref{sl:struct:clone-order} follows from column
trichotomy. A strict inequality between lifts to a common source
forces the corresponding original strict inequality: the opposite
weak inequality would lift monotonically in the opposite direction.
For transitivity of $x<y<z$, lift to the largest of their sources.
If that largest source is the first or last, the needed weak or
strict comparison follows directly by composing the two column
comparisons. If only the middle source is largest, the comparisons
there are
$R x\le y<R z$, so $R x<R z$. If the first and last sources
agree, this forces their strict comparison. If the first is larger,
it forces $x<R_{\rho(x)}z$; if the last is larger, it forces
$R_{\rho(z)}x\le z$. Otherwise monotonicity would contradict
the strict comparison at the middle source. These are exactly the
two clauses of \eqref{sl:struct:clone-order}. Thus the relation
is transitive. Its restriction to old points is
\eqref{sl:struct:directed-order}.

We give the positions of all new points, which also proves
antitonicity rather than assuming it. For an old isomorphism $a$
with source $u>e$, the point $a^-_1$ is immediately below $a$
relative to the old order, and $a^-_0$ is immediately above
$af_u$. For a quotient representative $a$, the point $a^+_0$
is immediately below the bottom of its old $\sigma$ fibre,
whereas $a^+_1$ is immediately above its top. Here ``relative to
the old order'' means that no old point lies in between.
To verify these statements against an old point $y$, separate
the cases $\rho(y)<u$, $\rho(y)=u$ and $\rho(y)>u$.
In the first case \eqref{sl:struct:clone-order} tests its lifted
old coordinate against a lower or upper copy; the lower-copy test
is strict, the upper-copy test weak. In the equal case use the
whole-column adjacent pairs of
Lemma~\ref{sl:struct:odd-columns}. In the last case both members
of an old split pair have the same strict lift, and so do $a$ and
$af_u$ by absorption. Substitution in
\eqref{sl:struct:clone-order} gives the stated comparisons in
each case.

Consequently every old point $x$ has one designated new lower
point $b(x)$ and upper point $a(x)$, as follows:
\[
\begin{array}{c|cc}
\text{kind of }x&b(x)&a(x)\\ \hline
E_J\text{ arrow}&x^-_1&(xf_{\rho(x)}^{-1})^-_0\\
\text{upper split point}&x^-_1&\sigma(x)^+_1\\
\text{lower split point, upper mate }h&x^+_0&h^-_0\\
\text{singleton quotient point}&x^+_0&x^+_1
\end{array}
\]
For the first row, $f_u^{-1}z=z$ at every $z>u$: multiply
$f_uz=z$ by the local inverse and use $uz=z$. Thus the previous
strict-lift comparison applies to its upper point too. Every new
point occurs exactly once in this table.

It remains to check $a(x)<b(y)$ for old $x<y$. An intermediate
old point already proves it, so assume $x<y$ is a cover. At equal
sources, either quotient coordinates differ, or the points form
a split pair, or they form a type-$J$ predecessor pair. The
underlying column order proves the first case, and $0<1$ proves
the latter two. If $\rho(x)<\rho(y)$, the cover forces
$x\rho(y)=y$, since otherwise this old strict lift lies between
them. Thus $y$ is an isomorphism, and the lift of $a(x)$ to
$\rho(y)^-$ is $y^-_1=b(y)$. The first clause of
\eqref{sl:struct:clone-order} gives the strict comparison.
If $\rho(x)>\rho(y)$, put $h=y\rho(x)$. The immediate
predecessor in the column gives
$x\le hf_{\rho(x)}<y$, so the cover forces $x=hf_{\rho(x)}$.
Now $a(x)=h^-_0$, while the lift of $b(y)$ to $\rho(x)^-$ is
$h^-_1$. The second clause gives the result. This proves the
lexicographic description in \eqref{sl:struct:three-copies}.

Prime exchanges $b(x)$ with $a(x')$. For type $J$ use
$x'=x^\dagger f_{\ell(x)}$. For a split pair it exchanges
the lower point with the upper point over the inverse coordinate:
the upper point $h$ has $h'=h^\dagger f_{\ell(h)}$, and
$(\sigma(h))'=h^\dagger$, by the residual formulas for a
matching isomorphism. On singleton quotient points it is
groupoid inversion. Together with
\eqref{sl:struct:new-matched}, these formulas prove
$(x,i)'=(x',-i)$, hence antitonicity.
\end{proof}

\begin{theorem}[The odd extension]
\label{sl:struct:odd-clone}
With the operations above, $B(A)$ is an odd cyclic involutive
residuated chain. The inclusion of $A$ is full, and every old
point has a new point on each side; every strict old pair receives
an intermediate point. Every old central unit remains central.
\end{theorem}
\begin{proof}
We prove cyclic duality. In a column with source $s>e$, comparison
with $e$ is strict comparison with its column identity $s$,
because $R^B_se=s$. In column $e$ it is the usual comparison
with $e$. For two arrows $a,b$ with the same source,
\begin{equation}\label{sl:struct:column-duality}
 a\le b\quad\Longleftrightarrow\quad a\odot b'\le e.
\end{equation}
In an old column this is old duality and old matching product.
In a new column, matching right translation by the inverse of
the underlying coordinate of $b$ is a strict order isomorphism
onto its target column. At unequal coordinates it turns the
comparison into strict comparison with the identity coordinate.
At equal coordinates the product mark is
$\min(i,1-j)$, which is $0$ exactly when $i\le j$.
This proves \eqref{sl:struct:column-duality}.

If $\rho(x)<\rho(y)$, then
$x\odot y'=(R^B_{\rho(y)}x)y'$, so
\eqref{sl:struct:clone-order} and
\eqref{sl:struct:column-duality} give
$x<y\Longleftrightarrow x\odot y'\le e$.
If $\rho(x)>\rho(y)$, set $b=R^B_{\rho(x)}y$.
It is a strict-lift isomorphism with target
$s=\alpha^B_y(\rho(x))>\ell(y)\ge e$.
Equation~\eqref{sl:struct:clone-inverse} gives
$x\odot y'=xb^\dagger$.
Right multiplication by $b^\dagger$ is a strict order
isomorphism on the whole relevant column: for an old column
this is Lemma~\ref{sl:struct:odd-columns}; for a new one the
underlying groupoid translation is an order isomorphism and the
upper inverse copy preserves each binary mark. Hence
\[
 x<b\ \Longleftrightarrow\ xb^\dagger<s
 \ \Longleftrightarrow\ x\odot y'\le e.
\]
Together with the equal-source case this proves
$x\le y\Longleftrightarrow x\odot y'\le e$ universally.
Antitonicity and involutivity now give
\[
 x\odot y\le e\quad\Longleftrightarrow\quad y\le x'
 \quad\Longleftrightarrow\quad y\odot x\le e.
\]
Associativity yields both residual adjunctions; for example
\[
 x\odot y\le z
 \ \Longleftrightarrow\ x\odot y\odot z'\le e
 \ \Longleftrightarrow\ y\le(z'\odot x)'.
\]
These adjunctions also imply monotonicity of multiplication.
Thus this is an odd cyclic involutive residuated chain and the
old inclusion is full. The assigned local endpoints are the
actual residual endpoints: matching multiplication gives
$xx'=\ell(x)'$ and $x'x=\rho(x)'$, including the lower
identity copy at a new object, and the residual formulas give
$x/x=\ell(x)$ and $x\backslash x=\rho(x)$.
Its type is $I$, because that lower identity is idempotent.
The order description proves the insertion assertions.

Finally, every old global unit $c$ has both endpoints $e$.
Conjugation by $c$ preserves old endpoints, types, quotient
projections and matching inverses. Extending it to the new
arrows by
\[
 a^\varepsilon_i\longmapsto
       (cac^{-1})^\varepsilon_i
\]
preserves all the defining formulas. This extended map equals
actual conjugation by $c$ in $B(A)$: multiplication by $c$ on
either side lifts it strictly to the relevant new object, using
an upper mark, which leaves $i$ unchanged under the minimum
rule. In the minus family its lifted representatives are $cv$
and $uc^{-1}$, giving $ca$ and $cac^{-1}$; in the plus family
they are $cf_v$ and $f_uc^{-1}$, and the quotient identity
equations give the same products. On old arrows this is old
conjugation. If $c$ was central, its old conjugation and therefore
its extended conjugation are the identity. Thus it stays central.
\end{proof}

\begin{corollary}[Odd densification]
\label{sl:struct:odd-dense}
Every countable odd cyclic involutive residuated chain has a
countable full dense extension without endpoints. All its old
central units remain central with their old inverses.
\end{corollary}
\begin{proof}
Iterate Theorem~\ref{sl:struct:odd-clone} countably many times.
Each step adds two points per old point, so is countable. A pair
or point in the union occurs at a finite stage and receives the
required intermediate or outer points at the next. Equations,
operations and residuals are preserved at every stage. Centrality
of each original central unit holds on every stage and hence on
the union.
\end{proof}

\begin{corollary}[Relative odd densification]
\label{sl:struct:relative-odd}
Suppose $A$ is countable and directly finite, meaning
$xy=e\Rightarrow yx=e$. There is a countable full odd extension
$A\hookrightarrow D$ with exactly the same unit group $G$ and
the following properties: $D$ is directly finite; every strict
pair with a nonunit endpoint has an intermediate point; every old
cover between two units remains a cover; and all old central
units remain central. If $G$ is nontrivial, $D$ has no endpoints.
\end{corollary}
\begin{proof}
In a directly finite odd chain there are no arrows between $e$
and another object. Indeed $\rho(x)=e$ gives $x'x=e$, and
direct finiteness gives $xx'=e$, so $\ell(x)=e$; the converse
is symmetric. Delete the newly introduced object $e^+$ and all
its arrows from $B(A)$. There are no arrows between this deleted
component and any other component, because the corresponding
old quotient arrows would have an endpoint $e$ and another
endpoint different from $e$.

The retained set is product- and prime-closed. More explicitly,
all its nonunits have both endpoints strictly above $e$, while
old units have endpoints $e$. Formula~\eqref{sl:struct:clone-product}
shows that a product involving a nonunit still has both endpoints
above $e$; a product of old units is its old product. Thus no
retained product falls in the deleted component. Prime exchanges
the two endpoints, so is retained as well. The resulting full
subalgebra has exactly the old global units. A product equal to
$e$ must have two old unit factors, which proves direct finiteness.
The deleted new points are precisely the two neighbours assigned
to each old unit in \eqref{sl:struct:three-copies}. Every nonunit
retains its two neighbours, and a cover between old units receives
no new point. Centrality follows by restriction from the full clone.
Iterate this restricted construction. Every strict pair with a
nonunit endpoint is separated at the next stage, and old unit
covers persist at every stage. A nonunit point gets points on
both sides. If $G$ is nontrivial, every unit has a smaller and a
larger old unit, by multiplication by a positive unit or its
inverse. Hence the union has no endpoints.
\end{proof}

\begin{corollary}[Relative densification with a central unit dualizer]
\label{sl:struct:relative-unit}
Let $J$ be a countable directly finite cyclic involutive residuated
chain with central invertible dualizer $q$. It has a countable
full extension $T$ with the same group of units, retaining every
old central unit and every old cover between units, in which every
strict pair with a nonunit endpoint has an intermediate point.
If $J$ has a central unit $n<r$ that immediately precedes its
identity $r$, the only covers of $T$ are $gn<g$ for $g\in U(J)$,
and $T$ has no endpoints.
\end{corollary}
\begin{proof}
Normalize to the odd prime $q^{-1}x^*$ as in
Lemma~\ref{sl:struct:unit-histories}, apply
Corollary~\ref{sl:struct:relative-odd}, and restore prime by
$x^*=q x^\circ$. The old $q$ remains central with its inverse,
so these operations are well defined and preserve the full old
signature. Direct finiteness and units concern the unchanged
ordered monoid. If $n<r$ is the indicated cover, multiplication
by a unit gives a cover $gn<g$ in $J$. Conversely every cover
between units must be of this form, because translating its
upper endpoint to $r$ gives the unique predecessor $n$.
Every remaining strict pair is separated either by a nonunit
already between its endpoints or by the relative construction.
The nontrivial unit group supplies the no-endpoints assertion.
\end{proof}

\section{Technical proofs for the lower-core construction}
\label{sl:lower:appendix}

\subsection{Old actions and their actual adjoints}
\label{sl:lower:appendix-actions}

\begin{proof}[Proof of Lemma~\ref{sl:lower:actions}]
Let $Rf(b)=b$ for an admissible return. If $z<f(b)$, then
$fR(z)\le z<f(b)$, whence $R(z)<b$ and $R(z)\le a$. Thus
$z\le R^*(a)$. On the other hand $R^*(a)<f(b)$, since
$Rf(b)=b>a$. Therefore $c=R^*(a)$ is the immediate predecessor
of $f(b)$, independently of the return chosen. Also $f(a)\le c$
because $Rf(a)\le a$. These inequalities give monotonicity across
the inserted gap. For distinct upper anchors $b<t$, the only possible
failure would be an active later gap with $f(b)=f(t)$. Its return
would give $Rf(b)=t>b$, contrary to $Rf\le\mathrm{id}$. This
proves monotonicity on the whole carrier.

For coefficient monotonicity, suppose $u\le\widetilde u$ and
$v\le\widetilde v$. The old upper values are ordered. A failure at
a new input could therefore occur only if they were equal, the larger
context were active, and the smaller inactive. But a return $(s,t)$
for the larger context is admissible for the smaller, and
\[
 b=s\widetilde u b\widetilde v t=su bvt\le b
\]
when their upper values agree. Thus the smaller context is active
as well, excluding that failure.

For composition write $f=(u,v)$ and $h=(w,z)$, so
$hf=(wu,vz)$. If $(s,t)$ returns $hf$ at $b$, then
$(sw,zt)$ returns $f$ at $b$, and $(us,tv)$ returns $h$ at
$f(b)$. Their contractive inequalities are cyclic rotations at $e$
of $swu\le e$ and $vzt\le e$. Their return equations follow
from $swubvzt=b$, and from this equality multiplied on the left
by $u$ and on the right by $v$. Conversely, if $(s,t)$ returns
$f$ and $(S,T)$ returns $h$ at its image, then $(sS,Tt)$
returns $hf$: for example $sSwu\le su\le e$ and
$vzTt\le vt\le e$, and the other two inequalities follow by
$e$-cyclicity. The endpoint equations compose. Thus $hf$ is
active exactly when both successive actions are active. In that
case both evaluations give the unique new point below $hf(b)$;
in every other case both give old $hf(b)$. This proves composition.
The identity label is active and fixes every inserted point.

At an old target $y$, a new input is eligible for $fX\le y$
exactly when its upper anchor $b$ satisfies $f(b)\le y$.
The greatest eligible input is consequently old $f^*(y)$. At a
new target $Y=X_{c,h}$ put $a=f^*(c)$. An eligible new input
above $a$ must map actively to $Y$. Its lower endpoint is eligible
for old $c$, and its upper endpoint is not, so that lower endpoint
must be $a$. There is at most one such gap. Its midpoint, if
eligible, is the maximum; otherwise $a$ is the maximum. These are
actual right adjoints on every target.

It remains to identify them with prime. On the old chain the dual
context $f^\dagger={} '\,f^*\,{}'$ is $(v,u)$. If $f$ actively
transports $(a,b)$ to $(c,h)$ with return $R$, then
$f^*(c)=a$ and $R^*(a)=c$. Indeed $f(a)\le c<f(b)$ and
adjunction determine the first equality, while the second was
proved above. Hence $f^\dagger$ sends $c'$ to $a'$, and the
dual return $R^\dagger=(t,s)$ sends $a'$ to $c'$; both composites
are contractive. It actively transports $(h',c')$ to $(b',a')$.
Dualizing again gives the converse. Thus the adjoint of $f$ at a
new target $Y$ is new exactly when $f^\dagger$ is active at $Y'$,
and in that case the two outputs are paired by prime. In the
remaining case it is old $f^*(c)$, whose prime is the inactive
value $f^\dagger(c')$. Therefore
$f^*={} '\,f^\dagger\,{}'$ on the entire carrier. The one-sided
cases give \eqref{sl:lower:action-laws}.
\end{proof}

\subsection{Attained coefficient comparisons}
\label{sl:lower:appendix-comparison}

\begin{proof}[Proof of Lemma~\ref{sl:lower:comparison}]
For any old $v$,
\[
 vX_J\le X_J\quad\Longleftrightarrow\quad v\le s_J:
 \tag{L9}\label{sl:lower:self-cut}
\]
if $v\le e$ the inequality is automatic, while for $v\ge e$
it is equivalent to equality. The action law shows $s_J^2$ fixes
$X_J$; since $s_J^2\ge s_J\ge e$, maximality gives idempotence.

Suppose $uX_J=X_K$ and put $v=u^-$. With $h=ub$ the activity
criterion gives $vh=b$, and $uvh=h$. Thus $u$ is an admissible
return for $v$ at $K$, and $vX_K=X_J$. Hence, for any old $w$,
\[
 wX_J\le X_K\iff(wv)X_K\le X_K
 \iff wv\le s_K\iff w\le s_K/v.
\]
Also $vZ\le X_J$ if and only if $Z\le X_K$: the reverse
implication is monotonicity, and if $Z>X_K$ then $Z\ge h$,
so $vZ\ge vh=b>X_J$. It follows that the same cut is
$\{w:vw\le s_J\}=\mathord\downarrow(v\backslash s_J)$.
Its maximum is at least $u$ and therefore transports $X_J$ to
$X_K$. This proves the fourth line.

For old $y$, the anchor formula gives $wX\le y$ exactly when
$w r(X)\le y$. For old $x$ and new $X_K$, the product $wx$
is old, so $wx\le X_K$ exactly when $wx\le c$. If both
arguments are new and there is no transport, $wX_J\le X_K$
means $wX_J\le c$, equivalently $wb\le c$. These observations
prove the other three lines.
\end{proof}

\subsection{Cyclic selection and old-factor equations}
\label{sl:lower:appendix-selection}

\begin{proof}[Proof of Lemma~\ref{sl:lower:selection}]
The old target $r(X)r(Y)$ is eligible for \eqref{sl:lower:k}, by
old cyclic adjunction and monotonicity of the old actions. Thus
$k(X,Y)\le r(X)r(Y)$, and
\[
 k(X,Y)\le z\iff z'X\le Y'\qquad(z\in A).
 \tag{L16}\label{sl:lower:k-cut}
\]
If either upper-action bound $r(X)Y$ or $Xr(Y)$ is at most
an old $z$, its cyclic adjunction implies $z'X\le Y'$. Hence
an upper-action bound below $k(X,Y)$ can only be the midpoint
immediately below $k(X,Y)$; a smaller anchor would be an eligible
old target below $k$. This proves \eqref{sl:lower:old-target}.
When an argument is old, one upper-action bound is the prescribed
product, the other is above it, and $k$ is its anchor. All prescribed
products and the unit are therefore retained. The cuts in
\eqref{sl:lower:k} show $k$ is monotone in both arguments; the
anchors and old actions are monotone too, proving monotonicity of $P$.

The formula for $P$ gives the following exact ternary comparison:
\[
 P(X,Y)\le Z'\iff
 \bigl(r(X)Y\le Z'\bigr)\ \text{or}\
 \bigl(r(Y)Z\le X'\bigr)\ \text{or}\
 \bigl(r(Z)X\le Y'\bigr).
 \tag{L17}\label{sl:lower:ternary}
\]
Indeed the term $k(X,Y)$ is old, so $k(X,Y)\le Z'$ is
$k(X,Y)\le r(Z)'$, which is the third clause by
\eqref{sl:lower:k-cut}. The second upper-action bound yields the
second clause by old cyclic adjunction; the other gives the first.
Cyclic rotation permutes the clauses. Since prime is an antitone
involution, this proves both adjunctions in \eqref{sl:lower:residuals}.
In particular the residual into $d=e'$ is the prescribed prime.

For completeness the old-factor equations follow from the cut,
not from an assumption of associativity. The action laws give
$r(uX)=ur(X)$ and $r(Xu)=r(X)u$, while comparisons with any
old $z$ give
\[
 k(uX,Y)=u k(X,Y),\quad k(X,Yu)=k(X,Y)u,
 \quad k(Xu,Y)=k(X,uY).
 \tag{L18}\label{sl:lower:k-actions}
\]
For the first identity, its left side is at most $z$ exactly when
$(z'u)X\le Y'$, equivalently $k(X,Y)\le(z'u)'$, equivalently
$u k(X,Y)\le z$. For the second, $z'X\le(Yu)'$ is equivalent
to $u(z'X)\le Y'$, and hence to $k(X,Y)u\le z$. For the
third, $(z'X)u\le Y'$ is equivalent to $z'X\le(uY)'$.
A monotone map of a chain preserves finite minima. Applying these
identities to the three terms in \eqref{sl:lower:P}, and using
commutation of left and right old actions, proves
\eqref{sl:lower:old-assoc}.

Finally the comparison of $k(P(X,Y),Z)$ with old $w$ is
$w'P(X,Y)\le Z'$. By \eqref{sl:lower:old-assoc} and cyclicity
this is equivalent to $P(w'X,Y)\le Z'$, then to
$P(Y,Z)\le(w'X)'$, and finally to
$w'X\le P(Y,Z)'$. By \eqref{sl:lower:k-cut} this says
$k(X,P(Y,Z))\le w$. The two old elements are equal, and
\eqref{sl:lower:old-target} proves \eqref{sl:lower:anchor-assoc}.
\end{proof}

\subsection{Transfer of all finite histories}
\label{sl:lower:appendix-history-transfer}

\begin{proof}[Proof of Lemma~\ref{sl:lower:history-transfer}]
The unit, paired involution and actual residuals were proved in
Lemma~\ref{sl:lower:selection}. The coefficient $d$ has no active
domain by Lemma~\ref{sl:lower:collapse}, so on every new point
$X$ both $dX$ and $Xd$ equal the old element $d r(X)=r(X)d$.
Thus $d$ remains central. The residual formulas give
$X\backslash e=(dX)'=(Xd)'=e/X$, proving $e$-cyclicity.
The old products, prime, constants and order are retained; the
residual formulas then retain both old residuals, as well as the
dual product and lattice operations. This proves fullness.

We verify every original finite history. An arbitrary central floor
power $c=q^j$ has local inverse $\bar c=q^{-j}$. It satisfies
\[
                         cc^-=c^-c=p.
 \tag{L22}\label{sl:lower:floor-return}
\]
Indeed $c\bar c=p\le e$ gives $\bar c\le c^-$ and hence
$p\le cc^-\le e$. Since $cc^-\in pA$, multiplying the upper
inequality by $p$ gives $cc^-\le p$; centrality gives the other
identity. The left and right activity tests for $c$ consequently
coincide, both being $pb=b$ at an upper endpoint $b$. Their
upper values agree, so every $q^j$ remains central in $D(A)$.
In particular $pZ\le Z$ for every point.

Fix a finite history $w$ with caps $C=c_w$, $H=h_w$, and put
$j=\#S(w)-\#T(w)$. Fullness of the old projection and the
history recursions in the local group give
\[
 c=pC=q^j,\qquad h=pH=q^{-j},\qquad ch=p.
 \tag{L23}\label{sl:lower:projected-history}
\]
Suppose $XY\le C$, $UV\le H$, and write $B=YX$, $E=VU$.
Multiplying the first premise by the old central local inverse
$\bar c$ gives $\bar cXY\le p\le e$. Unit cyclicity and
centrality rotate this to $\bar cYX\le e$; multiplication by
$c$ gives $pB\le c$. Similarly $pE\le h$.

If $w$ ends in $S$, then $C=dC_0$. Centrality of $d$ and
$d^2=q^2$ give $C^2=(pC)^2=c^2$. Therefore
\[
 B^3=Y(XY)^2X\le YC^2X=c^2B=c^2(pB)\le c^3.
\]
For the monotone map $F(Z)=BZE$, associativity and centrality give
\[
 F^3(Z)=B^3ZE^3\le c^3ZE^3
       =c^3Z(pE)^3\le c^3Zh^3=pZ\le Z.
\]
If $F(Z)>Z$, monotonicity implies $F^3(Z)\ge F(Z)>Z$,
a contradiction. Hence $YXZVU\le Z$, the required history.
If $w$ ends in $T$, then $H^2=h^2$; the same proof uses
$E^3\le h^3$ and $pB\le c$, without changing the order of
factors. For the empty history, $e$-cyclicity gives $B,E\le e$
immediately. Since $w$ was arbitrary, the entire history scheme
holds in the extension.
\end{proof}

\subsection{Preservation of the insertion hypotheses}
\label{sl:lower:appendix-invariants}

\begin{proof}[Proof of Lemma~\ref{sl:lower:invariants}]
For every central idempotent $p\le e$ one has $p^-p=pp^-=p$.
Indeed $p^-\ge e$ gives $p\le pp^-\le e$. If $v=pp^-$,
then $pv=v$ by idempotence, while $v\le e$ gives $pv\le p$;
thus $v=p$. Consequently the two activity tests for $p$ at a
gap with upper endpoint $b$ are $pb=b$. Internal gaps make
them inactive, so
\[
                         pX=pr(X)=Xp\in G
 \tag{L24}\label{sl:lower:floor-new}
\]
for every new $X$. Hence $p$ stays central and $pD(A)=G$.
Fullness retains the closed support term, $q=pd$, and $d^2=q^2$.
Thus the new chain is a lower complete core with this exact group.

Every gap of $D(A)$ is one of the halves $(a,X)$ or $(X,b)$
of an old gap $(a,b)$. The upper half has $pb<b$ as before.
For the lower half, $pX$ is old and $pX\le X$, so $pX<X$.
Thus all gaps remain internal.

To prove (PC), let $J=(a,b)$ be an old gap, $X=X_J$, and
let $Y$ be the formal midpoint of either half $H$ in $D(D(A))$.
Put
\[
 E_J=\{z\in A:z'X\le X\},\qquad
 E_H=\{z\in D(A):z'Y\le Y\}.
\]
We claim $E_H\cap A=E_J$. For $z\in A$ put $u=z'$ and
$v=u^-u$, whose values are retained by fullness. If $u\le e$,
both comparisons hold. For $u\ge e$, comparison is fixation.
At the old gap this means $vb=b$ and $ub=b$. For the upper
half the criterion is identical. For the lower half it means
$vX=X$ and $uX=X$. The latter equality already implies the
former: the active inverse $u^-$ returns the fixed midpoint to
itself, so $vX=u^-(uX)=X$. This proves the claim.

Let $k=\min E_J$, which exists by (PC) and prime reversal.
The set $E_H$ is an upper set and contains $k$. If $w\in E_H$
with $w<k$, the claim forces $w$ to be new relative to $A$.
Its old upper anchor belongs to $E_H\cap A$ and is at most
$k$, hence equals $k$. There is at most one such midpoint.
Thus $E_H$ has a least member, either $k$ or that midpoint.
By prime reversal its self-comparison cut has a maximum, which
is at least $e$ and fixes $Y$. This is (PC) for each new gap.
\end{proof}

\subsection{Canonical completion and support attainment}
\label{sl:lower:appendix-canonical}

We next attain (PC) on a changed carrier. We use canonical extension
only for bounded chains, after applying the full MacNeille completion
theorem. For a bounded residuated chain $A$, let $B=A^\delta$ be
the canonical extension of its lattice reduct. Write
\[
 \mathcal K(B)=\{\bigwedge F:F\subseteq A\},\qquad
 \mathcal O(B)=\{\bigvee I:I\subseteq A\}.
\]
Every point is a join of closed points and a meet of open points.
Compactness says that $\bigwedge F\le\bigvee I$ already follows
from finite subsets. The product is its $\sigma$-extension and the
residuals their $\pi$-extensions. These form a residuated monoid,
the embedding is full, and a named cyclic dualizer remains cyclic
and dualizing. The precise standard results used here are Gehrke,
\emph{Topological Duality and Algebraic Completions}, Proposition~5,
Corollary~2, and Propositions~7--8 (including the identification of
the extended linear negation in the proof of Proposition~8)~\cite{Gehrke2022};
and Kowalski--Metcalfe, \emph{Uniform Interpolation and Coherence},
12 January 2019 manuscript, Lemmas~5.9--5.10~\cite{KowalskiMetcalfe2019}.\footnote{The respective
source copies are \href{https://www.irif.fr/~mgehrke/TopologicalDualityandAlgebraicCompletions.pdf}{Gehrke\textquotesingle s author-hosted copy}
and \href{https://boris-portal.unibe.ch/bitstreams/050b56c5-fbe5-4176-94ee-732f56081fab/download}{the Kowalski--Metcalfe institutional copy}.
The date fixes the lemma numbering of the second source.}
In particular, for closed $k,l$,
\[
 kl=\bigwedge\{ab:a,b\in A,\ k\le a,\ l\le b\},
 \tag{L25}\label{sl:lower:closed-product}
\]
and multiplication preserves arbitrary joins separately.

\begin{proof}[Proof of Lemma~\ref{sl:lower:canonical}]
First, the closed points are linearly ordered, since filters of an
old chain are linearly ordered by inclusion. If $x\not\le y$,
join density supplies a closed $k\le x$ with $k\not\le y$.
Every closed point below $y$ is then at most $k$, so $y\le k\le x$.
Thus $B$ is a chain. The quoted completion results supply
associativity, residuals and involution. The old unit can also be
checked directly: \eqref{sl:lower:closed-product} gives $ek=k=ke$
for closed $k$, by taking the old approximant $e$ in that coordinate;
join density extends this to all points.

We prove the additional defining conditions of $\KR$ directly.
If $xy\le e$ and closed $k\le x,l\le y$ are given,
\eqref{sl:lower:closed-product} and compactness give old
$a\ge k,b\ge l$ with $ab\le e$. Indeed the finite list of
products supplied by compactness can be replaced by the product
of the least first and least second coordinates. If that finite
list is empty, the target cap is top, and we may take both old
coordinates equal to top. Old $e$-cyclicity
gives $ba\le e$, hence $lk\le e$. Join density gives $yx\le e$.

Fix any finite history $w$ and suppose $xy\le c_w$, $uv\le h_w$
in $B$. The caps are old by fullness. For closed approximants
$k_x\le x,k_y\le y,k_u\le u,k_v\le v,k_z\le z$ and an
open $o\ge z$, the same compactness argument supplies old
$a\ge k_x,b\ge k_y,s\ge k_u,t\ge k_v$ with
$ab\le c_w$, $st\le h_w$. Compactness between the closed
$k_z$ and open $o$ supplies old $r$ with $k_z\le r\le o$.
The old history gives
\[
 k_yk_xk_zk_vk_u\le bart s\le r\le o.
\]
Taking joins of the five closed lower approximants and then the
meet of the open upper approximants proves $yxzvu\le z$.
This works for each finite $w$ without a length bound.

We record the compactness facts needed for attainment. If $a\prec b$
in $B$, join density makes $b$ closed and meet density makes $a$
open. If $x<y<z$, choose closed $k$ and open $o$ with
$x<k\le y\le o<z$; compactness gives an old point $t$ with
$x<t<z$. Finally, for nonempty $R\subseteq A$ and closed $b$,
\[
 (\bigwedge R)b=\bigwedge_{r\in R}rb,
 \qquad b(\bigwedge R)=\bigwedge_{r\in R}br.
 \tag{L26}\label{sl:lower:closed-meet}
\]
For the nontrivial inequality in the first identity, put
$v=\bigwedge R$. Given old $a,c$ with $v\le a$, $b\le c$,
compactness yields some $r\in R$ with $r\le a$ (if the finite
subset is empty, any $r$ works). Thus
$\bigwedge_{r\in R}rb\le rb\le rc\le ac$.
Meet over all such $a,c$ and use \eqref{sl:lower:closed-product}.
The second identity is proved with the factors in reversed order.

For a gap $a\prec b$ of $B$, its left self-support set
$S=\{u\ge e:ub=b,\ u^-b=b\}$ contains $e$ and is an initial
interval of $B_{\ge e}$. In fact if $e\le v\le u\in S$,
then $b\le vb\le ub=b$ and
$b=u^-b\le v^-b\le b$. If $S$ had no maximum, for each
$x\in S$ choose $x<y<z$ in $S$. The old-point interpolation
just proved puts some old $t$ strictly between $x$ and $z$;
initiality puts $t$ in $S$. Thus $T=S\cap A$ is nonempty and
cofinal in $S$. Put $s=\bigvee T$. Complete join preservation and
residuation give $sb=b$ and $s^-=\bigwedge_{u\in T}u^-$.
The entries $u^-$ are old and $b$ is closed, so
\eqref{sl:lower:closed-meet} gives
\[
 s^-b=\bigwedge_{u\in T}u^-b=b.
\]
Hence $s\in S$ is its maximum, a contradiction. The right-support
proof uses the other identity of \eqref{sl:lower:closed-meet}.
The maximum is idempotent as well: supports are closed under
multiplication because $v^-u^-\le(uv)^-$, while $uv\ge e$;
thus $(uv)b=b$ and
$b=v^-u^-b\le(uv)^-b\le b$. The maximum $s$ therefore has
$s^2\le s$, and $s\ge e$ gives the reverse inequality.
\end{proof}

\section{Representation above the closed support}
\label{sl:upper:section}
We retain the closed terms $p,q,Q$ of
Theorem~\ref{sl:struct:support}.  The case $p=e$ has a central
invertible dualizer and is covered by the support densification theorem.
Throughout the construction below $p>e$.

\begin{lemma}[The upper corner and its two possibilities]
\label{sl:upper:corner}
Put $I=pC$, $f=p'$, and $n=q'$.  The set $I$ is closed under
product and the original prime and has local unit $p$, local dualizer
$f$, and the actual restricted residuals.  The elements $f,n$ are
central in $C$, and
\[
 n=p\backslash e=Qf<e<p,\qquad qn=f<d<q.
\]
The corner is unit cyclic, $n<p$ is a cover in $I$, and exactly one
of the following holds: $n^2=n$, or $n$ is invertible in $I$.
Every historical guard and its prime belongs to the central set
\[
 \{e,d\}\cup\{q^j,q^jf:j\in\mathbb Z\},\qquad q^0=p.
\]
\end{lemma}
\begin{proof}
If $x=px$, then $x(px')=xx'\le d$, so $px'\le x'$ by
residuation; the reverse inequality follows from $p\ge e$.
Thus $I$ is prime closed.  Products are closed by centrality and
idempotence of $p$.  The involutive residual formulas prove both
restricted residual assertions and identify the local prime with the
original prime, with local dualizer $p'=f$.
Residual currying gives $q'=p\backslash e$; taking the complement
of the local central unit $q$ gives $n=Qf$ and $qn=f$.
Normalizing the local dualizer $f$ by the central local unit $q$
makes $n$ another cyclic dualizer in $I$.  For $z\in I$,
$z\le e$ is equivalent to $z\le n$, since $pz=z$.

We verify unit cyclicity of $I$ from an actual history, rather than
assuming it follows merely from prime closure.  Suppose $x,y\in I$
and $xy\le p$.  The $SS$ history has guards
$c_{SS}=d^2=q^2$ and $h_{SS}=(d^3)'=Q^3f$.
Apply it to $(q^2x)y\le q^2$, to $h_{SS}e\le h_{SS}$,
and to middle input $e$.  It yields
\[
 yxq^2h_{SS}=yxn\le e.
\]
This product lies in $I$, so $yxn\le n$.  Since $n$ is a
cyclic dualizer of the local monoid, its residual into itself is $p$;
therefore $yx\le p$.  Unit cyclicity makes $f$ central in $I$ by
Lemma~\ref{sl:struct:centrality}.  Any element $a$ central in $I$
is central in $C$, because $ax=a(px)=(px)a=xa$.
Thus $f$ and $n=Qf$ are globally central.

No point of $I$ lies strictly between $n$ and $p$: a point $z\le e$
satisfies $z\le p\backslash e=n$, and a point $z\ge e$ satisfies
$z=pz\ge pe=p$.  Since $e\notin I$, prime closure gives
$d\notin I$, so $n<e<p$ and $f<d<q$.

Set $s=n\backslash p=(fn)'\in I$.  We have $p\le s$,
$ns\le p$, and $n\backslash n=(qn)'=f'=p$.
The cover implies either $ns=p$ or $ns\le n$.  In the first case
centrality gives $sn=p$.  In the second, residuation gives $s\le p$,
so $s=p$; taking primes in $(qn^2)'=s=p=(qn)'$ and multiplying
by $Q$ yields $n^2=n$.  Since $n<p$, an idempotent $n$ cannot
also be invertible.

Finally, on $I$ multiplication by $d$ is multiplication by $q$,
and $d\backslash x=Qx$: if $dy\le x\in I$, then
$q(py)\le x$ and $y\le py\le Qx$, with converse immediate.
Outside these two local orbits the guard recursion only uses
$de=d$, $dd=q^2$, $d\backslash e=Q^2f$, and
$d\backslash d=e$.  Since $e+x=d\backslash x$, induction on the
history word gives exactly the displayed central guard set; prime
exchanges the two local orbits and exchanges $e,d$.
\end{proof}

\begin{lemma}[The exact protected bands]
\label{sl:upper:exact-bands}
Let $G=U(I)$.  For $x\notin I$ put $a=px$, $b=px'$.  Then
\[
 ab=ba=q,\quad b=qa^{-1},\quad (px')'=an,
 \qquad an<x<a.
\]
The endpoints $an<a$ are consecutive in $I$, and
\[
 C\setminus I=\mathop{\dot\bigcup}_{a\in G}(an,a)_C.
\]
Prime sends $(an,a)$ to $(bn,b)$ for $b=qa^{-1}$.
Moreover $I$ is directly finite.  If $n$ is invertible, $G$ is
discretely ordered, with central predecessor $n$ of $p$, and
\[
 K=\{x:px\in G\}=G\cup(C\setminus I)
\]
is a full subalgebra, called the upper-invertible complete core.
Its fibre at $g$ is $(gn,g]$.
\end{lemma}
\begin{proof}
The closure $x\mapsto px$ gives the least point of $I$ above $x$;
by prime, $(px')'$ is the greatest point of $I$ below $x$.
For $x\notin I$ the two inequalities are strict and these two
support endpoints are consecutive.  From $x<a$ we obtain
$ax'>d$ by cyclic duality.  On the other hand
$ab=p(xx')\le pd=q$.  As $ab\in I$ and $q$ is the least point
of $I$ above $d$, it follows that $ab=q$.  Cyclicity gives
$ba>d$, and the same upper bound gives $ba=q$.  Thus $a,b$ are
local units and $b=qa^{-1}$.  Complementation of a local unit gives
$(px')'=b'=aq'=an$.  Translation by a unit carries the protected
cover $n<p$ to $an<a$, proving the converse and disjointness.
Prime gives the stated reversal of bands.

If $xy=p$ in $I$, unit cyclicity gives $yx\le p$.
Were $yx<p$, the cover would give $yx\le n$, whence
$p=(xy)^2=x(yx)y\le xny=nxy=n<p$, a contradiction.
Therefore $yx=p$.  If $n$ is a unit, translate its cover and its
inverse cover to obtain the discrete group order and all protected
covers.  Products of core points have group projection, while the
prime of a top $g$ is $fg^{-1}$ and the prime of an interior point
has projection $qg^{-1}$.  Thus $K$ is product and prime closed,
contains $e,d$, and is full.  The preceding endpoint calculation
identifies its fibres.
\end{proof}

\begin{lemma}[The interior projection in the idempotent case]
\label{sl:upper:interior}
If $n^2=n$, then $x\mapsto nx$ is the full interior projection
of Proposition~\ref{sl:struct:interior}.  Its local dualizer is
$nd=qn=f$, a central unit of $nC$ with inverse $Qn$.
Hence its complete core $K=\{x:nx\in U(nC)\}$ and its
nonunit support complement have the form and gluing interface of
Lemma~\ref{sl:struct:core}.
\end{lemma}
\begin{proof}
The preceding lemmas give centrality, idempotence and $n<e$.
Also $pn=n$, so $nd=npd=nq=f$, and
$(nd)(Qn)=n^2dQ=np=n$ in both orders.  The generic interior
projection and complete-core lemmas now apply with $r=n$.
\end{proof}

A useful history check will be used repeatedly.  If a constructed
cyclic involutive chain has central $d$, and every guard prime is
central, then $ab\le c_w$ implies $ba\le c_w$: the first inequality
is equivalent to $abc_w'\le d$, which rotates cyclically and uses
centrality of $c_w'$.  Rotating both premises and using
Lemma~\ref{sl:struct:centrality} yields
$bazvu\le c_wzh_w\le z$.  This proves all histories simultaneously.

\subsection{Dense coefficient and action extensions}
\label{sl:upper:coefficient-action-tools}

An ordered residuated monoid $U$ is \emph{unit cyclic} if
$u\backslash e=e/u$ for every $u$.  Write
\[
 u^-=u\backslash e=e/u,\qquad p_u=u^-u,\qquad q_u=uu^-.
\]
A \emph{cyclic action chain} over $U$ is a chain $L$ with an
order-reversing involution $x\mapsto x^\dagger$, commuting unital
left and right actions $\lambda_u x=ux$ and $\rho_u x=xu$, monotone
in both variables, and satisfying
\[
 ux\le y^\dagger\quad\Longleftrightarrow\quad
 yu\le x^\dagger.
\]
Thus the right adjoint of $\lambda_u$ is
$\dagger\rho_u\dagger$, and conversely.

\begin{lemma}[Return identities]\label{sl:upper:return-identities}
For $u,v\in U$,
\[
 p_{uv}\le p_v,\qquad vp_{uv}\le p_uv,\qquad
 v^-p_uv\le p_{uv}.
\]
Consequently, in every ordered left action and at every point $b$,
\[
 p_{uv}b=b\quad\Longleftrightarrow\quad
 p_vb=b\ \hbox{ and }\ p_u(vb)=vb.
\]
The right-handed statements follow by reversing multiplication.
\end{lemma}
\begin{proof}
All $p_u,q_u$ are at most $e$, and $u\le u^{--}$.
Put $t=(uv)^-$.  From $tuv,uvt\le e$ one obtains
$tu\le v^-$ and $vt\le u^-$, which give the first two
inequalities.  Also $uvv^-u^-\le uu^-\le e$, whence
$v^-u^-\le(uv)^-$ and the third inequality.  If $p_{uv}b=b$,
contractivity and the first two inequalities give both displayed
fixed-point conditions.  Conversely those conditions give
$b=v^-vb=v^-p_uvb\le p_{uv}b\le b$.
\end{proof}

\begin{lemma}[Filling action gaps]\label{sl:upper:action-gap}
Every cyclic action chain $L$ over a unit-cyclic residuated monoid
has a dense cyclic action extension obtained by inserting a copy of
$\mathbb Q$ in each adjacent old pair.  Every new point has old
lower and upper bounds.  Countability is preserved.

The same assertion holds for a selected collection of adjacent pairs
provided it is closed under involution and under the transports
defined in the proof.
\end{lemma}
\begin{proof}
Write an inserted point as $(a,b,r)$, where $a<b$ is old and adjacent
and $r\in\mathbb Q$, and set
\[
 (a,b,r)^\dagger=(b^\dagger,a^\dagger,-r).
\]
Retain the old actions and prescribe
\[
 u(a,b,r)=
 \begin{cases}
  (\operatorname{pred}(ub),ub,r),&p_ub=b,\\
  ub,&p_ub\ne b,
 \end{cases}
\]
where, in the active case $p_ub=b$, $\operatorname{pred}(ub)$ denotes
the immediate old predecessor of $ub$ whose existence is verified
in the next paragraph. On the right use the identical prescription
with $bq_u=b$.

Suppose $p_ub=b$ and put $d=ub$.  For every old $y$,
\[
 y<d\quad\Longleftrightarrow\quad u^-y<b.
\]
Indeed $y<d$ and $u^-y\ge b$ would imply
$y\ge uu^-y\ge ub=d$; the reverse implication follows from
$u^-d=b$.  If $g_-$ is the right adjoint of $\lambda_{u^-}$,
adjacency gives $y<d\iff u^-y\le a\iff y\le g_-(a)$.
Thus $c=g_-(a)$ is the required predecessor of $d$.
Moreover $ua<d$, since otherwise $p_ua\ge b>a$, contradicting
contractivity.  The new value lies between the old endpoint images,
so the action is monotone in its point argument.  Finally $q_ud=d$
and $p_{u^-}d=u^{--}u^-d\ge uu^-d=d$; contractivity gives
equality.  Therefore $u^-$ transports the target gap back to the
source gap.  The opposite monoid proves these assertions on the right.

By Lemma~\ref{sl:upper:return-identities}, $uv$ transports a gap
exactly when $v$ transports it and $u$ transports its image.  If both
transport, either composite has upper endpoint $uvb$ and coordinate
$r$; otherwise both collapse to that old point.  This proves the left
action law, and the opposite monoid gives the right action law.
For cross commutation, if $u$ transports at $b$ and $H=\rho_{q_v}$,
old commutation and the local inverse give
$H(ub)=ub\iff Hb=b$: apply $u^-$ for the forward implication.
The analogous assertion holds with left and right interchanged.
Both orders of the two actions consequently transport to the same gap
with the same coordinate, or collapse to the same old point.
The unit transports every gap identically.

For coefficient monotonicity let $u\le v$.  If the $v$-action
collapses, its old value bounds the $u$-action above.  If it transports
and $ub<vb$, then $ub$ lies below the predecessor of $vb$, giving
the same conclusion.  Otherwise $ub=vb$ and $p_vb=b$; hence
$b=v^-vb\le u^-ub=p_ub\le b$.  Both actions then transport
to the same gap with the same coordinate.  The right proof is identical.

It remains to verify cyclic adjunction.  For $k\le e$ and an old
adjacent pair $a<b$, the adjunction between
$F=\lambda_k\le\mathrm{id}$ and
$G=\dagger\rho_k\dagger\ge\mathrm{id}$ gives
\[
 F(b)=b\iff F(b)>a\iff b>G(a)\iff G(a)=a.
\]
Equivalently, $kb=b\iff a^\dagger k=a^\dagger$.
If $h=\lambda_u$ transports $a<b$ to $c<d$, its old right
adjoint $g$ satisfies $g(c)=a$.  Since $q_ud=d$, the preceding
equivalence says $c^\dagger q_u=c^\dagger$.  The right
$u$-action therefore transports $(d^\dagger,c^\dagger)$ back to
$(b^\dagger,a^\dagger)$.  Conversely every such right transport
arises this way.  Conjugating by dagger gives inverse transports for
the extended $h$ and $g$.

For new $X=(a,b,r)$ and $Z=(c,d,t)$, the inequality $hX\le Z$
holds exactly when $h(b)\le c$, or when $h$ transports the source
gap to the target gap and $r\le t$.  The inequality $X\le gZ$
has the corresponding alternatives $b\le g(c)$ and the inverse
transport with $r\le t$.  Old adjunction and the preceding paragraph
identify these alternatives.  If $Z$ is old the two inequalities reduce
to $h(b)\le Z\iff b\le g(Z)$; if $X$ is old they reduce to
$h(X)\le c\iff X\le g(c)$.  Thus adjunction holds everywhere.
All gaps have been filled densely.  The selected-gap variant uses only
the indicated transports and their duals, so the same verification applies.
\end{proof}

\begin{lemma}[Completing cyclic actions]\label{sl:upper:complete-action}
The Dedekind--MacNeille completion $\overline L$ of a cyclic action
chain carries extending cyclic actions and involution.  Its actions
preserve arbitrary joins.  If $L$ is dense, so is $\overline L$.
\end{lemma}
\begin{proof}
Old points are join and meet dense in the completion.  Set
\[
 X^\dagger=\bigwedge_{x\in L,\,x\le X}x^\dagger,
 \qquad
 uX=\bigvee_{x\in L,\,x\le X}ux,
 \qquad Xu=\bigvee_{x\in L,\,x\le X}xu.
\]
In the lower-cut presentation, applying dagger twice is exactly
MacNeille closure; thus it is an antitone involution.  Join density
and old cyclicity give
\[
\begin{aligned}
 uX\le Y^\dagger
 &\iff ux\le y^\dagger\quad(x\le X,y\le Y\text{ old})\\
 &\iff yu\le x^\dagger\quad(x\le X,y\le Y\text{ old})\\
 &\iff Yu\le X^\dagger.
\end{aligned}
\]
Consequently both actions have right adjoints and preserve all joins.
Their laws, commutation, coefficient monotonicity and unit law now follow
by expanding arbitrary points as joins of old points.  Every old value
is retained since its defining join includes that point.  If $X<Y$
were adjacent in the completion, join density would force $Y$ old and
meet density would force $X$ old, contrary to old density.
\end{proof}

\begin{lemma}[Dense graded coefficient extension]
\label{sl:upper:graded-coefficients}
Let $U$ be a unit-cyclic residuated chain and let
$\mu:U\twoheadrightarrow H$ preserve multiplication, both residuals,
order and unit, where $H$ is a biordered group.  There is a full
unit-residuated extension $U\hookrightarrow V$ and an extending full
onto degree map $V\twoheadrightarrow H$ whose degree fibres are dense
and have no least points.  An old greatest point in a fibre is retained
as its greatest point, and no greatest point is introduced in a fibre
which had none.  Countable data admit a countable such extension.

If the neutral fibre of $U$ has a greatest point, every fibre does.
If it has no greatest point, no fibre has either endpoint.
\end{lemma}
\begin{proof}
Put $J=\mathbb Q\cap(0,1]$ and order $F=U\overrightarrow\times J$
lexicographically.  Write $i(u)=(u,1)$ and $\pi(u,r)=u$.  Define
\[
 L(u,v)\iff p_uv=v,\qquad R(u,v)\iff uq_v=u,
\]
\[
 (u,r)(v,s)=
 \left(uv,\min\bigl(\{1\}\cup\{r:R(u,v)\}\cup\{s:L(u,v)\}\bigr)\right).
\]
The return identities and their opposites give
\[
 \begin{split}
 L(uv,w)&\iff L(v,w)\land L(u,vw),\\
 R(u,vw)&\iff R(u,v)\land R(uv,w).
 \end{split}
\]
Also
$L(u,v)\land R(uv,w)\iff R(v,w)\land L(u,vw)$.
For the forward implication, multiply $uvq_w=uv$ on the left
by $u^-$, using $p_uv=v$; the other assertion is immediate.
For the reverse implication, multiply $p_uvw=vw$ on the right
by $w^-$, using $vq_w=v$.  Hence each bracketing of a triple
retains exactly the same subset of its three rational coordinates,
and multiplication is associative.  Its unit is $i(e)$, and $i,\pi$
are monoid maps.

For monotonicity, the rational coordinates cause no difficulty.  If
$u<v$ and $uw=vw$, the condition $R(v,w)$ is impossible because
it would imply $u\ge uq_w=uw w^-=vw w^-=v$.  If $L(v,w)$
holds, then $w=v^-vw\le u^-uw=p_uw\le w$, so $L(u,w)$ holds.
Thus the product at the smaller base has rational coordinate no larger
than the product at the larger base.  The opposite argument handles
the other variable.  Moreover
\[
 (u,r)(v,s)\le i(e)\iff uv\le e\iff vu\le e,
 \qquad \pi^{-1}(\mathord\downarrow u)=\mathord\downarrow i(u).
\]
Thus $F$ is unit cyclic.  Each projection fibre has no least point,
and $F$ is dense: between different bases decrease the upper parameter.

Write $j$ for MacNeille closure of lower subsets of $F$.  Every
translation preimage of a principal lower ideal is closed.  Indeed,
in each projection fibre its section is empty, whole, or of the form
$\{(v,s):s\le t\}$.  Every excluded point therefore has a smaller
excluded point in the same fibre; that point bounds the whole lower
set and excludes the original point from its MacNeille closure.
Preimages of arbitrary closed lower sets are closed by intersection.
If $z$ bounds $AB$, first apply this fact to each left translation
by a point of $A$, then to each right translation by a point of $j(B)$.
It follows that $z$ bounds $j(A)j(B)$, and hence
\[
 j(A)j(B)\subseteq j(\mathord\downarrow AB).
\]
Consequently closed lower sets have associative product
$I\otimes K=j(\mathord\downarrow IK)$, principal unit
$\mathord\downarrow i(e)$, and residuals
\[
 I\backslash K=\{b:Ib\subseteq K\},\qquad
 K/I=\{b:bI\subseteq K\}.
\]
The residuals are closed by the preimage result, and their defining
equivalences prove residuation.  Associativity follows by removing each
interior closure using the displayed nucleus inequality and expansivity.

Give $(u,r)$ degree $\mu(u)$, and retain precisely the cuts
\[
 V_h=\{I=j(I):\exists a,b\in F_h\ (a\in I\subseteq\downarrow b)\},
 \qquad V=\bigcup_{h\in H}V_h.
\]
Every such cut has a unique degree.  Lower and upper same-degree bounds
show that products from $V_h,V_k$ lie in $V_{hk}$.
For the residual, bound $I\in V_h$ above by $i(u_+)$ and choose
$(v,t)\in K\cap F_k$.  If $v$ is not least in $U_k$, choose
$w<v$ in $U_k$; then $i(u_+\backslash w)$ is eligible in
$I\backslash K$ and has degree $h^{-1}k$.  If $v$ is least,
put $w=u_+^-v$.  Minimality gives $u_+w=v$, and thus
$p_{u_+}w=w$ and $i(u_+)(w,t)=(v,t)$.  This gives an eligible
point of the same required degree.  Conversely, choose
$(a,r)\in I\cap F_h$ and bound $K$ above by $i(k_+)$.
Every eligible $b$ satisfies $b\le i(a\backslash k_+)$,
an upper bound of degree $h^{-1}k$.  Therefore
$I\backslash K\in V_{h^{-1}k}$.  For $K/I$, the witnesses are
$i(w/u_+)$ or $(vu_+^-,t)$, and the upper bound is $i(k_+/a)$;
its degree is $kh^{-1}$.  Thus $V$ is residuated with the stated
full degree map.  Its unit cyclicity follows by testing all raw pairs
against $i(e)$.

For old coefficients,
\[
 \{b:i(u)b\le i(w)\}=\pi^{-1}(\downarrow(u\backslash w))
 =\downarrow i(u\backslash w),
\]
and likewise on the right.  The map $u\mapsto\downarrow i(u)$
therefore retains both residuals, as well as product, order and unit.
Within each $V_h$ the principal raw cuts are order dense: if
$I\subsetneq K$, choose $b\in K\setminus I$ in that degree
and then $c<b$ of the same degree outside $I$; such $c$ exists
by closedness of $I$ and density of $F_h$.  Then
$I\subsetneq\downarrow c\subsetneq K$.
There is no least point, because each raw fibre has none.  The
greatest-point assertions follow from the same-degree upper bound in
the definition of $V_h$.  In the countable case, the subalgebra
generated by all principal raw cuts is countable, full, and dense in
each fibre, with the same endpoint properties.

Finally, if $T$ is greatest in $U_1$ and $v\in U_h$, then for
$x\in U_h$ one has $xv^-\le T$ and therefore $x\le T/v^-$.
This last point has degree $h$ and is greatest in that fibre.
Conversely a greatest $T_h$ bounds every $uv$ with $u\in U_1$,
so $T_h/v\in U_1$ bounds the neutral fibre.  If $b$ is least
in $U_h$, every left multiplication by $u\in U_1$ preserves it:
the map has a right adjoint on that fibre and hence preserves its
least point.  Thus $ub=b$ for all such $u$, and $b/b\in U_1$
is greatest.  Absence of a neutral maximum excludes both endpoints
in every fibre.
\end{proof}

\subsection{Arbitrary multiplication in upper-support bands}
\label{sl:upper:graded-middle-section}

We use the following explicit form of the upper-invertible support
structure.  The support points are $t_g$ for a biordered group $G$,
whose identity is $1$.  There are central $\delta,\nu\in G$, with
$\nu$ the immediate predecessor of $1$, such that
\begin{gather*}
 K=\coprod_{g\in G}(t_{g\nu},t_g],\qquad
 \mu(x)=g\ \Longleftrightarrow\ t_{g\nu}<x\le t_g,\\
 \mu(xy)=\mu(x)\mu(y),\qquad
 t_gx=t_{g\mu(x)},\qquad xt_g=t_{\mu(x)g},\\
 t_g'=t_{g^{-1}\delta\nu},\qquad
 \mu(x')=\mu(x)^{-1}\delta\quad(x\notin\{t_g:g\in G\}),\\
 t_\nu<e<t_1,\qquad t_{\delta\nu}<d<t_\delta,\qquad
 e'=d,\qquad d^2=t_{\delta^2}.
\end{gather*}
Here $d$ is central.  Write $P=\{t_g:g\in G\}$ and
$N=K\setminus P$.  All references below to a support-band structure
mean these identities, with the actual cyclic involutive residuated
operations of $K$.

For the application to the upper-invertible branch, take the full
complete core $K$ from Lemma~\ref{sl:upper:exact-bands} and identify the
degree group with its actual support group $G=pK=\Un(pK)$. In this
identification the data are
\[
  1_G=p,\qquad t_g=g,\qquad \mu(x)=px,\qquad
  \delta=pd,\qquad \nu=(pd)'.
\]
Thus $t_1=p$ and $t_\nu=(pd)'$; the group identity $1_G$ is the
local unit $p$, not the ambient monoid unit $e$ of $K$.
The group products and inverses are taken with respect to this local
unit. Lemmas~\ref{sl:upper:corner} and~\ref{sl:upper:exact-bands}
give the displayed support-band identities for these data.
Since $K$ is a full subalgebra of the original member of $\KR$,
all its history conditions are inherited. The greatest neutral
coefficient required in Theorem~\ref{sl:upper:graded-dense} is supplied
by the extension of Theorem~\ref{sl:upper:neutral-maximum}.
It is therefore a property of a constructed extension, not an
additional hypothesis on the input of Theorem~\ref{sl:rep:standard}.

\begin{lemma}[Central guards imply every history]
\label{sl:upper:central-guards}
Let $L$ be a cyclic involutive residuated chain with central $d$.
If a product- and prime-closed set of central elements contains $e,d$,
then $L$ satisfies unit cyclicity and every finite history law.
\end{lemma}
\begin{proof}
The identities $(dx)'=x\backslash e$ and $(xd)'=e/x$ give unit
cyclicity.  Also $e+x=(dx')'=d\backslash x$.  All guards and their
primes therefore belong to the indicated central set.  For a guard
$c$ with central $c'$, cyclicity at $d$ gives
$ab\le c\iff c'ab\le d\iff bc'a\le d\iff ba\le c$.
Induction on the history word gives $c_wzh_w\le z$: the empty
case is equality, and the two steps use
\[
 (dc_w)z(d\backslash h_w)\le c_wzh_w,\qquad
 (d\backslash c_w)z(dh_w)\le c_wzh_w,
\]
by centrality of $d$ and residuation.  Thus $ab\le c_w$ and
$uv\le h_w$ imply $bazvu\le c_wzh_w\le z$.
\end{proof}

\begin{lemma}[Intrinsic coefficient and middle pieces]
\label{sl:upper:intrinsic-pieces}
\label{sl:upper:live-coefficients}
In a support-band structure put $x^-=x\backslash e=e/x$ and
\[
 U=\{x\in N:x^-\in N\},\qquad W=U',\qquad M=N\setminus(U\cup W).
\]
Then $U$ is a full unit-cyclic unit-residuated subchain.  Its degree
map is a full onto map to a subgroup $H\le G$.  The sets $U,W$ are
disjoint and respectively initial and final in each band, and $M'=M$.
Their product types are summarized below.  An entry contains every
product with its left factor in the row class and its right factor in
the column class; it does not assert equality of these sets.
\[
\begin{array}{c|cccc}
 \cdot & U & M & W & P \\
 \hline
 U & U & M & W & P \\
 M & M & M\cup W\cup P & P & P \\
 W & W & P & P & P \\
 P & P & P & P & P
\end{array}
\]
The entries involving $P$ use the support-factor identities from the
definition of a support-band structure.  This table describes types,
not commutativity: throughout, the degree of $xy$ is the ordered product
$\mu(x)\mu(y)$.
Every product of a coefficient with a nonsupport point is nonsupport.
For $u,v\in U$ one has
$uv'=(v/u)'$ and $v'u=(u\backslash v)'$.
\end{lemma}
\begin{proof}
For $x\in N$ of degree $g$, $t_{g^{-1}\nu}$ is eligible in
$x\backslash e$, while $t_{g^{-1}}$ and all higher degrees are not.
Consequently $x^-$ either is nonsupport of degree $g^{-1}$ or equals
$t_{g^{-1}\nu}$.  If $x\in U$, then $x\le x^{--}$ and the
ineligible upper support point show that $x^-\in U$.
If $u\in U,z\in N$ and $uz$ were support, then $u^-(uz)$
would be the support top of $z$'s degree, whereas $u^-uz\le z$.
The right-handed argument is the same.  Thus $xy$ is nonsupport for
$x,y\in U$, and the eligible nonsupport point $y^-x^-$ proves
$(xy)^-$ nonsupport.  Hence $UU\subseteq U$.

For $x,y\in U$ of degrees $g,h$, the nonsupport point $x^-y$ is
eligible in $x\backslash y$, at degree $g^{-1}h$.  The support
top there is ineligible since its product with $x$ is $t_h>y$.
Thus $z=x\backslash y$ is nonsupport of that degree.
From $xzy^-\le e$, unit cyclicity gives $z(y^-x)\le e$;
the eligible nonsupport second factor proves $z\in U$.
For $y/x$ use $yx^-$ and the inequality $(y/x)(xy^-)\le e$.
The unit belongs to $U$, so this proves fullness and the full group
degree formulas, and $H=\mu[U]$ is a subgroup.

Within a band, $x\le y\in U$ implies $x^-\ge y^-$, and the
same residual bounds force $x^-\in N$.  Thus $U$ is initial and
$W$ final.  The prime--product formulas and full residual closure
give the two reflection identities and $UW,WU\subseteq W$.
If $a\in U\cap W$, then $a^-a\in W$ and $a^-a\le e$ in the
neutral band.  Finality gives $e\in W$ and hence $d\in U$,
contradicting $d^2\in P$.  Thus $U\cap W=\varnothing$ and $M'=M$.

If $xy\in U$ for $x,y\in N$, then $y(xy)^-$ is nonsupport of
inverse degree to $x$ and $x[y(xy)^-]\le e$, so $x\in U$.
Rotating at $e$ gives $y\in U$ as well.  Next, take $w=u'\in W$
and $z\in N\setminus U$.  If $wz$ were nonsupport, then
$r=(wz)'=z\backslash u$ would be nonsupport and $zr\le u$.
The nonsupport point $ru^-$ would give $z(ru^-)\le e$, forcing
$z\in U$, a contradiction.  For $zw$ use
$r=(zw)'=u/z$ and $u^-r$, together with unit cyclicity.
Finally $um,mu$ for $u\in U,m\in M$ are nonsupport and not in
$U$.  If $um\in W$, then $u^-(um)\in W$ and $u^-(um)\le m$
at the same degree, contrary to finality and $m\notin W$.
The right proof uses $(mu)u^-$.  All asserted types follow.
\end{proof}

\begin{lemma}[Completing bands and retaining coefficient maxima]
\label{sl:upper:complete-bands}
Every support-band structure has a full extension with the same support
in which each closed band $[t_{g\nu},t_g]$ is complete.
If $U_1$ has a greatest point, all coefficient fibres have greatest
points $T_h$, and the extension's coefficient fibre at $h$ is exactly
$(t_{h\nu},T_h]$.  No new coefficient degree appears.
\end{lemma}
\begin{proof}
Take the MacNeille completion of the whole chain.  Its multiplication
is the closed lower ideal generated by old products, and its prime is
the antitone extension of old prime.  Here is the needed algebraic
verification.  Preimages of old principal lower ideals under old
translations are principal, by residuation; hence preimages of closed
lower sets are closed.  The two-successive-preimage argument in
Lemma~\ref{sl:upper:graded-coefficients} proves that MacNeille closure
is a multiplicative nucleus.  It gives associative multiplication and
preserves every old product and residual.  For completed $X,Y$,
\[
 XY\le d\iff ab\le d\quad(a\le X,b\le Y\text{ old})
 \iff Y\le X',
\]
so prime is cyclic dualizing and supplies both actual residuals.
Old central elements stay central by the join definition of product.
The old set $\{e,d\}\cup\{t_{\delta^r\nu^s}:r,s\in\mathbb Z\}$
is product- and prime-closed and central, so
Lemma~\ref{sl:upper:central-guards} gives all histories in the completion.

Retain only $K_0=\bigcup_g(t_{g\nu},t_g]$.  Join density gives
an old same-band lower witness below every point of $K_0$; hence
$px=xp=t_g$ there.  For $x,y$ in bands $g,h$, old such witnesses
give $t_{gh\nu}<xy\le t_{gh}$.  For interior $x$, old interior
witnesses above and below $x$ put $x'$ in the required interior
band; at a support top use its old prime.  Thus $K_0$ is product-
and prime-closed, contains $e,d$, and is full.  Its closed bands
inherit completeness and its support is exactly the old support.

Lemma~\ref{sl:upper:graded-coefficients} gives the old maxima $T_h$.
For $t_{h\nu}<x\le T_h$, its unit residual is at least $T_h^-$,
so it is nonsupport and $x$ is a coefficient.  If $x>T_h$,
join density supplies an old same-band point $a$ with $T_h<a\le x$.
Its supported unit residual bounds $x^-$ above, excluding $x$ from
the coefficient piece.  In an old band with no coefficient, any old
same-band lower witness gives the same argument.  This proves the
claims about coefficients and their maxima.
\end{proof}

\begin{theorem}[Dense bands over a greatest neutral coefficient]
\label{sl:upper:graded-dense}
Every countable support-band structure whose intrinsic $U_1$ has a
greatest point has a full countable dense extension, satisfying every
finite history, with exactly its old support group.  Every protected
closed band becomes densely ordered.
\end{theorem}
\noindent\emph{Formalization scope.}
The companion Lean development verifies the concrete closed-support
core instance used in Theorem~\ref{sl:rep:standard}, including the
preparation of a greatest neutral coefficient. An independent adapter
from arbitrary support-band data in the formulation above to that
concrete formal interface is not included; see
Section~\ref{sl:formal:scope}.
\begin{proof}
Apply Lemma~\ref{sl:upper:complete-bands}; the old algebra used during
the construction may now be uncountable.  Denote its coefficient
subgroup by $H$ and its coefficient maxima by $T_h$.
In degree $g$ the old middle is the open interval between
\[
 \ell_g=\begin{cases}T_g,&g\in H,\\t_{g\nu},&g\notin H,\end{cases}
 \qquad
 r_g=\begin{cases}T_{g^{-1}\delta}',&g^{-1}\delta\in H,\\
 t_g,&g^{-1}\delta\notin H.
 \end{cases}
\]
These endpoints are distinct.  Form an auxiliary complete interval
$\overline M_g=[\ell_g,r_g]$, keeping separate endpoint copies
where different intervals share an actual support endpoint.
Prime pairs the intervals at $g$ and $g^{-1}\delta$ and exchanges
their two endpoint types. For $h\in H$, membership of $g$ in $H$
is equivalent to membership of either $hg$ or $gh$ in $H$.
For the upper endpoint types, centrality of $\delta$ gives
\[
 (hg)^{-1}\delta=(g^{-1}\delta)h^{-1},\qquad
 (gh)^{-1}\delta=h^{-1}(g^{-1}\delta).
\]
Membership of these degrees in $H$ is therefore equivalent to
$g^{-1}\delta\in H$. These translations use only the subgroup
property of $H$ and centrality of $\delta$.

For $u\in U_h$, retain its old interior actions, map lower endpoints
to lower endpoints and upper endpoints to upper endpoints, at degrees
$hg$ on the left and $gh$ on the right.  The upper-endpoint values
are their actual old products, since
\[
 T_a/u=T_{ah^{-1}},\qquad u\backslash T_a=T_{h^{-1}a}.
\]
Indeed all coefficients of the indicated degree are eligible, and
the full degree map determines the residual's degree.  At support
endpoints use the support product rule.  Interior actions remain in
the middle by Lemma~\ref{sl:upper:intrinsic-pieces}.  Hence clamping
only the lower endpoints respects both action laws and commutation.
Cyclic adjunction on interiors is inherited.  At matching degrees,
if either tested point is a lower endpoint, both comparisons hold;
if one is an upper endpoint and the other is not lower, both fail,
because actions preserve interiors and endpoint types.  These exhaust
the endpoint cases.  Distinct
degrees decide both comparisons in the group order.  These are thus
cyclic action chains, with the indicated graded actions.  The set of
occupied middle degrees is invariant under left and right $H$-translation:
act by a coefficient of degree $h$ and then one of degree $h^{-1}$.

Apply Lemma~\ref{sl:upper:action-gap} only to adjacent pairs within
one auxiliary interval.  This family is invariant under prime and
transport: a transported upper endpoint remains interior or upper,
never lower, so its predecessor remains in the same target interval.
Complete each enlarged interval using
Lemma~\ref{sl:upper:complete-action}, whose proof applies separately
between the corresponding domain and target grades.  Delete its two
auxiliary endpoints and call the resulting dense interval $N_g$.
The old complete interval retains all its joins and meets.  Every new
point lies in one old adjacent gap: its old floor and ceiling otherwise
coincide, by old completeness, making it old.  Its ceiling
\[
 a(X)=\min\{m\in\overline M_g:X\le m\},\qquad X\in N_g,
\]
therefore belongs to $M_g\cup\{r_g\}$.  Endpoints used as anchors
always mean their actual old values.  The transport prescription gives
\begin{equation}
\begin{gathered}
 a(uX)=u\,a(X),\qquad a(Xu)=a(X)u,\\
 X\le m\iff a(X)\le m\quad(m\in\overline M_g).
\end{gathered}
\label{sl:upper:anchor-tests}
\end{equation}
These assertions survive completion within each inserted gap.  Actions
retain $N_g$: interior upper anchors have interior images, while every
degree-one return fixes an upper boundary, so its gap is transported
and never collapsed to the deleted endpoint.  This also handles empty
old middle intervals.  The ceiling preserves joins, being left adjoint
to the inclusion of the old complete interval.

Use the raw clone $F$ and the graded closed-cut coefficient extension
$V=\bigcup_{h\in H}V_h$ from
Lemma~\ref{sl:upper:graded-coefficients}.  For $A\in V_h$ define
\begin{equation}
 \lambda_AX=\bigvee_{f\in A\cap F_h}\pi(f)X,
 \qquad
 \rho_AX=\bigvee_{f\in A\cap F_h}X\pi(f).
\label{sl:upper:cut-actions}
\end{equation}
These joins are taken in the complete target auxiliary intervals.
They are nonempty and lie between one interior action value and the
interior value for $T_h$, so belong to the target $N$ interval.
For $X\in N_g,Z\in N_j$ and $b=jg^{-1}\in H$, let
\[
 E_b(X,Z)=\{f\in F_b:\pi(f)X\le Z\}.
\]
For the later comparison with $Z=Y'$ and $Y\in N_k$, its degree
is $j=k^{-1}\delta$, so centrality of $\delta$ gives
\[
 jg^{-1}=k^{-1}\delta g^{-1}=(gk)^{-1}\delta.
\]
Adjoin all lower-degree fibres to this boundary section.  The resulting
lower cut is MacNeille closed: each old projection fibre is whole or
empty, and every excluded point has a smaller excluded point in its
own fibre.  If $E_b$ is nonempty, the cut has a witness in $F_b$
and upper bound $i(T_b)$, so is an actual element
$\kappa_b(X,Z)\in V_b$.  An empty boundary is recorded separately,
and is not replaced by a coefficient of a smaller degree.

The maps in \eqref{sl:upper:cut-actions} preserve joins.  To verify their action law, test
$\lambda_A\lambda_BX\le Z$ in the target grade.  Expanding joins
tests every raw boundary product from $(A\cap F_h)(B\cap F_k)$
against the full comparison cut just described.  Lower-degree products
are automatically eligible.  Closedness makes this equivalent to the
same test on the MacNeille closure of the product, namely
$\lambda_{A\otimes B}X\le Z$.  The right law is identical.
Old commutation and independent joins give cross commutation.  Testing
each raw coefficient and using old cyclic adjunction gives
$\lambda_AX\le Y'\iff\rho_AY\le X'$.  The unit and the old
actions are retained because the principal cut of $i(u)$ contains
$i(u)$ and all its other same-degree projections are at most $u$.
For an old middle point these joins are old: they are bounded within
its old target middle by one image and the $T_h$ image, and old joins
are retained.  On an upper anchor they are its actual old action.
Join preservation of the ceiling consequently gives
\begin{equation}
 a(\lambda_AX)=\lambda_Aa(X),\qquad
 a(\rho_AX)=\rho_Aa(X).
\label{sl:upper:cut-anchor}
\end{equation}

In every degree $g$ now put the ordered pieces
\[
 [V_g\text{ if }g\in H]\ +N_g+
 [V_{g^{-1}\delta}^{*}\text{ if }g^{-1}\delta\in H]+\{t_g\}.
\]
Prime exchanges the two coefficient copies, retains interval prime,
and has the old support carry.  Keep coefficient products, the two
actions, support products, and the reflection rules
\[
 AB^*=(B/A)^*,\qquad A^*B=(B\backslash A)^*.
\]
Products of a reverse coefficient with a middle or reverse coefficient
are the support top of the product degree.  For $X\in N_g,Y\in N_k$
put $b=(gk)^{-1}\delta$ and define
\begin{equation}
 X\circ Y=
 \begin{cases}
 a(X)a(Y),&a(X)a(Y)\in M,\\
 \kappa_b(X,Y')^*,&\text{otherwise, }b\in H, E_b(X,Y')\ne\varnothing,\\
 t_{gk},&\text{otherwise}.
 \end{cases}
\label{sl:upper:middle-product}
\end{equation}
All comparison values in this rule are attained members of $V_b$.

If either input in \eqref{sl:upper:middle-product} is an old middle point, the output is exactly
the embedded old product of its two anchors.  Indeed \eqref{sl:upper:anchor-tests} gives
$uX\le m'\iff u\,a(X)\le m'$ for old $m$, and, by prime,
$um\le Y'\iff um\le a(Y)'$.  These tests also hold for upper
anchors.  A reverse old output gives precisely the principal lift of
its old prime as comparison cut; a support output gives an empty
boundary.  The internal case is already the first line of \eqref{sl:upper:middle-product}.
In particular every old product is retained.

We verify associativity, including the nonunit coefficient cases.
For three middle factors, an internal inner product is an old middle
point.  If both inner products are internal, both bracketings are the
old three-anchor product.  If one inner product is reverse or support,
its final product is support; the corresponding old anchor product
has that same type.  If the other inner product is internal, its
bracketing therefore gives the same old support value; if not, both
bracketings are support with the same degree.

For an outer coefficient $A\in V_h$, put $z=a(X)a(Y)$ and
choose $f_0\in A\cap F_h$, which is nonempty by the construction
of $V_h$. Every $f\in A\cap F_h$ has $\pi(f)\in U_h$ and
$\pi(f)\le T_h$. Equations \eqref{sl:upper:cut-actions}
and~\eqref{sl:upper:cut-anchor}, together with monotonicity of old
multiplication, give the two-sided bounds
\[
 \begin{aligned}
  \pi(f_0)z&\le a(\lambda_AX)a(Y)\le T_hz,\\
  z\pi(f_0)&\le a(X)a(\rho_AY)\le zT_h.
 \end{aligned}
\]
For $X\in N_g,Y\in N_k$, all three terms in the first line have
degree $hgk$, and those in the second have degree $gkh$.
If $z\in M$, both bounds in each line lie in the middle piece
by Lemma~\ref{sl:upper:intrinsic-pieces}; hence the intervening
anchor product is internal. If $z\in W$, the bounds lie in
the reverse piece by $UW\subseteq W$ and $WU\subseteq W$.
Convexity of that piece in its degree therefore keeps the anchor
product in $W$. In particular the lower witness excludes a fall
into the middle, while the upper bound excludes the support top.
If $z$ is support, both bounds equal the support point of the
target degree, so the anchor product is support as well.
Thus outer coefficient actions preserve the internal, reverse,
and support types. In the internal case,
\eqref{sl:upper:cut-actions}--\eqref{sl:upper:cut-anchor} and
preservation of old bounded joins by multiplication give the
required associativity.
For a coefficient between two middle factors, the two anchor products
agree without assuming their type in advance:
\[
\begin{aligned}
 a(\rho_AX)a(Y)
 &=\bigvee_{f\in A\cap F_h}(a(X)\pi(f))a(Y)\\
 &=\bigvee_{f\in A\cap F_h}a(X)(\pi(f)a(Y))
 =a(X)a(\lambda_AY).
\end{aligned}
\]
If this common product is internal, the equality settles the case.
In all noninternal cases the required identities are
\[
 \begin{split}
 \kappa_{bh^{-1}}(\lambda_AX,Y')&=\kappa_b(X,Y')/A,\\
 \kappa_{h^{-1}b}(X,(\rho_AY)')&=A\backslash\kappa_b(X,Y'),\\
 \kappa_c(\rho_AX,Y')&=\kappa_c(X,(\lambda_AY)'),
 \qquad c=(ghk)^{-1}\delta.
 \end{split}
\]
Each follows by testing an arbitrary coefficient in its displayed
boundary degree, using the action laws and cyclic adjunction.
If a boundary is empty, an eligible product would give an eligible
coefficient in that boundary, a contradiction.  Thus the identities
also cover support outputs.  If $b\notin H$, neither translated
boundary exists; the middle placement has the same boundary on both
sides. Together with the anchor-type preservation established above,
this completes all one-coefficient cases.
Three coefficients use associativity in $V$.  Two coefficients on a
middle use the action laws.  Two coefficients
on a reverse use $(C/B)/A=C/(AB)$, $(AB)\backslash C=
B\backslash(A\backslash C)$, and
$B\backslash(C/A)=(B\backslash C)/A$, all immediate from
residuation.  Every remaining triple has support output in both
bracketings, with the same degree.  The old unit is still a unit.

For monotonicity within a middle fibre, the anchor product is monotone
and its internal-output region is initial.  In that region use old
monotonicity; outside it, comparison sets decrease with either factor,
so their reversed values increase, with an empty boundary represented
by the support top.  Crossing the region boundary increases the type
from middle to reverse or support.  At a coefficient--middle interface,
if $g\in H$, $A\in V_g$, and $X\in N_g$, an internal output satisfies
\[
 AY\le T_gY\le T_ga(Y)\le a(X)a(Y)=X\circ Y.
\]
The right interface is identical.  Noninternal outputs lie above the
middle action values.  At a middle--reverse interface multiplication
by a middle has support upper values; coefficient multiplication
preserves the order of the pieces.  The reflection and support rules
settle the other interfaces.  Different degrees are strictly compared
by the biorder.  Thus multiplication is monotone everywhere.

For the cyclic cut it suffices to consider total degree $\delta$.
An internal middle anchor product is below $d$, and hence
$X\le a(X)\le a(Y)'\le Y'$.  Otherwise \eqref{sl:upper:middle-product} gives
\[
 X\circ Y\le d\iff e_V\le\kappa_1(X,Y')\iff X\le Y'.
\]
An empty boundary fails both conditions.  For coefficient--reverse
pairs the same cut is ordinary residuation in $V$; coefficient--
coefficient and coefficient--middle pairs satisfy both comparisons,
while reverse--middle, reverse--reverse and pairs with a support factor
fail both.  At other total degrees the group order decides the cut.
Thus prime is cyclic dualizing, and associativity supplies both actual
residuals by the prime--product formulas.  All old operations, constants
and prime were retained, so the old inclusion is full.

Finally $d=e_V^*$ is central: on coefficients this follows from unit
cyclicity of $V$, and on middle, reverse and support pieces both
products are support with equal degree, since $\delta$ is central.
Every $t_{\delta^r\nu^s}$ is central by the support product rule.
The set $\{e,d\}\cup\{t_{\delta^r\nu^s}:r,s\in\mathbb Z\}$
is closed under prime and product.  Lemma~\ref{sl:upper:central-guards}
therefore proves every finite history in the constructed algebra.
Multiplication by the retained closed support $p=t_1$ fixes precisely
the old support tops, so its support group is exactly $G$.

Every $N_g$ is dense with neither endpoint.  Each coefficient fibre
is dense without a least point and retains $T_h$ as greatest; its
reverse has no greatest point.  These facts fill every interface and
make every closed protected band dense.  They also make the whole
chain dense: between points of different degrees one may interpolate
in the upper point's band, using the protected predecessor there.
Starting with the original countable source, alternately take full
subalgebra closure and choose one interpolant in each strict pair.
The union is countable and dense, full in the constructed algebra,
and has exactly $G$, since all degrees are old and all old support
points were retained.  This proves the theorem.
\end{proof}

\subsection{A greatest neutral coefficient from proper filters}

\begin{lemma}[A bounded full extension before filter completion]
\label{sl:upper:bounded-cap}
Every nontrivial cyclic involutive residuated chain $K\in\KR$
has a full bounded extension $C\in\KR$ obtained by adjoining a new
least point $0$ and a new greatest point $1$. All original finite
histories hold simultaneously in $C$. The extension has no nonzero
zero divisors, regardless of whether the old chain has an absorbing
point in its own order.
\end{lemma}
\begin{proof}
On $C=\{0\}+K+\{1\}$ retain the old product, unit and dualizer;
make $0$ absorbing and $1$ absorb every nonzero point. Retain old
prime and put $0'=1$, $1'=0$. These rules are associative: a zero
factor makes either bracketing zero; without a zero, a factor $1$
makes either bracketing $1$; otherwise old associativity applies.
They are monotone and retain the old unit. For old factors the
cyclic cut is unchanged. A zero factor and the cases involving $1$
give directly
\[
                         xy\le d\iff y\le x'.             
\]
Thus prime is cyclic dualizing, and the prime--product formulas
supply both residuals. Old products and prime are retained, so the
old inclusion is full in the original signature. The new bottom
cannot be the product of two nonzero points.

Unit cyclicity follows by the same cases at the old threshold $e$.
Fix any finite history $w$. Its guards $c_w,h_w$ are old closed
terms. If a premise product at either guard is zero, one of its
two factors is zero, and the rotated five-factor conclusion is zero.
Otherwise both premise products are at most old guard values and
are nonzero; neither can involve $1$, so all four premise factors
belong to $K$. For an old middle factor $z$, use the old history.
For $z=0$ the conclusion is zero, and for $z=1$ the conclusion is
$1\le1$. This proves every history, with no bound on its length.
\end{proof}

\begin{lemma}[A proper-filter completion retaining every finite history]
\label{sl:upper:proper-filter}
Let $C\in\KR$ be nontrivial and bounded.  It has a full
bounded embedding into a complete chain $D\in\KR$.  If $C$ has
no nonzero zero divisors, the same is true of $D$.  In particular this one
completion retains every original finite history simultaneously.
The embedding need not preserve infinite joins or meets.
\end{lemma}
\begin{proof}
Write $C^+=C\setminus\{0\}$.
Let $P$ be the set of proper lattice filters of $C$: these are exactly
its nonempty upward closed subsets not containing $0$.  They are
linearly ordered by inclusion.  Let $D$ be the chain of all upward
closed subsets of $(P,\subseteq)$, including the empty set, ordered
by inclusion.  Arbitrary joins and meets in $D$ are unions and
intersections.  Set
\[
 i(a)=\{F\in P:a\in F\},\qquad
 j(F)=\{H\in P:F\subseteq H\}.                         
\]
Then $i$ is an injective bounded lattice homomorphism, and every
$X\in D$ is $\bigcup_{F\in X}j(F)$.  Injectivity follows, for
example, because $a<b$ is separated by the proper filter
$\mathord\uparrow b$.

For $F,G\in P$ define
\[
 F*G=\mathord\uparrow\{ab:a\in F,\ b\in G\},\qquad
 X\odot Y=\bigcup_{F\in X,\ G\in Y}j(F*G).             
\]
The filter $F*G$ is allowed to be the improper filter $C$.
For this one additional filter put $j(C)=\varnothing$.
On the set $P\cup\{C\}$, the improper filter is absorbing for
$*$, because the old bottom is absorbing.  Multiplication $*$ is
associative: taking upward closure before
or after a product gives the same result, by monotonicity of the
old multiplication.  Its unit is $\mathord\uparrow e$.
Also
\[
 j(F)\odot j(G)=j(F*G),\qquad
 i(a)\odot i(b)=i(ab).                                  
\]
For the second identity when $a,b>0$, use the principal filters
$\mathord\uparrow a,\mathord\uparrow b$ for one inclusion; for the
other, $a\in F,b\in G$ implies $ab\in F*G$.  Zero factors give
empty unions.  The product definition preserves arbitrary unions in each
coordinate.  Consequently the principal-set product identities and associativity of $*$ prove
associativity on all of $D$, with unit $i(e)$.

Define
\[
 \tau F=\{a\in C:a'\notin F\},\qquad
 X^\dagger=P\setminus\tau^{-1}(X).                       
\]
The map $\tau$ is an inclusion-reversing involution of $P$: its
values are nonempty upward closed sets containing $1$ and omitting
$0$, and $\tau^2F=F$.  Thus $\dagger$ is an order-reversing
involution of $D$ and $i(a)^\dagger=i(a')$.
For principal closed sets one has the exact equivalences
\[
 \begin{aligned}
 j(F)\odot j(G)\subseteq i(d)
 &\iff d\in F*G\\
 &\iff (\exists a\in F)(\exists b\in G)\ ab\le d\\
 &\iff (\exists b\in G)\ b'\in F\\
 &\iff j(F)\subseteq j(G)^\dagger .
 \end{aligned}                                           
\]
These equivalences also hold when $F*G=C$: its principal set is
empty, and the definition of upward closure supplies a zero-product
witness.  The third equivalence uses old cyclic duality and upward
closure of $F$.  Taking unions in the two arguments gives
\[
                 X\odot Y\subseteq i(d)
                     \iff X\subseteq Y^\dagger.          
\]
Associativity and the cyclic cut give the full residuals
$X\backslash Z=(Z^\dagger\odot X)^\dagger$ and
$Z/Y=(Y\odot Z^\dagger)^\dagger$.
Thus $i$ preserves both residuals as well as prime, product,
order and all constants.

Unit cyclicity is verified separately and exactly.  On principal
closed sets,
\[
 j(F)\odot j(G)\subseteq i(e)
 \iff (\exists a\in F)(\exists b\in G)\ ab\le e
 \iff j(G)\odot j(F)\subseteq i(e),                       
\]
using old unit cyclicity.  Arbitrary unions give the same equivalence
for all $X,Y$.

We verify all histories directly, without a square-top assumption.
Fix an arbitrary finite history $w$.  Its two guard constants
$c_w,h_w$ are old terms, so their new values are $i(c_w),i(h_w)$.
Suppose $X Y\subseteq i(c_w)$ and $R S\subseteq i(h_w)$.
For arbitrary $F\in X,G\in Y,H\in R,J\in S$,
the filter-product formulas supply $a\in F,b\in G,r\in H,s\in J$ with
$ab\le c_w$ and $rs\le h_w$.
For any $K\in Z$ and $k\in K$, the old history gives
$b a k s r\le k$.  Thus
\[
             G*F*K*J*H\supseteq K,
 \qquad j(G*F*K*J*H)\subseteq j(K)\subseteq Z.          
\]
Taking the independent unions over these five filters gives
$Y X Z S R\subseteq Z$.  If one of the relevant sets is empty,
the product is empty and the same conclusion is immediate.
This proves the chosen history, and $w$ was arbitrary.
If the old chain has no nonzero zero divisors, every product
$F*G$ of proper filters is proper.  A product of two nonempty
sets then contains a nonempty $j(F*G)$ and is nonempty.  This
proves the optional zero-freeness assertion.
\end{proof}

\begin{lemma}[The successor of an old point in the filter completion]
\label{sl:upper:filter-successor}
Let $C\in\KR$ be nontrivial and bounded, and let $D$ be its
proper-filter completion from
Lemma~\ref{sl:upper:proper-filter}, with embedding
$i:C\to D$. For every $a<1$ in $C$, put
\[
 F_a=\{c\in C:a<c\},\qquad B_a=j(F_a).
\]
Then $B_a$ is the immediate successor of $i(a)$ in $D$.
\end{lemma}
\begin{proof}
Use the proper-filter notation $P,i,j$ of
Lemma~\ref{sl:upper:proper-filter}.
The set $F_a$ is nonempty, upward closed and omits $0$, so it is
a member of $P$. Every filter containing $a$ contains
$F_a$, whereas $a\notin F_a$. Hence $i(a)\subsetneq B_a$.
If $X\in D$ strictly contains $i(a)$, choose
$H\in X\setminus i(a)$. Since $a\notin H$, upward closure
of $H$ gives $H\subseteq F_a$. Upward closure of $X$ in the
filter-inclusion order now gives $F_a\in X$ and $B_a\subseteq X$.
Thus no point lies strictly between $i(a)$ and $B_a$.
\end{proof}

\begin{theorem}[Every upper-invertible complete core acquires a neutral maximum]
\label{sl:upper:neutral-maximum}
Let $K\in\KR$ be a countable upper-invertible complete group
core, with its closed support $p>e$, support group $G=pK$,
$q=pd$, and $n=q'<e$. There is a full countable extension
$K\hookrightarrow A\in\KR$ such that
\[
 pA=G,\qquad
 \text{the intrinsic neutral coefficient fibre of $A$ has a
 greatest point.}                                      
\]
No condition is imposed on the original coefficient chain, its
support subgroup, its actions, or the middle multiplication.
\end{theorem}
\begin{proof}
First apply Lemma~\ref{sl:upper:bounded-cap}, and then
the proper-filter completion of
Lemma~\ref{sl:upper:proper-filter}. This gives
full embeddings $K\subseteq C\xrightarrow{i}D\in\KR$.
The completion preserves all old operations and every finite
history, although it may break old infinite suprema. Identify old
points with their images when writing algebraic products in $D$.

We use the principal-set product identity from
Lemma~\ref{sl:upper:proper-filter}:
\[
 j(F)\odot j(H)=j(F*H),\qquad
 F*H=\mathord\uparrow\{xy:x\in F,\ y\in H\},
\]
with $j(C)=\varnothing$ for an improper product filter.

Set $B=B_n=j(F_n)$, where $F_n=(n,1]_C$.
Lemma~\ref{sl:upper:filter-successor} gives
\[
                         n<B\le e.                       
\]
Here $e\in F_n$ proves $B\le i(e)$. We claim
\[
                         pB=Bp=p.                        
\]
Indeed, for every old $c>n$ in $K$, one has $pc\ge c>n$,
while $pc$ belongs to the old support group. The protected
support cover $n<p$ therefore gives $pc\ge p$; the same holds
on the right. For the new point $c=1$ the product is $1\ge p$.
Taking $c=e$ gives equality $pe=p$. Hence
$\mathord\uparrow(pF_n)=\mathord\uparrow p$ and likewise
$\mathord\uparrow(F_np)=\mathord\uparrow p$ in $C$.
The principal-set identity of Lemma~\ref{sl:upper:proper-filter},
with $i(p)=j(\mathord\uparrow p)$, now gives $pB=Bp=p$.

All closed support terms keep their old values in $D$. In
particular the old inverse equations make $q$ and $n$ local
units in $pD$. The upper support theory and
Lemma~\ref{sl:upper:exact-bands} therefore apply in this
ambient chain: $n<p$ remains a cover in $pD$, and each old
group point $g\in G$ bounds its protected band $(gn,g]_D$.
The subset
\[
                        E=\{x\in D:px\in G\}             
\]
is a full original-signature subalgebra containing $K$, with
$pE=G$. This uses only the ambient support identities, not density
or standardness of $D$. Explicitly products stay in $E$ because
$p$ is central and idempotent. For a support point $g\in G$,
its prime is the old support point $p'g^{-1}\in G$; for a
nonsupport point over $g$, Lemma~\ref{sl:upper:exact-bands} puts its prime over
$qg^{-1}\in G$. Thus $E$ is prime closed. It contains $e,d$,
so the prime--product residual formulas give full closure. Every
old $g\in G$ belongs to $K$, proving the exact support equality.
By the preceding product equality, $B\in E$.

Use the actual residual in $D$ to define
\[
                              T=e/B.                     
\]
Since $B\le e$, the point $e$ is eligible, so $T\ge e$.
Since $pB=p>e$, the point $p$ is ineligible; hence $T<p$.
Monotonicity gives $pT=p$, so $T\in E$ and it is nonsupport.
Moreover $TB\le e$ gives $B\le T\backslash e=T^-$.
On the other hand $Tp=p>e$ makes $p$ ineligible in
$T\backslash e$, so
\[
                         n<B\le T^-<p.                  
\]
There is no support point in $(n,p)_D$. Thus $T^-$ is
nonsupport and has neutral support degree. It follows that $T$
belongs to the intrinsic neutral coefficient fibre of $E$.

Now let $x$ be any intrinsic neutral coefficient in $E$.
Its unit residual $x^-$ is nonsupport and has neutral support
degree, by Lemma~\ref{sl:upper:live-coefficients}. Equivalently, $n$ is
eligible in $x\backslash e$, because $xn=n$, while $p$ is
ineligible because $xp=p>e$; nonsupport excludes equality at
$n$. Thus $n<x^-<p$. The successor property yields
$B\le x^-$, and therefore $xB\le e$. Actual residuation
in the definition of $T$ gives $x\le T$. This proves that $T$ is the greatest
intrinsic neutral coefficient of $E$.

Finally let $A$ be the full original-signature subalgebra of $E$
generated by the countable old $K$ together with $B,T$.
It is countable, lies in $\KR$, contains every old support point,
and has $pA=G$. Its unit residuals are inherited from $E$;
in particular $T^-$ belongs to $A$ and is nonsupport. The intrinsic
neutral coefficient fibre of $A$ therefore contains $T$, and every
one of its elements is at most $T$ by the preceding argument in
$E$. This proves the assertion.
\end{proof}

\subsection{Restoring the entire upper support}
\begin{lemma}[Protected top gluing]
\label{sl:upper:protected-glue}
Let $C\in\KR$ be countable with $p>e$ and invertible $n$ in
$I=pC$. Put $G=U(I)$ and $K=\{x:px\in G\}$.
Suppose $K$ has a countable full extension $\widehat K\in\KR$
with $p\widehat K=G$ and with every closed protected band
$[gn,g]$ densely ordered.  Then $C$ has a countable full dense
extension in $\KR$, and hence a full standard embedding.
\end{lemma}
\begin{proof}
The corner $I$ has local central unit dualizer $f=qn$, is directly
finite, and has the central unit predecessor $n$ of $p$.
The relative support extension of Corollary~\ref{sl:struct:relative-unit}
gives a countable full $I\subseteq J$ with exactly the old unit group
$G$, with $q,n$ still central, and with the pairs $gn<g$ as its only
adjacent pairs.  It is directly finite and has no endpoints.
In a directly finite monoid the nonunits form a two-sided ideal:
if $ab$ is a unit with inverse $c$, then $a(bc)=p$ makes $a$ a
unit, and then $b=a^{-1}(ab)$ is a unit.

Set $B=\widehat K\mathbin{\dot\cup}(J\setminus G)$, with projection
$\Pi:B\to J$ equal to $x\mapsto px$ on $\widehat K$ and to the
identity outside.  Order first by $\Pi$, then within a fibre by its
given order, identifying $g\in G$ with the top of its fibre.  Retain
the product of $\widehat K$ and put
\[
 xy=\Pi(x)\Pi(y)\quad\text{if one factor is outside }\widehat K.
\]
Such a mixed product is outside by the ideal property.  For every
top point $a\in J$ and every $x\in B$, these rules give
\[
 ax=a\Pi(x),\qquad xa=\Pi(x)a.
\]
For old unit tops this already holds in $\widehat K$ by centrality
of $p$.  Thus a triple with any outside factor has the same projected
top product under either bracketing; triples in $\widehat K$ use
its associativity.  The old $e$ is a unit.  Multiplication by a top
is monotone by the corresponding product in $J$.  Any nontop
multiplier has unit projection, whose translations strictly preserve
inequalities between distinct bases; equal-base comparisons use
monotonicity in $\widehat K$.  This proves both monotonicity laws.

Use prime from $\widehat K$ and the local prime from $J$ elsewhere.
The latter preserves nonunits since its dualizer $f$ is a unit.
Protected covers and Lemma~\ref{sl:upper:exact-bands} make this an
antitone involution.  In particular, for $y\in J\setminus G$ and
$x\in(gn,g]_{\widehat K}$,
\[
 y\le x'\quad\Longleftrightarrow\quad y\le g'_J.
\]
A mixed product is a nonunit support point.  The cover $f<q$ contains
$d$ strictly, so such a point is at most $d$ exactly when it is at
most $f$.  Therefore local cyclicity gives
\[
 xy\le d\iff\Pi(x)\Pi(y)\le f
 \iff\Pi(y)\le\Pi(x)'_J\iff y\le x'.
\]
The final equivalence uses the preceding protected-cover test if
$x\in\widehat K$; otherwise $x'$ is a singleton nonunit, whose
comparisons are comparisons of bases.  Two core inputs use their old
cyclicity.  Associativity and this cut give both actual residuals.

The embedding of $C$ retains every product, because an old supported
factor already satisfies $ax=a(px)$ and $xa=(px)a$, and retains prime
and constants.  It is full.  The old $d$ is central on the core and
acts on outside points through its central projection $q$.  Every
history guard and prime belongs to the central set in
Lemma~\ref{sl:upper:corner}; its elements remain central on $B$
because their projections are powers of $q$ and $f$.  The central-guard
argument following Lemma~\ref{sl:upper:interior} proves all histories.
Thus $B\in\KR$ and $pB=J$.

For density, equal-base strict pairs lie in a dense protected band.
Different bases which do not form a protected cover have an intervening
point of $J$.  If the bases are $gn<g$, a lower-fibre point below its
top is separated from the upper fibre by $gn$ itself; pairs starting
at $gn$ are separated inside $[gn,g]$.  These are all cases, since
the protected covers are the only gaps in $J$.  Its lack of endpoints
also gives that property to $B$.  The standard completion theorem
now gives a full standard embedding.
\end{proof}

\begin{corollary}[The upper-invertible branch]
\label{sl:upper:invertible-core}
Every countable upper-invertible complete core has a countable full
extension with exactly its old support group and with every protected
closed band densely ordered.  Every countable whole input in that
branch has a countable full dense extension and a full standard embedding.
\end{corollary}
\begin{proof}
Theorem~\ref{sl:upper:neutral-maximum} supplies a countable full core
extension with the same support group and a greatest neutral coefficient.
Theorem~\ref{sl:upper:graded-dense} gives the required exact-group dense
extension.  Lemma~\ref{sl:upper:protected-glue} restores the entire
original upper support and fills its remaining gaps.
\end{proof}

\subsection{The idempotent upper branch}

Throughout this subsection we write $\operatorname{Un}(M)=U(M)$ for
the group of units of a local monoid $M$, keeping the notation distinct
from the coefficient chains denoted by $U$ earlier in this section.

\begin{lemma}[Positive bands with separate group floors]
\label{sl:upper:separate-bands}
Let $A\in\KR$ be countable, with $p>e$, $n=(pd)'$, $n^2=n$.
Suppose $H=\operatorname{Un}(pA)$ and
\begin{equation}
 A=\bigcup_{a\in H}[na,a],\qquad
 [na,a]\cap pA=\{na,a\}.
 \label{sl:upper:unfold-1}
\end{equation}
Then $A$ has a full countable extension $A^+$ in which every
displayed closed band is densely ordered.  Its $p$-corner is exactly
the old $pA=nH\cup H$, and its $n$-floor group is exactly $nH$.
\end{lemma}
\begin{proof}
Put $q=pd\in Z(H)$ and $B_a=(na,a]$.  The map $a\mapsto na$
is injective: if $na=n$, multiplication by the local unit $a$
sends the cover $n<p$ in $pA$ to the cover $n<a$.  If $a>p$,
the point $p$ contradicts this cover; if $a<p$, then $a>n$
(otherwise multiplying by $a^{-1}$ gives $p\le n$), and $a$
contradicts the first cover.  Thus $a=p$, proving the assertion
by taking a quotient of two units.  Consequently the closed bands
are disjoint and ordered by the biordered group $H$.

For $x\in B_a,y\in B_b$, one has $p(xy)=ab$, so
$xy\in B_{ab}$.  Indeed a floor $nab$ is fixed by $p$ and is
different from the unit $ab$.  A top acts by
$ay=ab=xb$; a floor acts by $na\,y=nab=x\,nb$.
Prime maps the open band $(na,a)$ to $(nqa^{-1},qa^{-1})$
and exchanges $a$ with $nqa^{-1}$.

Form the biordered group $\widetilde H=H\overrightarrow\times
\mathbb Z$, with central predecessor $\nu=(p,-1)$ of $(p,0)$.
At degree $(a,k)$ put a labelled copy of $B_a$.  Write its top as
$t_{(a,k)}$, and order these copies lexicographically.  For positive
old labels define
\begin{equation}
 (a,k,x)(b,l,y)=(ab,k+l,xy).
 \label{sl:upper:unfold-2}
\end{equation}
The unit is $(p,0,e)$ and the dualizer is $(q,0,d)$.  Define prime
on an interior point and on a top respectively by
\begin{equation}
 (a,k,x)'=(qa^{-1},-k,x')\quad(x<a),\qquad
 t_{(a,k)}'=t_{(qa^{-1},-k-1)}.
 \label{sl:upper:unfold-3}
\end{equation}
Associativity, the unit and monotonicity follow from the old
positive-band multiplication and the biorder.  Prime is an
antitone involution.  For the cyclic cut, first compare $ab$ with
$q$; at equality compare $k+l$ with $0$.  At equality of both,
interior factors use precisely the old test $xy\le d\iff y\le x'$.
A top factor admits no positive point at that boundary degree,
and the largest eligible point is the support top one integer
level earlier, exactly as in \eqref{sl:upper:unfold-3}.  Thus prime is cyclic dualizing,
and both residuals are the actual prime--product maxima.

The dualizer is central.  The set consisting of $e,d$ and
$t_{(q^j,k)}$, for $j,k\in\mathbb Z$, is central and closed under
product and prime.  It contains all closed historical guards and
their primes.  Cyclic rotation at a guard with central prime gives
$xy\le c\iff yx\le c$; Lemma~\ref{sl:struct:centrality} then proves every
finite history.  Hence the unfolded chain is in $\KR$.
Direct calculation gives $d^2=t_{(q^2,0)}$,
$d^3\backslash d^2=t_{(q^{-1},0)}$, so its canonical support is
$t_{(p,0)}$ and its support group is exactly $\widetilde H$.
Its $(pd)'=t_{(p,-1)}$ is a local unit.  The unfolded chain is
therefore an actual upper-invertible complete group core.

Apply Corollary~\ref{sl:upper:invertible-core}.
In the resulting countable full extension $E$, the support is still
exactly $\widetilde H$, and all protected bands are dense.  Retain
only the positive bands
\begin{equation}
 P_a=(t_{(a,-1)},t_{(a,0)}]_E\qquad(a\in H).
 \label{sl:upper:unfold-4}
\end{equation}
They are product closed with $P_aP_b\subseteq P_{ab}$, because
their support degrees have integer coordinate $0$.  The prime of
an interior point is interior in $P_{qa^{-1}}$.  These assertions
are the universal upper-support identities, and do not require
that $E$ be obtained by a prescribed equivariant construction.

Adjoin a separate floor $s_a$ immediately below each open band
$P_a$, retaining its top $t_a=t_{(a,0)}$.  For $x\in P_b$ put
\begin{equation}
 s_ax=s_as_b=s_{ab},\qquad xs_a=s_{ba},\qquad
 s_a'=t_{qa^{-1}},\qquad t_a'=s_{qa^{-1}}.
 \label{sl:upper:unfold-5}
\end{equation}
Keep products of positive points, their interior primes, and the
old $e,d$ from $E$.  A floor factor makes either bracketing a floor
of the same group degree; without a floor, associativity is that
of $E$.  Monotonicity follows from the biorder and the endpoint
rules.  At product degree $q$, a floor is below $d$, and the
endpoint exchange in \eqref{sl:upper:unfold-5} gives exactly the cyclic cut.  On
positive pairs the old cut in $E$ gives that test.  At other degrees
the biorder decides both sides.  Thus this is a cyclic involutive
residuated chain with actual residual maxima.

Centrality of $d$ is retained, and the central product--prime closed
set is now
\begin{equation}
 \{e,d\}\cup\{s_{q^j},t_{q^j}:j\in\mathbb Z\}.
 \label{sl:upper:unfold-6}
\end{equation}
Centrality of the old $q$ and its local inverse is retained in $E$;
all the tops in \eqref{sl:upper:unfold-6} are their old powers.  The new floors at
these central degrees are central by \eqref{sl:upper:unfold-5}.  The same guard
argument proves every history.  Mapping old positive points to
their level-zero copies and $na$ to $s_a$ retains product, prime,
$e,d$, and hence the full signature.  The old closed terms $p,n$
are retained.  Multiplication by $p=t_p$ fixes precisely the new
floors and old tops, while multiplication by $n=s_p$ projects each
band to its floor.  The corner and floor assertions follow.
Finally every closed band $[s_a,t_a]$ is dense, since \eqref{sl:upper:unfold-4} is
a dense protected band with its lower endpoint renamed.  All
carriers used after applying the countable core theorem are countable.
\end{proof}

\begin{lemma}[Deleting the new unit-degree fibres]
\label{sl:upper:relative-corner}
Let $D\in\KR$ be countable, with central idempotent $r<e_D$,
a full projection $x\mapsto rx$ onto a group $G$, and dualizer
$d_D=s_\delta\in G$, where $s_1=r$ and $\delta\in Z(G)$.
Write its fibres as $D_g=[s_g,t_g]$ and suppose
$D_1=\{s_1,t_1\}$, where $t_1=e_D$.
Let $H$ be the image in $G$ of the global unit group of $D$.
Suppose $D$ has a full standard extension retaining a specified central
global unit $Q$ of degree $\delta$ as a central global unit.
Then $D$ has a countable full extension $D^+$ with exactly the old
floor group, unchanged fibres over $H$, dense fibres over $G\setminus H$,
and unchanged global unit group.  The element $Q$ remains central.
\end{lemma}
\begin{proof}
Every global unit $u$ of degree $h$ induces an order isomorphism
$D_1\to D_h$ by multiplication.  Its two values are $s_h$ and $u$;
thus $D_h=\{s_h,t_h\}$ and $u=t_h$.  In particular the degree map
is injective on units and $H$ is a subgroup.  Since $Q$ has degree
$\delta$, we have $\delta\in H$.

For $g\notin H$, the two products $t_gt_{g^{-1}}$ and
$t_{g^{-1}}t_g$ belong to $D_1$.  Neither can equal $t_1$.
Indeed, if $xy=t_1$ and $yx=s_1$, then floor absorption gives
\[
 t_1=(xy)^2=x(yx)y=xs_1y=s_1,
\]
a contradiction.  Thus either equality to $t_1$ would give a two-sided
inverse for $t_g$, contrary to $g\notin H$.  Hence
\begin{equation}\label{sl:upper:corner-return}
 t_gt_{g^{-1}}=s_1=t_{g^{-1}}t_g\qquad(g\notin H).
\end{equation}
If $g,h\notin H$ and $k=gh\in H$, the old global unit
$v=t_k$ carries the greatest point of one fibre to that of its
translated fibre.  Therefore $t_hv^{-1}=t_{g^{-1}}$, and
\begin{equation}\label{sl:upper:corner-cross-return}
 t_gt_hv^{-1}=s_1,\qquad t_gt_h=s_k.
\end{equation}
These are identities between old points, so they hold in every full
extension.

Identify $D$ with its image in the specified standard extension $S$.
We first check centrality of the floor in this target.  Since
$rQ=s_\delta=d_D$, the retained inverse equation gives
$r=d_DQ^{-1}$.  The target dualizer $d_D$ is central because
$S\in\KR$, and $Q$ remains a central global unit by hypothesis.
Thus $r$ is central in $S$; its idempotence and inequality $r<e_D$
are retained as well.  The interior projection therefore applies in $S$.
Pull back the old floor group:
\[
 E=\{x\in S:rx\in G\}.
\]
The full interior projection makes this a full subalgebra.  For
completeness, the old group $G$ is full in $rS$: multiplication by an
old group element and by its old inverse are mutually inverse order
isomorphisms, so its actual residuals are group division, and its local
prime is $g\mapsto s_\delta g^{-1}$.  Thus the inverse image is indeed
closed under all original operations.  Its fibre over $g$ is the whole
old-endpoint interval $[s_g,t_g]_S$.

Now retain only
\[
 E^\circ=D_H\ \cup\!\bigcup_{g\in G\setminus H}[s_g,t_g]_S,
 \qquad D_H=\bigcup_{h\in H}\{s_h,t_h\}.
\]
We check closure rather than assuming that deleting fibres commutes
with completion.  A product with a floor is the floor at its product
degree.  The old $H$-part is closed, and a product with exactly one
factor of degree outside $H$ again has degree outside $H$.  For two
outside factors $x\in E_g$, $y\in E_h$, the only case requiring a
check is $gh\in H$.  By monotonicity and
\eqref{sl:upper:corner-cross-return},
\[
 s_{gh}\le xy\le t_gt_h=s_{gh}.
\]
Thus this product is the retained old floor.  Prime sends a fibre of
degree $g$ to the fibre of degree $\delta g^{-1}$; since
$\delta\in H$, it preserves the selected degrees, and it preserves
$D_H$.  The constants lie in $D_H$.  Product and prime closure imply
closure under both actual residuals by the involutive residual formulas.
Consequently $E^\circ$ is a full subalgebra of $S$, belongs to $\KR$,
and retains every original finite history.  Centrality of $Q$ is
inherited from $S$.

No new unit has been retained.  If an outside point $x$ were a unit,
its inverse would have degree $g^{-1}$ and be bounded by
$t_{g^{-1}}$.  Equation~\eqref{sl:upper:corner-return} would force
$xx^{-1}\le s_1<t_1$, a contradiction.  Over $H$ only old points
remain, so the old unit group is exact.

Each outside fibre is a closed interval of the dense standard order,
so each of its strict pairs has an intermediate point in that same
fibre.  Starting with $D$, choose one such point for every strict
outside-fibre pair in the current countable subalgebra, take full
subalgebra closure inside $E^\circ$, and repeat countably many times.
The union $D^+$ is countable and full.  Every strict pair in any of
its outside fibres appears at a finite stage and is split at the next.
Its floor group and its entire $H$-part remain exact, and centrality of
$Q$ is retained.  This proves all assertions.
\end{proof}

\begin{theorem}[All upper-idempotent interior complete cores]
\label{sl:upper:idempotent-core}
Every countable $C\in\KR$ with $p>e$ and
$n=(pd)'$ satisfying $n^2=n$ has a full standard embedding.
Its complete core $K=\{x:nx\in\operatorname{Un}(nC)\}$ has a
full countable extension with exactly its old $n$-floor group and
densely ordered internal floor fibres.  All original operations,
actual residuals and every finite history are retained.
\end{theorem}
\begin{proof}
Work first in the complete core $K$.  Put $G=nK$, $D=pK$ and
$f=nd=p'$.  The interior projection theorem gives floor fibres
$[g,(f g^{-1})']$.  In $D$, the neutral floor fibre is exactly
$\{n,p\}$, by the protected cover $n<p$.  Its local unit is $p$
and its dualizer is $f=nq$, with $q=pd$ a central local unit.
Its canonical floor is $n$: for its own constants,
$f^3\backslash_D f^2=q^{-1}$ and
$f^2(q^{-1})^2=n$.  The residual assertion follows from its
floor degree $f^{-1}$: every point of that fibre is eligible,
and its greatest point is the old unit $q^{-1}$; higher degrees
are ineligible.  The old central set
$\{q^j,nq^j:j\in\mathbb Z\}$ is closed under the local product
and prime and contains the local constants $p,f$.  The central-guard
argument and Lemma~\ref{sl:struct:centrality} verify every local history, so
$D\in\KR$ with these local constants.  Thus $D$ is an integral-neutral regular core
of the type in the preceding lemma, with $\delta=f$ as a degree
of $G$ and central unit $q$ lifting that degree.

Let $H\le G$ be the image of $\operatorname{Un}(D)$ under $a\mapsto na$.
The injectivity argument in Lemma~\ref{sl:upper:separate-bands}
applies.  A unit $a$ of $D$ has floor fibre exactly $[na,a]$,
because $(f(na)^{-1})'=a$.  Its intersection with $D$ consists
only of its two endpoints, by the translated protected cover.
Conversely every $x\in K\setminus D$ has $px\in\operatorname{Un}(D)$
and $nx=n(px)$, by Lemma~\ref{sl:upper:exact-bands}.
Hence the full inverse image
\begin{equation}
 A=\{x\in K:nx\in H\}
 \label{sl:upper:glue-1}
\end{equation}
is precisely the union of these active closed bands, and every
point outside $A$ belongs to $D$.  Fullness follows from the
full floor projection and the subgroup $H$, which contains the
local dualizer degree.  Therefore $K=A\cup D$ and
$A\cap D=D_H$, with exactly two points in every $H$-fibre.

Apply Lemma~\ref{sl:upper:separate-bands} to obtain $A^+$ with all its fibres
dense and exactly the same corner $D_H$.  Apply Theorem~\ref{sl:lower:representation} and
Corollary~\ref{sl:lower:central-units}, retaining the old central global
unit $q$, to $D$.  Lemma~\ref{sl:upper:relative-corner}
then gives $D^+$,
with unchanged $D_H$ and dense fibres outside $H$.  Glue them
along this actual common subalgebra of the corner:
\begin{equation}
 B=A^+\cup_{D_H}D^+.
 \label{sl:upper:glue-2}
\end{equation}
Order by $G$ and by the retained orders within its fibres.  Keep
the product on $A^+$ and on $D^+$.  For $a\in A^+$ and
$x\in D^+\setminus D_H$ use
\begin{equation}
 ax=(pa)x,
 \qquad xa=x(pa),
 \label{sl:upper:glue-3}
\end{equation}
with the right sides evaluated in $D^+$.  The map $j(a)=pa$
is a multiplicative map from $A^+$ onto $D_H$, with $j(e)=p$.
For every $k\in D_H$, its endpoint multiplication in $A^+$ gives
$ak=j(a)k$ and $ka=kj(a)$, and these are again points of $D_H$.

These facts verify associativity in every mixed case.  With one
outside factor use multiplicativity of $j$ and associativity in
$D^+$.  With two outside factors, if their product returns to
$H$, it is a point of the common $D_H$ and the displayed endpoint
identities apply; otherwise all mixed products are evaluated in
$D^+$.  With three outside factors, an intermediate return is
fixed by $j$, so both bracketings remain those of $D^+$.
The old unit is a unit, and monotonicity follows within a fibre
from monotonicity of $j$ and the two original products; different
fibres are compared by the biorder.

Keep prime from $A^+$ on $H$ and from $D^+$ elsewhere.  It is an
antitone involution, because it maps degree $g$ to $\delta g^{-1}$
and $\delta\in H$.  If the total degree is $\delta$, either both
factors lie on $H$ or both lie outside it.  The former case uses
duality in $A^+$.  In the latter the product is in
$D_\delta=\{f,q\}$, and
\begin{equation}
 k\le d\quad\Longleftrightarrow\quad k\le f
 \qquad(k\in D_\delta),
 \label{sl:upper:glue-4}
\end{equation}
because $f<d<q$.  Duality in $D^+$ at its actual dualizer $f$
therefore gives exactly $xy\le d\iff y\le x'$.
At other total degrees the group order gives this equivalence.
Both residuals are consequently the actual prime--product maxima.

The old $d$ is central: it is central on $A^+$, and on the other
fibres its action is that of the retained central corner unit $q$.
The product--prime closed set
\begin{equation}
 \{e,d\}\cup\{q^j,nq^j:j\in\mathbb Z\}
 \label{sl:upper:glue-5}
\end{equation}
is central in $B$.  On $A^+$ this is \eqref{sl:upper:unfold-6}; on outside fibres
use centrality of $q$ in $D^+$ and the central floor degrees.
It contains every old closed guard and its prime.  Cyclic rotation
at each guard and Lemma~\ref{sl:struct:centrality} prove every finite history,
so $B\in\KR$ with no additional axiom. The original inclusion
retains products on each old piece. Mixed products agree with
\eqref{sl:upper:glue-3}, since $p$ was already central and every
old point outside $A$ was $p$-fixed. For old outside points
$x,y\in K\setminus A$, the product $xy$ lies in $D=pK$; if its
$n$-degree $n(xy)$ belongs to $H$, then \eqref{sl:upper:glue-1}
also gives $xy\in A$, and hence
\[
 xy\in D\cap A=pK\cap A=D_H.
\]
Thus the two descriptions agree even when an outside product returns
to an $H$-fibre. Prime is retained on each old piece, as are $e,d$,
so the inclusion preserves the entire original signature.

Let $\iota:K\hookrightarrow B$ be the full embedding just constructed,
and write
\[
 n_K=(p_Kd_K)',\qquad n_B=(p_Bd_B)'.
\]
Fullness preserves these closed terms, so $n_B=\iota(n_K)$ and
$n_B\iota(x)=\iota(n_Kx)$ for every $x\in K$.
The assertion that the floor group is unchanged is the set equality
\begin{equation}
 \{n_Bz:z\in B\}=\iota[\{n_Kx:x\in K\}].
 \label{sl:upper:glue-floor-exact}
\end{equation}
Indeed, the floor of every point of $A^+$ is an old floor over $H$
by Lemma~\ref{sl:upper:separate-bands}, while the floor of every
point of $D^+$ is an old floor in $G$ by
Lemma~\ref{sl:upper:relative-corner}. These old floor points are
identified consistently in the gluing. Since $B=A^+\cup_{D_H}D^+$,
this proves the inclusion from left to right in
\eqref{sl:upper:glue-floor-exact}. Conversely, each point
$\iota(n_Kx)$ is $n_B\iota(x)$, which proves the reverse inclusion.
Here $n_KD=n_KK=G$, since $D=p_KK$ and $n_Kp_K=n_K$.

Thus this intermediate extension $B$ is countable, has exactly the
old $n$-floor group, and has dense fibres with their old endpoints.
The equality \eqref{sl:upper:glue-floor-exact} concerns $B$ before
the subsequent support densification. The central-unit
support densification of Theorem~\ref{sl:struct:unit-dense}, followed
by the interior gluing of Theorem~\ref{sl:struct:gluing} at $r=n$,
fills all cross-fibre gaps while retaining these dense fibres.
This gives a countable full dense extension of the complete core,
and its standard completion is a full standard embedding.
The whole-input reconstruction of
Corollary~\ref{sl:struct:core-standard-reconstruction} now restores
$C$.  Its applicability at $r=n$ follows from
Lemma~\ref{sl:upper:interior}, which also holds in the standard target
because the closed terms and the retained equations $p>e$, $n^2=n$
are unchanged.  Thus the whole input has both a countable full dense
extension and a full standard embedding.
\end{proof}

\section{Scope of the Lean verification}
\label{sl:formal:scope}

The companion Lean~4 development v1.1.0 verifies the full standard
representation theorem, the relative comparison construction,
countable-chain completeness relative to arbitrary theories,
and strong standard completeness for $\HR$.
It also verifies effective Den elimination for finite proofs
without additional assumption leaves.
The definitions preserve the noncommutative signature and all
finite histories. The standard carrier is $[0,1]$ with its usual
order, and the algebraic unit is not identified by definition
with the upper endpoint.
The transitive axiom audit permits only
\texttt{propext}, \texttt{Classical.choice}, and \texttt{Quot.sound}.

For IUL and $\mathsf{HpsUL}^{*}$, the bounded embeddings and
semantic transfers are verified, but ordinary chain completeness
for arbitrary theories remains an explicit external parameter;
it is not independently derived from the syntax of these calculi.
The pointed pseudo-uninorm transfer includes a genuine additional
constant $f$, interpreted at the image of its original value
independently of the auxiliary involutive dualizer.

The graded construction verifies the concrete closed-support core,
including its greatest neutral coefficient, used by
Theorem~\ref{sl:rep:standard}.
An adapter from arbitrary support-band data in the general
formulation of Theorem~\ref{sl:upper:graded-dense}
to this concrete interface is not included.
The idempotent-branch equality
\eqref{sl:upper:glue-floor-exact}
is verified for the intermediate gluing, before subsequent
support densification.

The exported Den elimination functions are executable.
They process ordinary rule nodes recursively and invoke verified
enumeration at Den nodes.
Decidable certificate checks, together with the proved existence
of a successful code, establish termination on derivable targets.
Ordinary Cut is retained; no complexity bound or decision procedure
for underivability is claimed.
Runtime tests of the exported functions and the transitive axiom
audit are recorded separately from static source checks.

\paragraph{Development workflow and availability.}
OpenAI ChatGPT and API-based assistance were used to draft and revise
Lean code and verification scripts. Candidate changes were checked
with the pinned Lean toolchain and project verification scripts;
compiler feedback was used in subsequent revision rounds.
AI-generated suggestions were not treated as proof certificates.
Kernel checking, transitive axiom auditing and correspondence with
the paper's mathematical statements are distinct verification tasks.
The versioned source archive and reproduction instructions are available
on Zenodo~\cite{WangLeanArchive}.

\section*{Declaration of generative AI and AI-assisted technologies
in the manuscript preparation process}
During the preparation of this work, the author used OpenAI ChatGPT
to assist with manuscript revision, organization, language editing
and LaTeX preparation. AI assistance in the formalization workflow
is described in Appendix~\ref{sl:formal:scope}.
The author takes responsibility for the final mathematical statements,
proofs, references and correspondence between the manuscript and
the formal development.


\begin{thebibliography}{99}

\bibitem{MetcalfeMontagna2007}
G. Metcalfe, F. Montagna, Substructural fuzzy logics. The Journal of Symbolic Logic 72(3) (2007) 834--864. \href{https://doi.org/10.2178/jsl/1191333844}{\nolinkurl{doi:10.2178/jsl/1191333844}}

\bibitem{CiabattoniMetcalfe2008}
A. Ciabattoni, G. Metcalfe, Density elimination. Theoretical Computer Science 403(2-3) (2008) 328--346. \href{https://doi.org/10.1016/j.tcs.2008.05.019}{\nolinkurl{doi:10.1016/j.tcs.2008.05.019}}

\bibitem{MetcalfeTsinakis2017}
G. Metcalfe, C. Tsinakis, Density revisited. Soft Computing 21(1) (2017) 175--189. \href{https://doi.org/10.1007/s00500-016-2420-7}{\nolinkurl{doi:10.1007/s00500-016-2420-7}}

\bibitem{BaldiTerui2016}
P. Baldi, K. Terui, Densification of FL chains via residuated frames. Algebra Universalis 75(2) (2016) 169--195. \href{https://doi.org/10.1007/s00012-016-0372-5}{\nolinkurl{doi:10.1007/s00012-016-0372-5}}

\bibitem{GalatosHorcik2022}
N. Galatos, R. Hor\v{c}\'ik, Densification via polynomials, languages, and frames. Journal of Pure and Applied Algebra 226(3) (2022) 106852. \href{https://doi.org/10.1016/j.jpaa.2021.106852}{\nolinkurl{doi:10.1016/j.jpaa.2021.106852}}

\bibitem{WangIUL2019}
S. Wang, A Proof of the Standard Completeness for the Involutive Uninorm Logic. Symmetry 11(4) (2019) 445. \href{https://doi.org/10.3390/sym11040445}{\nolinkurl{doi:10.3390/sym11040445}}

\bibitem{WangPseudo2019}
S. Wang, The Logic of Pseudo-Uninorms and Their Residua. Symmetry 11(3) (2019) 368. \href{https://doi.org/10.3390/sym11030368}{\nolinkurl{doi:10.3390/sym11030368}}

\bibitem{JeneiFixedPoint2023}
S. Jenei, Involutive Uninorm Logic with Fixed Point enjoys finite strong standard completeness. Archive for Mathematical Logic 62 (2023) 67--86. \href{https://doi.org/10.1007/s00153-022-00839-1}{\nolinkurl{doi:10.1007/s00153-022-00839-1}}

\bibitem{JeneiGroups2022}
S. Jenei, Group Representation for Even and Odd Involutive Commutative Residuated Chains. Studia Logica 110(4) (2022) 881--922. \href{https://doi.org/10.1007/s11225-021-09981-y}{\nolinkurl{doi:10.1007/s11225-021-09981-y}}

\bibitem{JeneiGroupsCorrection2022}
S. Jenei, Correction to: Group Representation for Even and Odd Involutive Commutative Residuated Chains. Studia Logica 110(4) (2022) 923--924. \href{https://doi.org/10.1007/s11225-022-09991-4}{\nolinkurl{doi:10.1007/s11225-022-09991-4}}

\bibitem{JeneiGroupsCorrection2023}
S. Jenei, Correction to: Group Representation for Even and Odd Involutive Commutative Residuated Chains. Studia Logica 111(5) (2023) 897--898. \href{https://doi.org/10.1007/s11225-023-10059-0}{\nolinkurl{doi:10.1007/s11225-023-10059-0}}

\bibitem{JeneiDensification2025}
S. Jenei, Densification in Classes of Involutive Commutative Residuated Lattices. Studia Logica (2025). Published online 7 July 2025. \href{https://doi.org/10.1007/s11225-025-10190-0}{\nolinkurl{doi:10.1007/s11225-025-10190-0}}

\bibitem{WangLeanArchive}
S.~Wang,
\emph{CyclicInvolutive: Lean 4 Companion to
``Strong Standard Completeness of a Hypersequent Calculus
for Cyclic Involutive Residuated Lattices''},
software and accompanying materials, Zenodo, 2026.
\href{https://doi.org/10.5281/zenodo.24242771}{\nolinkurl{doi:10.5281/zenodo.24242771}}.

\bibitem{GilFerezEtAl2017Preprint}
J. Gil-F\'erez, L. Spada, C. Tsinakis, H. Zhou, Join-completions of ordered algebras. arXiv:1708.04990v1 (2017). \url{https://arxiv.org/abs/1708.04990v1}

\bibitem{GilFerezEtAl2020}
J. Gil-F\'erez, L. Spada, C. Tsinakis, H. Zhou, Join-completions of partially ordered algebras. Annals of Pure and Applied Logic 171(10) (2020) 102842. \href{https://doi.org/10.1016/j.apal.2020.102842}{\nolinkurl{doi:10.1016/j.apal.2020.102842}}

\bibitem{WangFinite2018}
S. Wang, Logics for Finite UL and IUL-Algebras Are Substructural Fuzzy Logics. Symmetry 10(12) (2018) 755. \href{https://doi.org/10.3390/sym10120755}{\nolinkurl{doi:10.3390/sym10120755}}

\bibitem{GalatosRaftery2004}
N. Galatos, J.G. Raftery, Adding Involution to Residuated Structures. Studia Logica 77(2) (2004) 181--207. \href{https://doi.org/10.1023/B:STUD.0000037126.29193.09}{\nolinkurl{doi:10.1023/B:STUD.0000037126.29193.09}}

\bibitem{Gehrke2022}
M. Gehrke, Topological Duality and Algebraic Completions. In: Galatos, N., Terui, K. (eds.) Hiroakira Ono on Substructural Logics, Outstanding Contributions to Logic, vol. 23, pp. 143--184. Springer, Cham (2022). \href{https://doi.org/10.1007/978-3-030-76920-8_4}{\nolinkurl{doi:10.1007/978-3-030-76920-8_4}}

\bibitem{KowalskiMetcalfe2019}
T. Kowalski, G. Metcalfe, Uniform interpolation and coherence. Annals of Pure and Applied Logic 170(7) (2019) 825--841. \href{https://doi.org/10.1016/j.apal.2019.02.004}{\nolinkurl{doi:10.1016/j.apal.2019.02.004}}

\end{thebibliography}
\end{document}